\documentclass[aos,preprint]{imsart}

\RequirePackage[OT1]{fontenc}
\RequirePackage{amsthm,amsmath,natbib}
\RequirePackage[colorlinks,citecolor=blue,urlcolor=blue]{hyperref}

\usepackage{color}  
\definecolor{AliceBlue}{rgb}{0.94,0.97,1.00}
\definecolor{AntiqueWhite1}{rgb}{1.00,0.94,0.86}
\definecolor{AntiqueWhite2}{rgb}{0.93,0.87,0.80}
\definecolor{AntiqueWhite3}{rgb}{0.80,0.75,0.69}
\definecolor{AntiqueWhite4}{rgb}{0.55,0.51,0.47}
\definecolor{AntiqueWhite}{rgb}{0.98,0.92,0.84}
\definecolor{BlanchedAlmond}{rgb}{1.00,0.92,0.80}
\definecolor{BlueViolet}{rgb}{0.54,0.17,0.89}
\definecolor{CadetBlue1}{rgb}{0.60,0.96,1.00}
\definecolor{CadetBlue2}{rgb}{0.56,0.90,0.93}
\definecolor{CadetBlue3}{rgb}{0.48,0.77,0.80}
\definecolor{CadetBlue4}{rgb}{0.33,0.53,0.55}
\definecolor{CadetBlue}{rgb}{0.37,0.62,0.63}
\definecolor{CornflowerBlue}{rgb}{0.39,0.58,0.93}
\definecolor{DarkBlue}{rgb}{0.00,0.00,0.55}
\definecolor{DarkCyan}{rgb}{0.00,0.55,0.55}
\definecolor{DarkGoldenrod1}{rgb}{1.00,0.73,0.06}
\definecolor{DarkGoldenrod2}{rgb}{0.93,0.68,0.05}
\definecolor{DarkGoldenrod3}{rgb}{0.80,0.58,0.05}
\definecolor{DarkGoldenrod4}{rgb}{0.55,0.40,0.03}
\definecolor{DarkGoldenrod}{rgb}{0.72,0.53,0.04}
\definecolor{DarkGray}{rgb}{0.66,0.66,0.66}
\definecolor{DarkGreen}{rgb}{0.00,0.39,0.00}
\definecolor{DarkGrey}{rgb}{0.66,0.66,0.66}
\definecolor{DarkKhaki}{rgb}{0.74,0.72,0.42}
\definecolor{DarkMagenta}{rgb}{0.55,0.00,0.55}
\definecolor{DarkOliveGreen1}{rgb}{0.79,1.00,0.44}
\definecolor{DarkOliveGreen2}{rgb}{0.74,0.93,0.41}
\definecolor{DarkOliveGreen3}{rgb}{0.64,0.80,0.35}
\definecolor{DarkOliveGreen4}{rgb}{0.43,0.55,0.24}
\definecolor{DarkOliveGreen}{rgb}{0.33,0.42,0.18}
\definecolor{DarkOrange1}{rgb}{1.00,0.50,0.00}
\definecolor{DarkOrange2}{rgb}{0.93,0.46,0.00}
\definecolor{DarkOrange3}{rgb}{0.80,0.40,0.00}
\definecolor{DarkOrange4}{rgb}{0.55,0.27,0.00}
\definecolor{DarkOrange}{rgb}{1.00,0.55,0.00}
\definecolor{DarkOrchid1}{rgb}{0.75,0.24,1.00}
\definecolor{DarkOrchid2}{rgb}{0.70,0.23,0.93}
\definecolor{DarkOrchid3}{rgb}{0.60,0.20,0.80}
\definecolor{DarkOrchid4}{rgb}{0.41,0.13,0.55}
\definecolor{DarkOrchid}{rgb}{0.60,0.20,0.80}
\definecolor{DarkRed}{rgb}{0.55,0.00,0.00}
\definecolor{DarkSalmon}{rgb}{0.91,0.59,0.48}
\definecolor{DarkSeaGreen1}{rgb}{0.76,1.00,0.76}
\definecolor{DarkSeaGreen2}{rgb}{0.71,0.93,0.71}
\definecolor{DarkSeaGreen3}{rgb}{0.61,0.80,0.61}
\definecolor{DarkSeaGreen4}{rgb}{0.41,0.55,0.41}
\definecolor{DarkSeaGreen}{rgb}{0.56,0.74,0.56}
\definecolor{DarkSlateBlue}{rgb}{0.28,0.24,0.55}
\definecolor{DarkSlateGray1}{rgb}{0.59,1.00,1.00}
\definecolor{DarkSlateGray2}{rgb}{0.55,0.93,0.93}
\definecolor{DarkSlateGray3}{rgb}{0.47,0.80,0.80}
\definecolor{DarkSlateGray4}{rgb}{0.32,0.55,0.55}
\definecolor{DarkSlateGray}{rgb}{0.18,0.31,0.31}
\definecolor{DarkSlateGrey}{rgb}{0.18,0.31,0.31}
\definecolor{DarkTurquoise}{rgb}{0.00,0.81,0.82}
\definecolor{DarkViolet}{rgb}{0.58,0.00,0.83}
\definecolor{DeepPink1}{rgb}{1.00,0.08,0.58}
\definecolor{DeepPink2}{rgb}{0.93,0.07,0.54}
\definecolor{DeepPink3}{rgb}{0.80,0.06,0.46}
\definecolor{DeepPink4}{rgb}{0.55,0.04,0.31}
\definecolor{DeepPink}{rgb}{1.00,0.08,0.58}
\definecolor{DeepSkyBlue1}{rgb}{0.00,0.75,1.00}
\definecolor{DeepSkyBlue2}{rgb}{0.00,0.70,0.93}
\definecolor{DeepSkyBlue3}{rgb}{0.00,0.60,0.80}
\definecolor{DeepSkyBlue4}{rgb}{0.00,0.41,0.55}
\definecolor{DeepSkyBlue}{rgb}{0.00,0.75,1.00}
\definecolor{DimGray}{rgb}{0.41,0.41,0.41}
\definecolor{DimGrey}{rgb}{0.41,0.41,0.41}
\definecolor{DodgerBlue1}{rgb}{0.12,0.56,1.00}
\definecolor{DodgerBlue2}{rgb}{0.11,0.53,0.93}
\definecolor{DodgerBlue3}{rgb}{0.09,0.45,0.80}
\definecolor{DodgerBlue4}{rgb}{0.06,0.31,0.55}
\definecolor{DodgerBlue}{rgb}{0.12,0.56,1.00}
\definecolor{FloralWhite}{rgb}{1.00,0.98,0.94}
\definecolor{ForestGreen}{rgb}{0.13,0.55,0.13}
\definecolor{GhostWhite}{rgb}{0.97,0.97,1.00}
\definecolor{GreenYellow}{rgb}{0.68,1.00,0.18}
\definecolor{HotPink1}{rgb}{1.00,0.43,0.71}
\definecolor{HotPink2}{rgb}{0.93,0.42,0.65}
\definecolor{HotPink3}{rgb}{0.80,0.38,0.56}
\definecolor{HotPink4}{rgb}{0.55,0.23,0.38}
\definecolor{HotPink}{rgb}{1.00,0.41,0.71}
\definecolor{IndianRed1}{rgb}{1.00,0.42,0.42}
\definecolor{IndianRed2}{rgb}{0.93,0.39,0.39}
\definecolor{IndianRed3}{rgb}{0.80,0.33,0.33}
\definecolor{IndianRed4}{rgb}{0.55,0.23,0.23}
\definecolor{IndianRed}{rgb}{0.80,0.36,0.36}
\definecolor{LavenderBlush1}{rgb}{1.00,0.94,0.96}
\definecolor{LavenderBlush2}{rgb}{0.93,0.88,0.90}
\definecolor{LavenderBlush3}{rgb}{0.80,0.76,0.77}
\definecolor{LavenderBlush4}{rgb}{0.55,0.51,0.53}
\definecolor{LavenderBlush}{rgb}{1.00,0.94,0.96}
\definecolor{LawnGreen}{rgb}{0.49,0.99,0.00}
\definecolor{LemonChiffon1}{rgb}{1.00,0.98,0.80}
\definecolor{LemonChiffon2}{rgb}{0.93,0.91,0.75}
\definecolor{LemonChiffon3}{rgb}{0.80,0.79,0.65}
\definecolor{LemonChiffon4}{rgb}{0.55,0.54,0.44}
\definecolor{LemonChiffon}{rgb}{1.00,0.98,0.80}
\definecolor{LightBlue1}{rgb}{0.75,0.94,1.00}
\definecolor{LightBlue2}{rgb}{0.70,0.87,0.93}
\definecolor{LightBlue3}{rgb}{0.60,0.75,0.80}
\definecolor{LightBlue4}{rgb}{0.41,0.51,0.55}
\definecolor{LightBlue}{rgb}{0.68,0.85,0.90}
\definecolor{LightCoral}{rgb}{0.94,0.50,0.50}
\definecolor{LightCyan1}{rgb}{0.88,1.00,1.00}
\definecolor{LightCyan2}{rgb}{0.82,0.93,0.93}
\definecolor{LightCyan3}{rgb}{0.71,0.80,0.80}
\definecolor{LightCyan4}{rgb}{0.48,0.55,0.55}
\definecolor{LightCyan}{rgb}{0.88,1.00,1.00}
\definecolor{LightGoldenrod1}{rgb}{1.00,0.93,0.55}
\definecolor{LightGoldenrod2}{rgb}{0.93,0.86,0.51}
\definecolor{LightGoldenrod3}{rgb}{0.80,0.75,0.44}
\definecolor{LightGoldenrod4}{rgb}{0.55,0.51,0.30}
\definecolor{LightGoldenrodYellow}{rgb}{0.98,0.98,0.82}
\definecolor{LightGoldenrod}{rgb}{0.93,0.87,0.51}
\definecolor{LightGray}{rgb}{0.83,0.83,0.83}
\definecolor{LightGreen}{rgb}{0.56,0.93,0.56}
\definecolor{LightGrey}{rgb}{0.83,0.83,0.83}
\definecolor{LightPink1}{rgb}{1.00,0.68,0.73}
\definecolor{LightPink2}{rgb}{0.93,0.64,0.68}
\definecolor{LightPink3}{rgb}{0.80,0.55,0.58}
\definecolor{LightPink4}{rgb}{0.55,0.37,0.40}
\definecolor{LightPink}{rgb}{1.00,0.71,0.76}
\definecolor{LightSalmon1}{rgb}{1.00,0.63,0.48}
\definecolor{LightSalmon2}{rgb}{0.93,0.58,0.45}
\definecolor{LightSalmon3}{rgb}{0.80,0.51,0.38}
\definecolor{LightSalmon4}{rgb}{0.55,0.34,0.26}
\definecolor{LightSalmon}{rgb}{1.00,0.63,0.48}
\definecolor{LightSeaGreen}{rgb}{0.13,0.70,0.67}
\definecolor{LightSkyBlue1}{rgb}{0.69,0.89,1.00}
\definecolor{LightSkyBlue2}{rgb}{0.64,0.83,0.93}
\definecolor{LightSkyBlue3}{rgb}{0.55,0.71,0.80}
\definecolor{LightSkyBlue4}{rgb}{0.38,0.48,0.55}
\definecolor{LightSkyBlue}{rgb}{0.53,0.81,0.98}
\definecolor{LightSlateBlue}{rgb}{0.52,0.44,1.00}
\definecolor{LightSlateGray}{rgb}{0.47,0.53,0.60}
\definecolor{LightSlateGrey}{rgb}{0.47,0.53,0.60}
\definecolor{LightSteelBlue1}{rgb}{0.79,0.88,1.00}
\definecolor{LightSteelBlue2}{rgb}{0.74,0.82,0.93}
\definecolor{LightSteelBlue3}{rgb}{0.64,0.71,0.80}
\definecolor{LightSteelBlue4}{rgb}{0.43,0.48,0.55}
\definecolor{LightSteelBlue}{rgb}{0.69,0.77,0.87}
\definecolor{LightYellow1}{rgb}{1.00,1.00,0.88}
\definecolor{LightYellow2}{rgb}{0.93,0.93,0.82}
\definecolor{LightYellow3}{rgb}{0.80,0.80,0.71}
\definecolor{LightYellow4}{rgb}{0.55,0.55,0.48}
\definecolor{LightYellow}{rgb}{1.00,1.00,0.88}
\definecolor{LimeGreen}{rgb}{0.20,0.80,0.20}
\definecolor{MediumAquamarine}{rgb}{0.40,0.80,0.67}
\definecolor{MediumBlue}{rgb}{0.00,0.00,0.80}
\definecolor{MediumOrchid1}{rgb}{0.88,0.40,1.00}
\definecolor{MediumOrchid2}{rgb}{0.82,0.37,0.93}
\definecolor{MediumOrchid3}{rgb}{0.71,0.32,0.80}
\definecolor{MediumOrchid4}{rgb}{0.48,0.22,0.55}
\definecolor{MediumOrchid}{rgb}{0.73,0.33,0.83}
\definecolor{MediumPurple1}{rgb}{0.67,0.51,1.00}
\definecolor{MediumPurple2}{rgb}{0.62,0.47,0.93}
\definecolor{MediumPurple3}{rgb}{0.54,0.41,0.80}
\definecolor{MediumPurple4}{rgb}{0.36,0.28,0.55}
\definecolor{MediumPurple}{rgb}{0.58,0.44,0.86}
\definecolor{MediumSeaGreen}{rgb}{0.24,0.70,0.44}
\definecolor{MediumSlateBlue}{rgb}{0.48,0.41,0.93}
\definecolor{MediumSpringGreen}{rgb}{0.00,0.98,0.60}
\definecolor{MediumTurquoise}{rgb}{0.28,0.82,0.80}
\definecolor{MediumVioletRed}{rgb}{0.78,0.08,0.52}
\definecolor{MidnightBlue}{rgb}{0.10,0.10,0.44}
\definecolor{MintCream}{rgb}{0.96,1.00,0.98}
\definecolor{MistyRose1}{rgb}{1.00,0.89,0.88}
\definecolor{MistyRose2}{rgb}{0.93,0.84,0.82}
\definecolor{MistyRose3}{rgb}{0.80,0.72,0.71}
\definecolor{MistyRose4}{rgb}{0.55,0.49,0.48}
\definecolor{MistyRose}{rgb}{1.00,0.89,0.88}
\definecolor{NavajoWhite1}{rgb}{1.00,0.87,0.68}
\definecolor{NavajoWhite2}{rgb}{0.93,0.81,0.63}
\definecolor{NavajoWhite3}{rgb}{0.80,0.70,0.55}
\definecolor{NavajoWhite4}{rgb}{0.55,0.47,0.37}
\definecolor{NavajoWhite}{rgb}{1.00,0.87,0.68}
\definecolor{NavyBlue}{rgb}{0.00,0.00,0.50}
\definecolor{OldLace}{rgb}{0.99,0.96,0.90}
\definecolor{OliveDrab1}{rgb}{0.75,1.00,0.24}
\definecolor{OliveDrab2}{rgb}{0.70,0.93,0.23}
\definecolor{OliveDrab3}{rgb}{0.60,0.80,0.20}
\definecolor{OliveDrab4}{rgb}{0.41,0.55,0.13}
\definecolor{OliveDrab}{rgb}{0.42,0.56,0.14}
\definecolor{OrangeRed1}{rgb}{1.00,0.27,0.00}
\definecolor{OrangeRed2}{rgb}{0.93,0.25,0.00}
\definecolor{OrangeRed3}{rgb}{0.80,0.22,0.00}
\definecolor{OrangeRed4}{rgb}{0.55,0.15,0.00}
\definecolor{OrangeRed}{rgb}{1.00,0.27,0.00}
\definecolor{PaleGoldenrod}{rgb}{0.93,0.91,0.67}
\definecolor{PaleGreen1}{rgb}{0.60,1.00,0.60}
\definecolor{PaleGreen2}{rgb}{0.56,0.93,0.56}
\definecolor{PaleGreen3}{rgb}{0.49,0.80,0.49}
\definecolor{PaleGreen4}{rgb}{0.33,0.55,0.33}
\definecolor{PaleGreen}{rgb}{0.60,0.98,0.60}
\definecolor{PaleTurquoise1}{rgb}{0.73,1.00,1.00}
\definecolor{PaleTurquoise2}{rgb}{0.68,0.93,0.93}
\definecolor{PaleTurquoise3}{rgb}{0.59,0.80,0.80}
\definecolor{PaleTurquoise4}{rgb}{0.40,0.55,0.55}
\definecolor{PaleTurquoise}{rgb}{0.69,0.93,0.93}
\definecolor{PaleVioletRed1}{rgb}{1.00,0.51,0.67}
\definecolor{PaleVioletRed2}{rgb}{0.93,0.47,0.62}
\definecolor{PaleVioletRed3}{rgb}{0.80,0.41,0.54}
\definecolor{PaleVioletRed4}{rgb}{0.55,0.28,0.36}
\definecolor{PaleVioletRed}{rgb}{0.86,0.44,0.58}
\definecolor{PapayaWhip}{rgb}{1.00,0.94,0.84}
\definecolor{PeachPuff1}{rgb}{1.00,0.85,0.73}
\definecolor{PeachPuff2}{rgb}{0.93,0.80,0.68}
\definecolor{PeachPuff3}{rgb}{0.80,0.69,0.58}
\definecolor{PeachPuff4}{rgb}{0.55,0.47,0.40}
\definecolor{PeachPuff}{rgb}{1.00,0.85,0.73}
\definecolor{PowderBlue}{rgb}{0.69,0.88,0.90}
\definecolor{RosyBrown1}{rgb}{1.00,0.76,0.76}
\definecolor{RosyBrown2}{rgb}{0.93,0.71,0.71}
\definecolor{RosyBrown3}{rgb}{0.80,0.61,0.61}
\definecolor{RosyBrown4}{rgb}{0.55,0.41,0.41}
\definecolor{RosyBrown}{rgb}{0.74,0.56,0.56}
\definecolor{RoyalBlue1}{rgb}{0.28,0.46,1.00}
\definecolor{RoyalBlue2}{rgb}{0.26,0.43,0.93}
\definecolor{RoyalBlue3}{rgb}{0.23,0.37,0.80}
\definecolor{RoyalBlue4}{rgb}{0.15,0.25,0.55}
\definecolor{RoyalBlue}{rgb}{0.25,0.41,0.88}
\definecolor{SaddleBrown}{rgb}{0.55,0.27,0.07}
\definecolor{SandyBrown}{rgb}{0.96,0.64,0.38}
\definecolor{SeaGreen1}{rgb}{0.33,1.00,0.62}
\definecolor{SeaGreen2}{rgb}{0.31,0.93,0.58}
\definecolor{SeaGreen3}{rgb}{0.26,0.80,0.50}
\definecolor{SeaGreen4}{rgb}{0.18,0.55,0.34}
\definecolor{SeaGreen}{rgb}{0.18,0.55,0.34}
\definecolor{SkyBlue1}{rgb}{0.53,0.81,1.00}
\definecolor{SkyBlue2}{rgb}{0.49,0.75,0.93}
\definecolor{SkyBlue3}{rgb}{0.42,0.65,0.80}
\definecolor{SkyBlue4}{rgb}{0.29,0.44,0.55}
\definecolor{SkyBlue}{rgb}{0.53,0.81,0.92}
\definecolor{SlateBlue1}{rgb}{0.51,0.44,1.00}
\definecolor{SlateBlue2}{rgb}{0.48,0.40,0.93}
\definecolor{SlateBlue3}{rgb}{0.41,0.35,0.80}
\definecolor{SlateBlue4}{rgb}{0.28,0.24,0.55}
\definecolor{SlateBlue}{rgb}{0.42,0.35,0.80}
\definecolor{SlateGray1}{rgb}{0.78,0.89,1.00}
\definecolor{SlateGray2}{rgb}{0.73,0.83,0.93}
\definecolor{SlateGray3}{rgb}{0.62,0.71,0.80}
\definecolor{SlateGray4}{rgb}{0.42,0.48,0.55}
\definecolor{SlateGray}{rgb}{0.44,0.50,0.56}
\definecolor{SlateGrey}{rgb}{0.44,0.50,0.56}
\definecolor{SpringGreen1}{rgb}{0.00,1.00,0.50}
\definecolor{SpringGreen2}{rgb}{0.00,0.93,0.46}
\definecolor{SpringGreen3}{rgb}{0.00,0.80,0.40}
\definecolor{SpringGreen4}{rgb}{0.00,0.55,0.27}
\definecolor{SpringGreen}{rgb}{0.00,1.00,0.50}
\definecolor{SteelBlue1}{rgb}{0.39,0.72,1.00}
\definecolor{SteelBlue2}{rgb}{0.36,0.67,0.93}
\definecolor{SteelBlue3}{rgb}{0.31,0.58,0.80}
\definecolor{SteelBlue4}{rgb}{0.21,0.39,0.55}
\definecolor{SteelBlue}{rgb}{0.27,0.51,0.71}
\definecolor{VioletRed1}{rgb}{1.00,0.24,0.59}
\definecolor{VioletRed2}{rgb}{0.93,0.23,0.55}
\definecolor{VioletRed3}{rgb}{0.80,0.20,0.47}
\definecolor{VioletRed4}{rgb}{0.55,0.13,0.32}
\definecolor{VioletRed}{rgb}{0.82,0.13,0.56}
\definecolor{WhiteSmoke}{rgb}{0.96,0.96,0.96}
\definecolor{YellowGreen}{rgb}{0.60,0.80,0.20}
\definecolor{aliceblue}{rgb}{0.94,0.97,1.00}
\definecolor{antiquewhite}{rgb}{0.98,0.92,0.84}
\definecolor{aquamarine1}{rgb}{0.50,1.00,0.83}
\definecolor{aquamarine2}{rgb}{0.46,0.93,0.78}
\definecolor{aquamarine3}{rgb}{0.40,0.80,0.67}
\definecolor{aquamarine4}{rgb}{0.27,0.55,0.45}
\definecolor{aquamarine}{rgb}{0.50,1.00,0.83}
\definecolor{azure1}{rgb}{0.94,1.00,1.00}
\definecolor{azure2}{rgb}{0.88,0.93,0.93}
\definecolor{azure3}{rgb}{0.76,0.80,0.80}
\definecolor{azure4}{rgb}{0.51,0.55,0.55}
\definecolor{azure}{rgb}{0.94,1.00,1.00}
\definecolor{beige}{rgb}{0.96,0.96,0.86}
\definecolor{bisque1}{rgb}{1.00,0.89,0.77}
\definecolor{bisque2}{rgb}{0.93,0.84,0.72}
\definecolor{bisque3}{rgb}{0.80,0.72,0.62}
\definecolor{bisque4}{rgb}{0.55,0.49,0.42}
\definecolor{bisque}{rgb}{1.00,0.89,0.77}
\definecolor{black}{rgb}{0.00,0.00,0.00}
\definecolor{blanchedalmond}{rgb}{1.00,0.92,0.80}
\definecolor{blue1}{rgb}{0.00,0.00,1.00}
\definecolor{blue2}{rgb}{0.00,0.00,0.93}
\definecolor{blue3}{rgb}{0.00,0.00,0.80}
\definecolor{blue4}{rgb}{0.00,0.00,0.55}
\definecolor{blueviolet}{rgb}{0.54,0.17,0.89}
\definecolor{blue}{rgb}{0.00,0.00,1.00}
\definecolor{brown1}{rgb}{1.00,0.25,0.25}
\definecolor{brown2}{rgb}{0.93,0.23,0.23}
\definecolor{brown3}{rgb}{0.80,0.20,0.20}
\definecolor{brown4}{rgb}{0.55,0.14,0.14}
\definecolor{brown}{rgb}{0.65,0.16,0.16}
\definecolor{burlywood1}{rgb}{1.00,0.83,0.61}
\definecolor{burlywood2}{rgb}{0.93,0.77,0.57}
\definecolor{burlywood3}{rgb}{0.80,0.67,0.49}
\definecolor{burlywood4}{rgb}{0.55,0.45,0.33}
\definecolor{burlywood}{rgb}{0.87,0.72,0.53}
\definecolor{cadetblue}{rgb}{0.37,0.62,0.63}
\definecolor{chartreuse1}{rgb}{0.50,1.00,0.00}
\definecolor{chartreuse2}{rgb}{0.46,0.93,0.00}
\definecolor{chartreuse3}{rgb}{0.40,0.80,0.00}
\definecolor{chartreuse4}{rgb}{0.27,0.55,0.00}
\definecolor{chartreuse}{rgb}{0.50,1.00,0.00}
\definecolor{chocolate1}{rgb}{1.00,0.50,0.14}
\definecolor{chocolate2}{rgb}{0.93,0.46,0.13}
\definecolor{chocolate3}{rgb}{0.80,0.40,0.11}
\definecolor{chocolate4}{rgb}{0.55,0.27,0.07}
\definecolor{chocolate}{rgb}{0.82,0.41,0.12}
\definecolor{coral1}{rgb}{1.00,0.45,0.34}
\definecolor{coral2}{rgb}{0.93,0.42,0.31}
\definecolor{coral3}{rgb}{0.80,0.36,0.27}
\definecolor{coral4}{rgb}{0.55,0.24,0.18}
\definecolor{coral}{rgb}{1.00,0.50,0.31}
\definecolor{cornflowerblue}{rgb}{0.39,0.58,0.93}
\definecolor{cornsilk1}{rgb}{1.00,0.97,0.86}
\definecolor{cornsilk2}{rgb}{0.93,0.91,0.80}
\definecolor{cornsilk3}{rgb}{0.80,0.78,0.69}
\definecolor{cornsilk4}{rgb}{0.55,0.53,0.47}
\definecolor{cornsilk}{rgb}{1.00,0.97,0.86}
\definecolor{cyan1}{rgb}{0.00,1.00,1.00}
\definecolor{cyan2}{rgb}{0.00,0.93,0.93}
\definecolor{cyan3}{rgb}{0.00,0.80,0.80}
\definecolor{cyan4}{rgb}{0.00,0.55,0.55}
\definecolor{cyan}{rgb}{0.00,1.00,1.00}
\definecolor{darkblue}{rgb}{0.00,0.00,0.55}
\definecolor{darkcyan}{rgb}{0.00,0.55,0.55}
\definecolor{darkgoldenrod}{rgb}{0.72,0.53,0.04}
\definecolor{darkgray}{rgb}{0.66,0.66,0.66}
\definecolor{darkgreen}{rgb}{0.00,0.39,0.00}
\definecolor{darkgrey}{rgb}{0.66,0.66,0.66}
\definecolor{darkkhaki}{rgb}{0.74,0.72,0.42}
\definecolor{darkmagenta}{rgb}{0.55,0.00,0.55}
\definecolor{darkolive}{rgb}{0.33,0.42,0.18}
\definecolor{darkorange}{rgb}{1.00,0.55,0.00}
\definecolor{darkorchid}{rgb}{0.60,0.20,0.80}
\definecolor{darkred}{rgb}{0.55,0.00,0.00}
\definecolor{darksalmon}{rgb}{0.91,0.59,0.48}
\definecolor{darksea}{rgb}{0.56,0.74,0.56}
\definecolor{darkslate}{rgb}{0.18,0.31,0.31}
\definecolor{darkslate}{rgb}{0.18,0.31,0.31}
\definecolor{darkslate}{rgb}{0.28,0.24,0.55}
\definecolor{darkturquoise}{rgb}{0.00,0.81,0.82}
\definecolor{darkviolet}{rgb}{0.58,0.00,0.83}
\definecolor{deeppink}{rgb}{1.00,0.08,0.58}
\definecolor{deepsky}{rgb}{0.00,0.75,1.00}
\definecolor{dimgray}{rgb}{0.41,0.41,0.41}
\definecolor{dimgrey}{rgb}{0.41,0.41,0.41}
\definecolor{dodgerblue}{rgb}{0.12,0.56,1.00}
\definecolor{firebrick1}{rgb}{1.00,0.19,0.19}
\definecolor{firebrick2}{rgb}{0.93,0.17,0.17}
\definecolor{firebrick3}{rgb}{0.80,0.15,0.15}
\definecolor{firebrick4}{rgb}{0.55,0.10,0.10}
\definecolor{firebrick}{rgb}{0.70,0.13,0.13}
\definecolor{floralwhite}{rgb}{1.00,0.98,0.94}
\definecolor{forestgreen}{rgb}{0.13,0.55,0.13}
\definecolor{gainsboro}{rgb}{0.86,0.86,0.86}
\definecolor{ghostwhite}{rgb}{0.97,0.97,1.00}
\definecolor{gold1}{rgb}{1.00,0.84,0.00}
\definecolor{gold2}{rgb}{0.93,0.79,0.00}
\definecolor{gold3}{rgb}{0.80,0.68,0.00}
\definecolor{gold4}{rgb}{0.55,0.46,0.00}
\definecolor{goldenrod1}{rgb}{1.00,0.76,0.15}
\definecolor{goldenrod2}{rgb}{0.93,0.71,0.13}
\definecolor{goldenrod3}{rgb}{0.80,0.61,0.11}
\definecolor{goldenrod4}{rgb}{0.55,0.41,0.08}
\definecolor{goldenrod}{rgb}{0.85,0.65,0.13}
\definecolor{gold}{rgb}{1.00,0.84,0.00}
\definecolor{gray0}{rgb}{0.00,0.00,0.00}
\definecolor{gray100}{rgb}{1.00,1.00,1.00}
\definecolor{gray10}{rgb}{0.10,0.10,0.10}
\definecolor{gray11}{rgb}{0.11,0.11,0.11}
\definecolor{gray12}{rgb}{0.12,0.12,0.12}
\definecolor{gray13}{rgb}{0.13,0.13,0.13}
\definecolor{gray14}{rgb}{0.14,0.14,0.14}
\definecolor{gray15}{rgb}{0.15,0.15,0.15}
\definecolor{gray16}{rgb}{0.16,0.16,0.16}
\definecolor{gray17}{rgb}{0.17,0.17,0.17}
\definecolor{gray18}{rgb}{0.18,0.18,0.18}
\definecolor{gray19}{rgb}{0.19,0.19,0.19}
\definecolor{gray1}{rgb}{0.01,0.01,0.01}
\definecolor{gray20}{rgb}{0.20,0.20,0.20}
\definecolor{gray21}{rgb}{0.21,0.21,0.21}
\definecolor{gray22}{rgb}{0.22,0.22,0.22}
\definecolor{gray23}{rgb}{0.23,0.23,0.23}
\definecolor{gray24}{rgb}{0.24,0.24,0.24}
\definecolor{gray25}{rgb}{0.25,0.25,0.25}
\definecolor{gray26}{rgb}{0.26,0.26,0.26}
\definecolor{gray27}{rgb}{0.27,0.27,0.27}
\definecolor{gray28}{rgb}{0.28,0.28,0.28}
\definecolor{gray29}{rgb}{0.29,0.29,0.29}
\definecolor{gray2}{rgb}{0.02,0.02,0.02}
\definecolor{gray30}{rgb}{0.30,0.30,0.30}
\definecolor{gray31}{rgb}{0.31,0.31,0.31}
\definecolor{gray32}{rgb}{0.32,0.32,0.32}
\definecolor{gray33}{rgb}{0.33,0.33,0.33}
\definecolor{gray34}{rgb}{0.34,0.34,0.34}
\definecolor{gray35}{rgb}{0.35,0.35,0.35}
\definecolor{gray36}{rgb}{0.36,0.36,0.36}
\definecolor{gray37}{rgb}{0.37,0.37,0.37}
\definecolor{gray38}{rgb}{0.38,0.38,0.38}
\definecolor{gray39}{rgb}{0.39,0.39,0.39}
\definecolor{gray3}{rgb}{0.03,0.03,0.03}
\definecolor{gray40}{rgb}{0.40,0.40,0.40}
\definecolor{gray41}{rgb}{0.41,0.41,0.41}
\definecolor{gray42}{rgb}{0.42,0.42,0.42}
\definecolor{gray43}{rgb}{0.43,0.43,0.43}
\definecolor{gray44}{rgb}{0.44,0.44,0.44}
\definecolor{gray45}{rgb}{0.45,0.45,0.45}
\definecolor{gray46}{rgb}{0.46,0.46,0.46}
\definecolor{gray47}{rgb}{0.47,0.47,0.47}
\definecolor{gray48}{rgb}{0.48,0.48,0.48}
\definecolor{gray49}{rgb}{0.49,0.49,0.49}
\definecolor{gray4}{rgb}{0.04,0.04,0.04}
\definecolor{gray50}{rgb}{0.50,0.50,0.50}
\definecolor{gray51}{rgb}{0.51,0.51,0.51}
\definecolor{gray52}{rgb}{0.52,0.52,0.52}
\definecolor{gray53}{rgb}{0.53,0.53,0.53}
\definecolor{gray54}{rgb}{0.54,0.54,0.54}
\definecolor{gray55}{rgb}{0.55,0.55,0.55}
\definecolor{gray56}{rgb}{0.56,0.56,0.56}
\definecolor{gray57}{rgb}{0.57,0.57,0.57}
\definecolor{gray58}{rgb}{0.58,0.58,0.58}
\definecolor{gray59}{rgb}{0.59,0.59,0.59}
\definecolor{gray5}{rgb}{0.05,0.05,0.05}
\definecolor{gray60}{rgb}{0.60,0.60,0.60}
\definecolor{gray61}{rgb}{0.61,0.61,0.61}
\definecolor{gray62}{rgb}{0.62,0.62,0.62}
\definecolor{gray63}{rgb}{0.63,0.63,0.63}
\definecolor{gray64}{rgb}{0.64,0.64,0.64}
\definecolor{gray65}{rgb}{0.65,0.65,0.65}
\definecolor{gray66}{rgb}{0.66,0.66,0.66}
\definecolor{gray67}{rgb}{0.67,0.67,0.67}
\definecolor{gray68}{rgb}{0.68,0.68,0.68}
\definecolor{gray69}{rgb}{0.69,0.69,0.69}
\definecolor{gray6}{rgb}{0.06,0.06,0.06}
\definecolor{gray70}{rgb}{0.70,0.70,0.70}
\definecolor{gray71}{rgb}{0.71,0.71,0.71}
\definecolor{gray72}{rgb}{0.72,0.72,0.72}
\definecolor{gray73}{rgb}{0.73,0.73,0.73}
\definecolor{gray74}{rgb}{0.74,0.74,0.74}
\definecolor{gray75}{rgb}{0.75,0.75,0.75}
\definecolor{gray76}{rgb}{0.76,0.76,0.76}
\definecolor{gray77}{rgb}{0.77,0.77,0.77}
\definecolor{gray78}{rgb}{0.78,0.78,0.78}
\definecolor{gray79}{rgb}{0.79,0.79,0.79}
\definecolor{gray7}{rgb}{0.07,0.07,0.07}
\definecolor{gray80}{rgb}{0.80,0.80,0.80}
\definecolor{gray81}{rgb}{0.81,0.81,0.81}
\definecolor{gray82}{rgb}{0.82,0.82,0.82}
\definecolor{gray83}{rgb}{0.83,0.83,0.83}
\definecolor{gray84}{rgb}{0.84,0.84,0.84}
\definecolor{gray85}{rgb}{0.85,0.85,0.85}
\definecolor{gray86}{rgb}{0.86,0.86,0.86}
\definecolor{gray87}{rgb}{0.87,0.87,0.87}
\definecolor{gray88}{rgb}{0.88,0.88,0.88}
\definecolor{gray89}{rgb}{0.89,0.89,0.89}
\definecolor{gray8}{rgb}{0.08,0.08,0.08}
\definecolor{gray90}{rgb}{0.90,0.90,0.90}
\definecolor{gray91}{rgb}{0.91,0.91,0.91}
\definecolor{gray92}{rgb}{0.92,0.92,0.92}
\definecolor{gray93}{rgb}{0.93,0.93,0.93}
\definecolor{gray94}{rgb}{0.94,0.94,0.94}
\definecolor{gray95}{rgb}{0.95,0.95,0.95}
\definecolor{gray96}{rgb}{0.96,0.96,0.96}
\definecolor{gray97}{rgb}{0.97,0.97,0.97}
\definecolor{gray98}{rgb}{0.98,0.98,0.98}
\definecolor{gray99}{rgb}{0.99,0.99,0.99}
\definecolor{gray9}{rgb}{0.09,0.09,0.09}
\definecolor{gray}{rgb}{0.75,0.75,0.75}
\definecolor{green1}{rgb}{0.00,1.00,0.00}
\definecolor{green2}{rgb}{0.00,0.93,0.00}
\definecolor{green3}{rgb}{0.00,0.80,0.00}
\definecolor{green4}{rgb}{0.00,0.55,0.00}
\definecolor{greenyellow}{rgb}{0.68,1.00,0.18}
\definecolor{green}{rgb}{0.00,1.00,0.00}
\definecolor{grey0}{rgb}{0.00,0.00,0.00}
\definecolor{grey100}{rgb}{1.00,1.00,1.00}
\definecolor{grey10}{rgb}{0.10,0.10,0.10}
\definecolor{grey11}{rgb}{0.11,0.11,0.11}
\definecolor{grey12}{rgb}{0.12,0.12,0.12}
\definecolor{grey13}{rgb}{0.13,0.13,0.13}
\definecolor{grey14}{rgb}{0.14,0.14,0.14}
\definecolor{grey15}{rgb}{0.15,0.15,0.15}
\definecolor{grey16}{rgb}{0.16,0.16,0.16}
\definecolor{grey17}{rgb}{0.17,0.17,0.17}
\definecolor{grey18}{rgb}{0.18,0.18,0.18}
\definecolor{grey19}{rgb}{0.19,0.19,0.19}
\definecolor{grey1}{rgb}{0.01,0.01,0.01}
\definecolor{grey20}{rgb}{0.20,0.20,0.20}
\definecolor{grey21}{rgb}{0.21,0.21,0.21}
\definecolor{grey22}{rgb}{0.22,0.22,0.22}
\definecolor{grey23}{rgb}{0.23,0.23,0.23}
\definecolor{grey24}{rgb}{0.24,0.24,0.24}
\definecolor{grey25}{rgb}{0.25,0.25,0.25}
\definecolor{grey26}{rgb}{0.26,0.26,0.26}
\definecolor{grey27}{rgb}{0.27,0.27,0.27}
\definecolor{grey28}{rgb}{0.28,0.28,0.28}
\definecolor{grey29}{rgb}{0.29,0.29,0.29}
\definecolor{grey2}{rgb}{0.02,0.02,0.02}
\definecolor{grey30}{rgb}{0.30,0.30,0.30}
\definecolor{grey31}{rgb}{0.31,0.31,0.31}
\definecolor{grey32}{rgb}{0.32,0.32,0.32}
\definecolor{grey33}{rgb}{0.33,0.33,0.33}
\definecolor{grey34}{rgb}{0.34,0.34,0.34}
\definecolor{grey35}{rgb}{0.35,0.35,0.35}
\definecolor{grey36}{rgb}{0.36,0.36,0.36}
\definecolor{grey37}{rgb}{0.37,0.37,0.37}
\definecolor{grey38}{rgb}{0.38,0.38,0.38}
\definecolor{grey39}{rgb}{0.39,0.39,0.39}
\definecolor{grey3}{rgb}{0.03,0.03,0.03}
\definecolor{grey40}{rgb}{0.40,0.40,0.40}
\definecolor{grey41}{rgb}{0.41,0.41,0.41}
\definecolor{grey42}{rgb}{0.42,0.42,0.42}
\definecolor{grey43}{rgb}{0.43,0.43,0.43}
\definecolor{grey44}{rgb}{0.44,0.44,0.44}
\definecolor{grey45}{rgb}{0.45,0.45,0.45}
\definecolor{grey46}{rgb}{0.46,0.46,0.46}
\definecolor{grey47}{rgb}{0.47,0.47,0.47}
\definecolor{grey48}{rgb}{0.48,0.48,0.48}
\definecolor{grey49}{rgb}{0.49,0.49,0.49}
\definecolor{grey4}{rgb}{0.04,0.04,0.04}
\definecolor{grey50}{rgb}{0.50,0.50,0.50}
\definecolor{grey51}{rgb}{0.51,0.51,0.51}
\definecolor{grey52}{rgb}{0.52,0.52,0.52}
\definecolor{grey53}{rgb}{0.53,0.53,0.53}
\definecolor{grey54}{rgb}{0.54,0.54,0.54}
\definecolor{grey55}{rgb}{0.55,0.55,0.55}
\definecolor{grey56}{rgb}{0.56,0.56,0.56}
\definecolor{grey57}{rgb}{0.57,0.57,0.57}
\definecolor{grey58}{rgb}{0.58,0.58,0.58}
\definecolor{grey59}{rgb}{0.59,0.59,0.59}
\definecolor{grey5}{rgb}{0.05,0.05,0.05}
\definecolor{grey60}{rgb}{0.60,0.60,0.60}
\definecolor{grey61}{rgb}{0.61,0.61,0.61}
\definecolor{grey62}{rgb}{0.62,0.62,0.62}
\definecolor{grey63}{rgb}{0.63,0.63,0.63}
\definecolor{grey64}{rgb}{0.64,0.64,0.64}
\definecolor{grey65}{rgb}{0.65,0.65,0.65}
\definecolor{grey66}{rgb}{0.66,0.66,0.66}
\definecolor{grey67}{rgb}{0.67,0.67,0.67}
\definecolor{grey68}{rgb}{0.68,0.68,0.68}
\definecolor{grey69}{rgb}{0.69,0.69,0.69}
\definecolor{grey6}{rgb}{0.06,0.06,0.06}
\definecolor{grey70}{rgb}{0.70,0.70,0.70}
\definecolor{grey71}{rgb}{0.71,0.71,0.71}
\definecolor{grey72}{rgb}{0.72,0.72,0.72}
\definecolor{grey73}{rgb}{0.73,0.73,0.73}
\definecolor{grey74}{rgb}{0.74,0.74,0.74}
\definecolor{grey75}{rgb}{0.75,0.75,0.75}
\definecolor{grey76}{rgb}{0.76,0.76,0.76}
\definecolor{grey77}{rgb}{0.77,0.77,0.77}
\definecolor{grey78}{rgb}{0.78,0.78,0.78}
\definecolor{grey79}{rgb}{0.79,0.79,0.79}
\definecolor{grey7}{rgb}{0.07,0.07,0.07}
\definecolor{grey80}{rgb}{0.80,0.80,0.80}
\definecolor{grey81}{rgb}{0.81,0.81,0.81}
\definecolor{grey82}{rgb}{0.82,0.82,0.82}
\definecolor{grey83}{rgb}{0.83,0.83,0.83}
\definecolor{grey84}{rgb}{0.84,0.84,0.84}
\definecolor{grey85}{rgb}{0.85,0.85,0.85}
\definecolor{grey86}{rgb}{0.86,0.86,0.86}
\definecolor{grey87}{rgb}{0.87,0.87,0.87}
\definecolor{grey88}{rgb}{0.88,0.88,0.88}
\definecolor{grey89}{rgb}{0.89,0.89,0.89}
\definecolor{grey8}{rgb}{0.08,0.08,0.08}
\definecolor{grey90}{rgb}{0.90,0.90,0.90}
\definecolor{grey91}{rgb}{0.91,0.91,0.91}
\definecolor{grey92}{rgb}{0.92,0.92,0.92}
\definecolor{grey93}{rgb}{0.93,0.93,0.93}
\definecolor{grey94}{rgb}{0.94,0.94,0.94}
\definecolor{grey95}{rgb}{0.95,0.95,0.95}
\definecolor{grey96}{rgb}{0.96,0.96,0.96}
\definecolor{grey97}{rgb}{0.97,0.97,0.97}
\definecolor{grey98}{rgb}{0.98,0.98,0.98}
\definecolor{grey99}{rgb}{0.99,0.99,0.99}
\definecolor{grey9}{rgb}{0.09,0.09,0.09}
\definecolor{grey}{rgb}{0.75,0.75,0.75}
\definecolor{honeydew1}{rgb}{0.94,1.00,0.94}
\definecolor{honeydew2}{rgb}{0.88,0.93,0.88}
\definecolor{honeydew3}{rgb}{0.76,0.80,0.76}
\definecolor{honeydew4}{rgb}{0.51,0.55,0.51}
\definecolor{honeydew}{rgb}{0.94,1.00,0.94}
\definecolor{hotpink}{rgb}{1.00,0.41,0.71}
\definecolor{indianred}{rgb}{0.80,0.36,0.36}
\definecolor{ivory1}{rgb}{1.00,1.00,0.94}
\definecolor{ivory2}{rgb}{0.93,0.93,0.88}
\definecolor{ivory3}{rgb}{0.80,0.80,0.76}
\definecolor{ivory4}{rgb}{0.55,0.55,0.51}
\definecolor{ivory}{rgb}{1.00,1.00,0.94}
\definecolor{khaki1}{rgb}{1.00,0.96,0.56}
\definecolor{khaki2}{rgb}{0.93,0.90,0.52}
\definecolor{khaki3}{rgb}{0.80,0.78,0.45}
\definecolor{khaki4}{rgb}{0.55,0.53,0.31}
\definecolor{khaki}{rgb}{0.94,0.90,0.55}
\definecolor{lavenderblush}{rgb}{1.00,0.94,0.96}
\definecolor{lavender}{rgb}{0.90,0.90,0.98}
\definecolor{lawngreen}{rgb}{0.49,0.99,0.00}
\definecolor{lemonchiffon}{rgb}{1.00,0.98,0.80}
\definecolor{lightblue}{rgb}{0.68,0.85,0.90}
\definecolor{lightcoral}{rgb}{0.94,0.50,0.50}
\definecolor{lightcyan}{rgb}{0.88,1.00,1.00}
\definecolor{lightgoldenrod}{rgb}{0.93,0.87,0.51}
\definecolor{lightgoldenrod}{rgb}{0.98,0.98,0.82}
\definecolor{lightgray}{rgb}{0.83,0.83,0.83}
\definecolor{lightgreen}{rgb}{0.56,0.93,0.56}
\definecolor{lightgrey}{rgb}{0.83,0.83,0.83}
\definecolor{lightpink}{rgb}{1.00,0.71,0.76}
\definecolor{lightsalmon}{rgb}{1.00,0.63,0.48}
\definecolor{lightsea}{rgb}{0.13,0.70,0.67}
\definecolor{lightsky}{rgb}{0.53,0.81,0.98}
\definecolor{lightslate}{rgb}{0.47,0.53,0.60}
\definecolor{lightslate}{rgb}{0.47,0.53,0.60}
\definecolor{lightslate}{rgb}{0.52,0.44,1.00}
\definecolor{lightsteel}{rgb}{0.69,0.77,0.87}
\definecolor{lightyellow}{rgb}{1.00,1.00,0.88}
\definecolor{limegreen}{rgb}{0.20,0.80,0.20}
\definecolor{linen}{rgb}{0.98,0.94,0.90}
\definecolor{magenta1}{rgb}{1.00,0.00,1.00}
\definecolor{magenta2}{rgb}{0.93,0.00,0.93}
\definecolor{magenta3}{rgb}{0.80,0.00,0.80}
\definecolor{magenta4}{rgb}{0.55,0.00,0.55}
\definecolor{magenta}{rgb}{1.00,0.00,1.00}
\definecolor{maroon1}{rgb}{1.00,0.20,0.70}
\definecolor{maroon2}{rgb}{0.93,0.19,0.65}
\definecolor{maroon3}{rgb}{0.80,0.16,0.56}
\definecolor{maroon4}{rgb}{0.55,0.11,0.38}
\definecolor{maroon}{rgb}{0.69,0.19,0.38}
\definecolor{mediumaquamarine}{rgb}{0.40,0.80,0.67}
\definecolor{mediumblue}{rgb}{0.00,0.00,0.80}
\definecolor{mediumorchid}{rgb}{0.73,0.33,0.83}
\definecolor{mediumpurple}{rgb}{0.58,0.44,0.86}
\definecolor{mediumsea}{rgb}{0.24,0.70,0.44}
\definecolor{mediumslate}{rgb}{0.48,0.41,0.93}
\definecolor{mediumspring}{rgb}{0.00,0.98,0.60}
\definecolor{mediumturquoise}{rgb}{0.28,0.82,0.80}
\definecolor{mediumviolet}{rgb}{0.78,0.08,0.52}
\definecolor{midnightblue}{rgb}{0.10,0.10,0.44}
\definecolor{mintcream}{rgb}{0.96,1.00,0.98}
\definecolor{mistyrose}{rgb}{1.00,0.89,0.88}
\definecolor{moccasin}{rgb}{1.00,0.89,0.71}
\definecolor{navajowhite}{rgb}{1.00,0.87,0.68}
\definecolor{navyblue}{rgb}{0.00,0.00,0.50}
\definecolor{navy}{rgb}{0.00,0.00,0.50}
\definecolor{oldlace}{rgb}{0.99,0.96,0.90}
\definecolor{olivedrab}{rgb}{0.42,0.56,0.14}
\definecolor{orange1}{rgb}{1.00,0.65,0.00}
\definecolor{orange2}{rgb}{0.93,0.60,0.00}
\definecolor{orange3}{rgb}{0.80,0.52,0.00}
\definecolor{orange4}{rgb}{0.55,0.35,0.00}
\definecolor{orangered}{rgb}{1.00,0.27,0.00}
\definecolor{orange}{rgb}{1.00,0.65,0.00}
\definecolor{orchid1}{rgb}{1.00,0.51,0.98}
\definecolor{orchid2}{rgb}{0.93,0.48,0.91}
\definecolor{orchid3}{rgb}{0.80,0.41,0.79}
\definecolor{orchid4}{rgb}{0.55,0.28,0.54}
\definecolor{orchid}{rgb}{0.85,0.44,0.84}
\definecolor{palegoldenrod}{rgb}{0.93,0.91,0.67}
\definecolor{palegreen}{rgb}{0.60,0.98,0.60}
\definecolor{paleturquoise}{rgb}{0.69,0.93,0.93}
\definecolor{paleviolet}{rgb}{0.86,0.44,0.58}
\definecolor{papayawhip}{rgb}{1.00,0.94,0.84}
\definecolor{peachpuff}{rgb}{1.00,0.85,0.73}
\definecolor{peru}{rgb}{0.80,0.52,0.25}
\definecolor{pink1}{rgb}{1.00,0.71,0.77}
\definecolor{pink2}{rgb}{0.93,0.66,0.72}
\definecolor{pink3}{rgb}{0.80,0.57,0.62}
\definecolor{pink4}{rgb}{0.55,0.39,0.42}
\definecolor{pink}{rgb}{1.00,0.75,0.80}
\definecolor{plum1}{rgb}{1.00,0.73,1.00}
\definecolor{plum2}{rgb}{0.93,0.68,0.93}
\definecolor{plum3}{rgb}{0.80,0.59,0.80}
\definecolor{plum4}{rgb}{0.55,0.40,0.55}
\definecolor{plum}{rgb}{0.87,0.63,0.87}
\definecolor{powderblue}{rgb}{0.69,0.88,0.90}
\definecolor{purple1}{rgb}{0.61,0.19,1.00}
\definecolor{purple2}{rgb}{0.57,0.17,0.93}
\definecolor{purple3}{rgb}{0.49,0.15,0.80}
\definecolor{purple4}{rgb}{0.33,0.10,0.55}
\definecolor{purple}{rgb}{0.63,0.13,0.94}
\definecolor{red1}{rgb}{1.00,0.00,0.00}
\definecolor{red2}{rgb}{0.93,0.00,0.00}
\definecolor{red3}{rgb}{0.80,0.00,0.00}
\definecolor{red4}{rgb}{0.55,0.00,0.00}
\definecolor{red}{rgb}{1.00,0.00,0.00}
\definecolor{rosybrown}{rgb}{0.74,0.56,0.56}
\definecolor{royalblue}{rgb}{0.25,0.41,0.88}
\definecolor{saddlebrown}{rgb}{0.55,0.27,0.07}
\definecolor{salmon1}{rgb}{1.00,0.55,0.41}
\definecolor{salmon2}{rgb}{0.93,0.51,0.38}
\definecolor{salmon3}{rgb}{0.80,0.44,0.33}
\definecolor{salmon4}{rgb}{0.55,0.30,0.22}
\definecolor{salmon}{rgb}{0.98,0.50,0.45}
\definecolor{sandybrown}{rgb}{0.96,0.64,0.38}
\definecolor{seagreen}{rgb}{0.18,0.55,0.34}
\definecolor{seashell1}{rgb}{1.00,0.96,0.93}
\definecolor{seashell2}{rgb}{0.93,0.90,0.87}
\definecolor{seashell3}{rgb}{0.80,0.77,0.75}
\definecolor{seashell4}{rgb}{0.55,0.53,0.51}
\definecolor{seashell}{rgb}{1.00,0.96,0.93}
\definecolor{sienna1}{rgb}{1.00,0.51,0.28}
\definecolor{sienna2}{rgb}{0.93,0.47,0.26}
\definecolor{sienna3}{rgb}{0.80,0.41,0.22}
\definecolor{sienna4}{rgb}{0.55,0.28,0.15}
\definecolor{sienna}{rgb}{0.63,0.32,0.18}
\definecolor{skyblue}{rgb}{0.53,0.81,0.92}
\definecolor{slateblue}{rgb}{0.42,0.35,0.80}
\definecolor{slategray}{rgb}{0.44,0.50,0.56}
\definecolor{slategrey}{rgb}{0.44,0.50,0.56}
\definecolor{snow1}{rgb}{1.00,0.98,0.98}
\definecolor{snow2}{rgb}{0.93,0.91,0.91}
\definecolor{snow3}{rgb}{0.80,0.79,0.79}
\definecolor{snow4}{rgb}{0.55,0.54,0.54}
\definecolor{snow}{rgb}{1.00,0.98,0.98}
\definecolor{springgreen}{rgb}{0.00,1.00,0.50}
\definecolor{steelblue}{rgb}{0.27,0.51,0.71}
\definecolor{tan1}{rgb}{1.00,0.65,0.31}
\definecolor{tan2}{rgb}{0.93,0.60,0.29}
\definecolor{tan3}{rgb}{0.80,0.52,0.25}
\definecolor{tan4}{rgb}{0.55,0.35,0.17}
\definecolor{tan}{rgb}{0.82,0.71,0.55}
\definecolor{thistle1}{rgb}{1.00,0.88,1.00}
\definecolor{thistle2}{rgb}{0.93,0.82,0.93}
\definecolor{thistle3}{rgb}{0.80,0.71,0.80}
\definecolor{thistle4}{rgb}{0.55,0.48,0.55}
\definecolor{thistle}{rgb}{0.85,0.75,0.85}
\definecolor{tomato1}{rgb}{1.00,0.39,0.28}
\definecolor{tomato2}{rgb}{0.93,0.36,0.26}
\definecolor{tomato3}{rgb}{0.80,0.31,0.22}
\definecolor{tomato4}{rgb}{0.55,0.21,0.15}
\definecolor{tomato}{rgb}{1.00,0.39,0.28}
\definecolor{turquoise1}{rgb}{0.00,0.96,1.00}
\definecolor{turquoise2}{rgb}{0.00,0.90,0.93}
\definecolor{turquoise3}{rgb}{0.00,0.77,0.80}
\definecolor{turquoise4}{rgb}{0.00,0.53,0.55}
\definecolor{turquoise}{rgb}{0.25,0.88,0.82}
\definecolor{violetred}{rgb}{0.82,0.13,0.56}
\definecolor{violet}{rgb}{0.93,0.51,0.93}
\definecolor{wheat1}{rgb}{1.00,0.91,0.73}
\definecolor{wheat2}{rgb}{0.93,0.85,0.68}
\definecolor{wheat3}{rgb}{0.80,0.73,0.59}
\definecolor{wheat4}{rgb}{0.55,0.49,0.40}
\definecolor{wheat}{rgb}{0.96,0.87,0.70}
\definecolor{whitesmoke}{rgb}{0.96,0.96,0.96}
\definecolor{white}{rgb}{1.00,1.00,1.00}
\definecolor{yellow1}{rgb}{1.00,1.00,0.00}
\definecolor{yellow2}{rgb}{0.93,0.93,0.00}
\definecolor{yellow3}{rgb}{0.80,0.80,0.00}
\definecolor{yellow4}{rgb}{0.55,0.55,0.00}
\definecolor{yellowgreen}{rgb}{0.60,0.80,0.20}
\definecolor{yellow}{rgb}{1.00,1.00,0.00}
 
\usepackage[usenames,dvipsnames]{xcolor}
\usepackage[psamsfonts]{amsfonts}
\usepackage{graphicx}
\usepackage{amssymb}

\usepackage{hyperref}
\hypersetup{pdfpagemode=FullScreen,   
backref,
colorlinks=true,
citecolor=Bittersweet,
linkcolor=Bittersweet,
urlcolor=Bittersweet   
}

\setattribute{journal}{name}{}

\startlocaldefs

\theoremstyle{plain}
\newtheorem{Lem}{Lemma}[section]

\newtheorem{Theor}{Theorem}[section] 
\newtheorem{Corol}{Corollary}[section]

\numberwithin{equation}{section}

\newcommand{\n}{^{(n)}}

\newcommand{\cqfd}{\hfill $\square$}

\newcommand{\R}{\mathbb R}

\newcommand{\Xb}{\mathbf{X}}
\newcommand{\Sb}{\mathbf{S}}

\newcommand{\Vb}{\mathbf{V}}
\newcommand{\wb}{\mathbf{w}}

\newcommand{\Zb}{\mathbf{Z}}
\newcommand{\zb}{\mathbf{z}}

\newcommand{\eb}{\ensuremath{\mathbf{e}}}
\newcommand{\vb}{\ensuremath{\mathbf{v}}}
\newcommand{\xb}{\ensuremath{\mathbf{x}}}

\newcommand{\Ab}{\ensuremath{\mathbf{A}}}
\newcommand{\Bb}{\ensuremath{\mathbf{B}}}

\newcommand{\Db}{\ensuremath{\mathbf{D}}}

\newcommand{\Wb}{\ensuremath{\mathbf{W}}}
\newcommand{\Yb}{\ensuremath{\mathbf{Y}}}

\newcommand{\Jb}{\ensuremath{\mathbf{J}}}

\newcommand{\thetab}{{\pmb \theta}}

\newcommand{\Lamb}{{\pmb \Lambda}}
\newcommand{\Sigb}{{\pmb \Sigma}}

\newcommand{\Deltab}{{\pmb \Delta}}
\newcommand{\taub}{{\pmb \tau}}

\newcommand{\Gamb}{{\pmb \Gamma}}

\newcommand{\pr}{^{\prime}}

\newcommand{\ny}{n\rightarrow\infty}

\endlocaldefs

\begin{document} 

\begin{frontmatter}
\title{Testing for Principal Component Directions under Weak Identifiability}
\runtitle{Testing for PC Directions under Weak Identifiability}

\begin{aug}
\author{\fnms{Davy} \snm{Paindaveine}\thanksref{t1}\ead[label=e1]{dpaindav@ulb.ac.be}\ead[label=u1,url]{http://homepages.ulb.ac.be/dpaindav}},
\author{\fnms{Julien} \snm{Remy}\ead[label=e2]{Julien.Remy@vub.ac.be}
\ead[label=e2]{Julien.Remy@vub.ac.be}}
\and
\author{\fnms{ Thomas} \snm{Verdebout}\thanksref{t2}
\ead[label=e3]{tverdebo@ulb.ac.be}
\ead[label=u3,url]{http://tverdebo.ulb.ac.be}}

\thankstext{t1}{Corresponding author. Davy Paindaveine's research is supported by a research fellowship from the Francqui Foundation and by the Program of Concerted Research Actions (ARC) of the Universit\'{e} libre de Bruxelles.}
\thankstext{t2}{Thomas Verdebout's research is supported by the Cr\'{e}dit de Recherche J.0134.18 of the FNRS (Fonds National pour la Recherche Scientifique), Communaut\'{e} Fran\c{c}aise de Belgique, and by the aforementioned ARC program of the Universit\'{e} libre de Bruxelles.}
\runauthor{D. Paindaveine, J. Remy and Th. Verdebout}

\affiliation{Universit\'{e} libre de Bruxelles}

\address{Universit\'{e} libre de Bruxelles\\
ECARES and D\'{e}partement de Math\'{e}matique\\
Avenue F.D. Roosevelt, 50\\ 
ECARES, CP114/04\\
B-1050, Brussels\\ 
Belgium\\
\printead{e1}\\  
\printead{u1}\\
}

\address{Universit\'{e} libre de Bruxelles\\
ECARES and D\'{e}partement de Math\'{e}matique\\
Boulevard du Triomphe, CP210\\
B-1050, Brussels\\
Belgium\\ 
\printead{e2}\\ 
\phantom{E-mail:\ }\printead*{e3}\\ 
\printead{u3}\\
}

\end{aug}
\vspace{3mm}

\begin{abstract}
We consider the problem of testing, on the basis of a $p$-variate Gaussian random sample, the null hypothesis~${\cal H}_0: \thetab_1= \thetab_1^0$ against the alternative~${\cal H}_1: \thetab_1 \neq \thetab_1^0$, where~$\thetab_1$ is the ``first" eigenvector of the underlying covariance matrix and~$\thetab_1^0$ is a fixed unit $p$-vector. In the classical setup where eigenvalues~$\lambda_1>\lambda_2\geq \ldots\geq \lambda_p$ are fixed, the \cite{And63} likelihood ratio test (LRT) and the \cite{HPV10} Le Cam optimal test for this problem are asymptotically equivalent under the null hypothesis, hence also under sequences of contiguous alternatives. We show that this equivalence does not survive asymptotic scenarios where~$\lambda_{n1}/\lambda_{n2}=1+O(r_n)$ with~$r_n=O(1/\sqrt{n})$. For such scenarios, the Le Cam optimal test still asymptotically meets the nominal level constraint, whereas the LRT severely overrejects the null hypothesis. Consequently, the former test should be favored over the latter one whenever the two largest sample eigenvalues are close to each other. By relying on the Le Cam's asymptotic theory of statistical experiments, we study the non-null and optimality properties of the Le Cam optimal test in the aforementioned asymptotic scenarios and show that the null robustness of this test is not obtained at the expense of power. Our asymptotic investigation is extensive in the sense that it allows~$r_n$ to converge to zero at an arbitrary rate. While we restrict to single-spiked spectra of the form~$\lambda_{n1}>\lambda_{n2}=\ldots=\lambda_{np}$ to make our results as striking as possible, we extend our results to the more general elliptical case. Finally, we present an illustrative real data example.  
 \end{abstract}



\end{frontmatter}

\section{Introduction} 
\label{sec:intro} 

Principal Component Analysis (PCA) is one of the most classical tools in multivariate statistics. For a random $p$-vector~$\Xb$ with mean zero and a covariance matrix~$\Sigb$ admitting the spectral decomposition~$\Sigb=\sum_{j=1}^p \lambda_j \thetab_j\thetab_j'$ ($\lambda_1\geq \ldots\geq \lambda_p$), the $j$th principal component is~$\thetab_j\pr \Xb$, that is, the projection of~$\Xb$ onto the $j$th unit eigenvector~$\thetab_j$ of~$\Sigb$. In practice,~$\Sigb$ is usually unknown, so that one of the key issues in PCA is to perform inference on eigenvectors. The seminal paper \cite{And63} focused on the multinormal case and derived asymptotic results for the maximum likelihood estimators of the~$\thetab_j$'s and~$\lambda_j$'s. Later, \cite{Tyl81,Tyl83} extended those results to the elliptical case, where, to avoid moment assumptions,~$\Sigb$ is then the corresponding ``scatter" matrix rather than the covariance matrix. Still under ellipticity assumptions, \cite{HPV10} obtained Le Cam optimal tests on eigenvectors and eigenvalues, whereas \cite{HPV14} developed efficient R-estimators for eigenvectors. \cite{Gent00}, \cite{Mia05} and \cite{He11} proposed various robust methods for PCA. Recently, \cite{JoLu09}, \cite{BeRi14} and \cite{FaLiu14} considered inference on eigenvectors of~$\Sigb$ in sparse high-dimensional situations. PCA has also been extensively considered in the functional case; see, e.g., \cite{boente2000}, \cite{bali2011} or the review paper \cite{Cue2014}.
 
In this work, we focus on the problem of testing the null hypothesis ${\cal H}_0: \thetab_1=\thetab_{1}^0$ against the alternative~${\cal H}_1: \thetab_1 \neq \thetab_{1}^0$, where~$\thetab_{1}^0$ is a given unit vector of~$\R^p$. While, strictly speaking, the fact that~$\thetab_1$ will below be defined up to a sign only should lead us to formulate the null hypothesis as~${\cal H}_0: \thetab_1\in\{\thetab_{1}^0,-\thetab_{1}^0\}$, we will stick to the formulation above, which is the traditional one in the literature; we refer to the many references provided below. We restrict to~$\thetab_1$ for the sake of simplicity only; our results could indeed be extended to null hypotheses of the form~$\mathcal{H}_0:\thetab_j=\thetab_j^0$ for any other~$j$.
While the emphasis in PCA is usually more on point estimation, the testing problems above are also of high practical relevance. For instance, they are of paramount importance in \emph{confirmatory PCA}, that is, when it comes to testing that~$\thetab_1$ (or any other~$\thetab_j$) coincides with an eigenvector obtained from an earlier real data analysis (``historical data") or with an eigenvector resulting from a theory or model. In line with this, tests for  the null hypothesis ${\cal H}_0: \thetab_1= \thetab_{1}^0$ have been used in, among others, \cite{Jac05} to analyze the concentration of a chemical component in a solution and in \cite{Sy08} for the study of the geometric similarity in modern humans. 

More specifically, we want to test~${\cal H}_0: \thetab_1= \thetab_{1}^0$ against~${\cal H}_1: \thetab_1 \neq \thetab_{1}^0$ on the basis of a random sample~$\Xb_1,\ldots,\Xb_n$ from the $p$-variate normal distribution with mean~${\pmb \mu}$ and covariance matrix~$\Sigb$ (the extension to elliptical distributions will also be considered). Denoting as $\hat{\lambda}_{1} \geq \hat{\lambda}_{2} \geq \ldots\geq \hat{\lambda}_{p}$ the eigenvalues of the sample covariance matrix~$
{\bf S}
:=
\frac{1}{n}
\sum_{i=1}^n
\,
(\Xb_i-\bar{\Xb}) (\Xb_i-\bar{\Xb})'
$
(as usual, $\bar{\Xb}:=\frac{1}{n}\sum_{i=1}^n \Xb_i$ here), the classical test for this problem is the \cite{And63} likelihood ratio test, $\phi_{\rm A}$ say, rejecting the null hypothesis at asymptotic level~$\alpha$ when 
$$
Q_{\rm A}
:=
n \big( \hat{\lambda}_{1} \thetab_{1}^{0\prime} {\bf S}^{-1} \thetab_{1}^0 +  \hat{\lambda}_{1}^{-1} \thetab_{1}^{0\prime} {\bf S}\, \thetab_{1}^0-2\big) > \chi_{p-1, 1-\alpha}^2, 
$$
where~$\chi_{\ell, 1-\alpha}^2$ stands for the $\alpha$-upper quantile of the chi-square distribution with $\ell$ degrees of freedom. Various extensions of this test have been proposed in the literature: to mention only a few, \cite{Jo84} considered a small-sample test, \cite{Flu88} proposed an extension to more eigenvectors, \cite{Tyl81,Tyl83} robustified the test to possible (elliptical) departures from multinormality, while \cite{Sc08} considered extensions to the case of Gaussian random matrices. More recently, \cite{HPV10} obtained the Le Cam optimal test for the problem above. This test, $\phi_{\rm HPV}$ say, rejects  the null hypothesis at asymptotic level $\alpha$ when
$$
Q_{\rm HPV}
:=
 \frac{n}{\hat{\lambda}_{1}}  \sum_{j=2}^p \hat{\lambda}_{j}^{-1} 
 \big(\tilde{\thetab}_j\pr {\bf S} \thetab_{1}^0\big)^2 > \chi_{p-1, 1-\alpha}^2,
$$
where~$\tilde{\thetab}_j$, $j=2,\ldots,p$, defined recursively through 
\begin{equation}
\label{GS}
\tilde{\thetab}_j
:=
\frac{({\bf I}_p- \thetab_{1}^0\thetab_{1}^{0\prime}- \sum_{k=2}^{j-1}\tilde{\thetab}_k \tilde{\thetab}_k\pr)\hat{\thetab}_j}{\| ({\bf I}_p- \thetab_{1}^0\thetab_{1}^{0\prime}- \sum_{k=2}^{j-1}\tilde{\thetab}_k \tilde{\thetab}_k\pr)\hat{\thetab}_j \|}
\end{equation}
(with summation over an empty collection of indices being equal to zero), result from a Gram-Schmidt orthogonalization of~$\thetab_{1}^0,\hat{\thetab}_2, \ldots, \hat{\thetab}_p$, where~$\hat\thetab_j$ is a unit eigenvector of~$\Sb$ associated with the eigenvalue~$\hat\lambda_j$, $j=2,\ldots,p$. When the eigenvalues of~$\Sigb$ are fixed and satisfy $\lambda_1>\lambda_2 \geq \lambda_3\geq \ldots \geq \lambda_p$ (the minimal condition under which~$\thetab_1$ is identified---up to an unimportant sign, as already mentioned), both tests above are asymptotically equivalent under the null hypothesis, hence also under sequences of contiguous alternatives, which implies that~$\phi_{\rm A}$ is also Le Cam optimal; see \cite{HPV10}. The tests~$\phi_{\rm A}$ and $\phi_{\rm HPV}$ can therefore be considered perfectly equivalent, at least asymptotically so. 
 
In the present paper, we compare the asymptotic behaviors of these tests in a non-standard asymptotic framework where eigenvalues may depend on~$n$ and where~$\lambda_{n1}/\lambda_{n2}$ converges to~$1$ as~$n$ diverges to infinity. Such asymptotic scenarios provide \emph{weak identifiability} since the first eigenvector~$\thetab_1$ is not properly identified in the limit. To make our results as striking as possible, we will 
restrict to single-spiked spectra of the form~$\lambda_{n1}>\lambda_{n2}=\ldots=\lambda_{np}$. In other words, we will consider triangular arrays of observations~$\Xb_{ni}$, $i=1,\ldots,n$, $n=1,2,\ldots$, where $\Xb_{n1},\ldots,\Xb_{nn}$ form a random sample from the $p$-variate normal distribution with mean~$\pmb \mu_n$ and covariance matrix 
\begin{eqnarray} 
\label{localspiked}
\Sigb_n 
&\!\!:=\!\!&
 \sigma^2_n ({\bf I}_p+ r_n v\, \thetab_1\thetab_1\pr) \nonumber \\
&\!\!=\!\!&
 \sigma^2_n (1+r_n v) \thetab_1 \thetab_1\pr + \sigma^2_n ({\bf I}_p-\thetab_1\thetab_1\pr)
,
\end{eqnarray}
where~$v$ is a positive real number, $(\sigma_n)$ is a positive real sequence, $(r_n)$ is a bounded positive real sequence, and~${\bf I}_\ell$ denotes the $\ell$-dimensional identity matrix (again, the multinormality assumption will be relaxed later in the paper). The eigenvalues of the covariance matrix~$\Sigb_n$ are then~$\lambda_{n1}=\sigma^2_n(1+r_n v)$ (with corresponding eigenvector~$\thetab_1$) and~$\lambda_{n2}=\ldots=\lambda_{np}=\sigma^2_n$ (with corresponding eigenspace being the orthogonal complement of~$\thetab_1$ in~$\R^p$). If~$r_n\equiv 1$ (or more generally if~$r_n$ stays away from~$0$ as $\ny$), then this setup is similar to the classical one where the first eigenvector~$\thetab_1$ remains identified in the limit. In contrast,
\vspace{-.5mm}
  if~$r_n=o(1)$, then the resulting weak identifiability
intuitively makes the problem of testing~${\cal H}_0\n: \thetab_1= \thetab_{1}^0$ against~${\cal H}_1\n: \thetab_1 \neq \thetab_{1}^0$ increasingly hard as~$n$ diverges to infinity.

Our results show that, while they are, as mentioned above, equivalent in the standard asymptotic scenario associated with~$r_n\equiv 1$, the tests~$\phi_{\rm HPV}$ and~$\phi_{\rm A}$ actually exhibit very different behaviors under weak identifiability. More precisely, we show that this asymptotic equivalence survives scenarios where~$r_n=o(1)$ with~$\sqrt{n}r_n\to \infty$, but not scenarios where~$r_n=O(1/\sqrt{n})$. Irrespective of the asymptotic scenario considered, the test~$\phi_{\rm HPV}$ asymptotically meets the nominal level constraint, hence may be considered robust to weak identifiability. On the contrary, in scenarios where~$r_n=O(1/\sqrt{n})$, the test~$\phi_{\rm A}$ dramatically overrejects the null hypothesis. Consequently, despite the asymptotic equivalence of these tests in standard asymptotic scenarios, the test~$\phi_{\rm HPV}$ should be favored over~$\phi_{\rm A}$.
 
Of course, this nice robustness property of~$\phi_{\rm HPV}$ only refers to the null asymptotic behavior of this test, and it is of interest to investigate whether or not this null robustness is obtained at the expense of power.  In order to do so, we study the non-null and optimality properties of~$\phi_{\rm HPV}$ under suitable local alternatives. This is done by exploiting the Le Cam's asymptotic theory of statistical experiments. In every asymptotic scenario considered, we show that the corresponding sequence of experiments converges to a limiting experiment in the Le Cam sense. Interestingly, (i) the corresponding contiguity rate crucially depends on the underlying asymptotic  scenario and (ii) the resulting limiting experiment is not always a Gaussian shift experiment, such as in the standard \emph{local asymptotic normality (LAN)} setup. In all cases, however, we can derive the asymptotic non-null distribution of~$Q_{\rm HPV}$ under contiguous alternatives by resorting to the Le Cam third lemma, and we can establish that this test enjoys excellent optimality properties.     

The problem we consider in this paper is characterized by the fact that the parameter of interest (here, the first eigenvector) is unidentified when a nuisance parameter is equal to some given value (here, when the ratio of both largest eigenvalues is equal to one). Such situations have already been  considered in the statistics and econometrics literatures; we refer, e.g.,  to \cite{Duf1997}, \cite{Pot2002}, \cite{ForHil2003}, \cite{Duf2006}, or \cite{PV17}. To the best of our knowledge, however, no results have been obtained in PCA under weak identifiability. We think that, far from being of academic interest only, our results are also crucial for practitioners: they indeed provide a clear warning that, when the underlying distribution is close to spherical (more generally, when both largest sample eigenvalues are nearly equal), the daily-practice Gaussian test~$\phi_{\rm A}$ tends to overreject the null hypothesis, hence may lead to wrong conclusions (false positives) with very high probability, whereas the test~$\phi_{\rm HPV}$ remains a reliable procedure in such cases. We provide an illustrative real data example that shows the practical relevance of our results. 

The paper is organized as follows. In Section~\ref{sec:preliminaries}, we introduce the distributional setup and notation to be used throughout and we derive preliminary results on the asymptotic behavior of sample eigenvalues/eigenvectors. In Section~\ref{sec:Null}, we show that the null asymptotic distribution of~$Q_{\rm HPV}$ is~$\chi^2_{p-1}$ under all asymptotic scenarios, whereas that of~$Q_{\rm A}$ is~$\chi^2_{p-1}$ only if~$\sqrt{n}r_n\to\infty$. We also explicitly provide the null asymptotic distribution of~$Q_{\rm A}$ when~$r_n=O(1/\sqrt{n})$. In Section~\ref{sec:LAN}, we show that, in all asymptotic scenarios, the sequence of experiments considered converges to a limiting experiment. Then, this is used to study the non-null and optimality properties of~$\phi_{\rm HPV}$. In Section~\ref{sec:ellipt}, we extend our results to the more general elliptical case. Theoretical findings in Sections~\ref{sec:Null} to~\ref{sec:ellipt} are illustrated through Monte Carlo exercises. We treat a real data illustration in Section~\ref{sec:real}. Finally, we wrap up and shortly discuss research perspectives in Section~\ref{sec:wrapup}. All proofs are provided in the appendix.  


\section{Preliminary results}
\label{sec:preliminaries}

As mentioned above, we will consider throughout triangular arrays of observations~$\Xb_{ni}$, $i=1,\ldots,n$,  $n=1,2,\ldots$, where $\Xb_{n1},\ldots,\Xb_{nn}$ form a random sample from the $p$-variate normal distribution with mean~${\pmb \mu}_n$ and covariance matrix~$\Sigb_n=\sigma^2_n ({\bf I}_p+ r_n v\, \thetab_1\thetab_1\pr)$, where~$\thetab_1$ is
 a unit $p$-vector and~$\sigma_n,r_n$ and~$v$ are
\vspace{-.3mm}
 positive real numbers. The resulting hypothesis will be denoted 
\vspace{-.2mm}
as~${\rm P}_{{\pmb \mu}_n,\sigma_n,\thetab_1,r_n,v}={\rm P}\n_{{\pmb \mu}_n,\sigma_n,\thetab_1,r_n,v}$ (the superscript~$^{(n)}$ will be dropped in the sequel).
Throughout,~$\bar{\Xb}_n:= \frac{1}{n} \sum_{i=1}^n \Xb_{ni}$ and~${\bf S}_n:= \frac{1}{n} \sum_{i=1}^n (\Xb_{ni}-\bar{\Xb}_n) (\Xb_{ni}-\bar{\Xb}_n)\pr$ will denote the sample average and sample covariance matrix of~$\Xb_{n1},\ldots,\Xb_{nn}$. For any~$j=1,\ldots,p$, the $j$th largest eigenvalue of~$\Sb_n$ and ``the" corresponding unit eigenvector will be
  denoted as~$\hat{\lambda}_{nj}$ and~$\hat{\thetab}_{nj}$, respectively 
\vspace{-.8mm}
(identifiability is discussed at the end of this paragraph). With this notation, the tests~$\phi_{\rm A}=\phi_{\rm A}\n$ and~$\phi_{\rm HPV}=\phi_{\rm HPV}\n$ from the introduction reject the null hypothesis at asymptotic level~$\alpha$ when 
\begin{eqnarray} 
\nonumber
Q_{\rm A}
=
Q_{\rm A}\n
&\!\!\!=\!\!\!&
n \big( \hat{\lambda}_{n1} \thetab_{1}^{0\prime} {\bf S}_n^{-1} \thetab_{1}^0 +  \hat{\lambda}_{n1}^{-1} \thetab_{1}^{0\prime} {\bf S}_n \thetab_{1}^0 -2\big) \\[2mm]
&\!\!\!=\!\!\!&
\frac{n}{\hat{\lambda}_{n1}}
 \sum_{j=2}^p \hat{\lambda}_{nj}^{-1}(\hat{\lambda}_{n1}-\hat{\lambda}_{nj})^2  \big(\hat{\thetab}_{nj}\pr \thetab_1^0 \big)^2 
> \chi_{p-1, 1-\alpha}^2 
\label{Andstat}
\end{eqnarray}
and
\begin{equation}
\label{HPV2010} 
Q_{\rm HPV}
=
Q_{\rm HPV}\n
=
 \frac{n}{\hat{\lambda}_{n1}}  \sum_{j=2}^p \hat{\lambda}_{nj}^{-1} 
 \big(\tilde{\thetab}_{nj}\pr {\bf S}_n \thetab_{1}^0 \big)^2
  > \chi_{p-1, 1-\alpha}^2
,
\end{equation}
respectively, 
where~$\tilde{\thetab}_{n2},\ldots,\tilde{\thetab}_{np}$ result from the Gram-Schmidt orthogonalization in~(\ref{GS}) 
\vspace{-.5mm}
 applied to~$\thetab_{1}^0,\hat{\thetab}_{n2}, \ldots, \hat{\thetab}_{np}$. Under~${\rm P}_{{\pmb \mu}_n,\sigma_n,\thetab_1,r_n,v}$, the sample eigenvalue~$\hat{\lambda}_{nj}$ is uniquely defined with probability one, but~$\hat{\thetab}_{nj}$ is, still with probability one, defined up to a sign only. Clearly, this sign does not play any role in~(\ref{Andstat})--(\ref{HPV2010}), hence will be fixed arbitrarily. At a few places below, however, this sign will need to be fixed in an appropriate way.

For obvious reasons, the asymptotic behavior of~$\Sb_n$ will play a crucial role when investigating the asymptotic properties of the tests above. To describe this behavior, we need to introduce the following notation. For an $\ell\times \ell$ matrix~$\Ab$, denote as~${\rm vec} ({\bf A})$ the vector obtained by stacking the columns of~${\bf A}$ on top of each other. We will let~$\Ab^{\otimes 2}:=\Ab\otimes\Ab$, where ${\bf A} \otimes {\bf B}$ is the Kronecker product of~${\bf A}$ and~${\bf B}$. The \emph{commutation matrix}~${\bf K}_{k,\ell}$, that is such that~${\bf K}_{k,\ell}({\rm vec}\,\Ab)={\rm vec}(\Ab')$ for any~$k\times \ell$ matrix~$\Ab$, satisfies ${\bf K}_{p,k}({\bf A} \otimes {\bf B})=({\bf B} \otimes {\bf A}){\bf K}_{q,\ell}$ for any~$k\times \ell$ matrix~${\bf A}$ and $p\times q$ matrix~${\bf B}$; see, e.g., \citet{MagNeu2007}. If~$\Xb$ is $p$-variate standard normal, then the covariance matrix of~${\rm vec}(\Xb\Xb')$ is~$\mathbf{I}_{p^2}+{\bf K}_p$, with~${\bf K}_p:={\bf K}_{p,p}$; the Levy-Lindeberg central limit theorem then easily provides the following result.

\begin{Lem}
\label{LemLL}
Fix a unit $p$-vector~$\thetab_1$, $v>0$ and a bounded positive real sequence~$(r_n)$. Then, under ${\rm P}_{{\bf 0},1,\thetab_1,r_n,v}$, 
	$\sqrt{n}(\Sigb_n^{-1/2})^{\otimes 2} {\rm vec}
	({\bf S}_n- \Sigb_n)$ is asymptotically normal with mean zero and covariance matrix~${\bf I}_{p^2}+ {\bf K}_p$. In particular, (i) if $r_n\equiv 1$, then $\sqrt{n} \, {\rm vec}
	({\bf S}_n- \Sigb_n)$ is asymptotically normal with mean zero and covariance matrix~$({\bf I}_{p^2}+ {\bf K}_p)(\Sigb(v))^{\otimes 2}$, with~$\Sigb(v):={\bf I}_p+ v\, \thetab_1\thetab_1\pr$; 
	(ii) if~$r_n$ is~$o(1)$, then~$\sqrt{n} \, {\rm vec}
	({\bf S}_n- \Sigb_n)$ is asymptotically normal with mean zero and covariance matrix~${\bf I}_{p^2}+ {\bf K}_p$.
\end{Lem}

Clearly, the tests~$\phi_{\rm A}$ and~$\phi_{\rm HPV}$ above are invariant under translations  and scale transformations, that is, respectively, under transformations of the form $(\Xb_{n1},\ldots,\Xb_{nn}) \mapsto (\Xb_{n1}+{\bf t},\ldots,\Xb_{nn}+{\bf t})$, with~${\bf t}\in\R^p$, and $(\Xb_{n1},\ldots,\Xb_{nn}) \mapsto (s\Xb_{n1},\ldots,s\Xb_{nn})$, with~$s>0$. This implies that, when investigating the behavior of these tests, we may assume without loss of generality that~${\pmb \mu}_n\equiv {\bf 0}$ and~$\sigma_n\equiv 1$, that is, we may restrict to hypotheses of the form~${\rm P}_{\thetab_1,r_n,v}:={\rm P}_{{\bf 0},1,\thetab_1,r_n,v}$, as we already did in Lemma~\ref{LemLL}. We therefore restrict to such hypotheses in the rest of the paper. 

The tests~$\phi_{\rm A}$ and~$\phi_{\rm HPV}$ are based on statistics that do not only involve the sample covariance matrix~$\Sb_n$, but also the corresponding sample eigenvalues and eigenvectors. It is therefore no surprise that investigating the asymptotic behavior of these tests under weak identifiability will require controlling the asymptotic behaviors of sample eigenvalues and eigenvectors. For eigenvalues, we have the following result (throughout,~${\rm diag}({\bf A}_1, \ldots, {\bf A}_m)$ stands for the block-diagonal matrix with diagonal blocks~${\bf A}_1, \ldots, {\bf A}_m$).

\begin{Lem} 
\label{Lemeigenvalues}
Fix a unit $p$-vector~$\thetab_1$, $v> 0$ and a bounded  positive real sequence~$(r_n)$. 
Let~$\Zb(v)$ be a $p\times p$ random matrix such that~${\rm vec}({\bf Z}(v))\sim {\cal N}({\bf 0}, ({\bf I}_{p^2}+ {\bf K}_p)(\Lamb(v))^{\otimes 2})$, with~$\Lamb(v):={\rm diag}(1+v,1,\ldots,1)$, and let~$\Zb_{22}(v)$ be the matrix obtained from~$\Zb(v)$ by deleting its first row and first column. Write~$\Zb:=\Zb(0)$ and~$\Zb_{22}:=\Zb_{22}(0)$. Then, under~${\rm P}_{\thetab_1,r_n,v}$, 
$$
{\pmb \ell}_n:=\big(\sqrt{n}(\hat{\lambda}_{n1}- (1+r_nv)), \sqrt{n}(\hat{\lambda}_{n2}-1), \ldots,  \sqrt{n}(\hat{\lambda}_{np}-1) \big)\pr
\stackrel{\mathcal{D}}{\to}
{\pmb \ell}=(\ell_1, \ldots,\ell_p)'
,
$$
where~$\stackrel{\mathcal{D}}{\to}$ denotes weak convergence and where~${\pmb \ell}$ is as follows:
\begin{itemize}
\item[(i)] if $r_n\equiv 1$, then~$\ell_1$ and $(\ell_2, \ldots, \ell_p)'$ are mutually independent,~$\ell_1$ is normal with mean zero and variance~$2(1+v)^2$, and $\ell_2\geq  \ldots\geq \ell_p$ are the eigenvalues of~${\bf Z}_{22}(v)$;
\vspace{1mm} 
\item[(ii)] if $r_n$ is $o(1)$ with $\sqrt{n} r_n  \to \infty$, then~$\ell_1$ and $(\ell_2, \ldots, \ell_p)'$ are mutually independent,~$\ell_1$ is normal with mean zero and variance~$2$, and $\ell_2\geq  \ldots\geq \ell_p$ are the eigenvalues of~${\bf Z}_{22}$;
\vspace{1mm} 
\item[(iii)] if $r_n= 1/\sqrt{n}$, then~$\ell_1$ is the largest eigenvalue of~${\bf Z}-{\rm diag}(0, v,\ldots,v)$ and $\ell_2\geq \ldots\geq \ell_p$ are the $p-1$ smallest eigenvalues of~${\bf Z}+{\rm diag}(v, 0,
\linebreak
 \ldots, 0)$;
\vspace{1mm} 
\item[(iv)] If $r_n=o(1/\sqrt{n})$, then~${\pmb \ell}$ is the vector of eigenvalues of~${\bf Z}$ $($in decreasing order$)$, hence has density 
\begin{equation}
\label{explidensl}	
(\ell_1, \ldots, \ell_p)' 
\mapsto 
b_p \; 
{\exp} \bigg( \!\! -\frac{1}{4} \sum_{j=1}^p \ell_j^2 \bigg) 
\bigg( \prod_{1\leq k<j \leq p} (\ell_k- \ell_j) \bigg)
\,
\mathbb{I}[\ell_1
\geq\ldots\geq \ell_p],
\end{equation}
where~$b_p$ is a normalizing constant and where~$\mathbb{I}[A]$ is the indicator function of~$A$.
\end{itemize}
\end{Lem} 

Lemma~\ref{Lemeigenvalues} shows that, unlike the sample covariance matrix~$\Sb_n$, sample eigenvalues 
exhibit an asymptotic behavior that crucially depends on~($r_n$). The important threshold, associated with $r_n=1/\sqrt{n}$, provides sequences of hypotheses~${\rm P}_{\thetab_1,r_n,v}$ that are contiguous to the spherical hypotheses~${\rm P}_{\thetab_1,0,v}$ under which the first eigenvector~$\thetab_1$ is unidentified (contiguity follows, e.g., from Proposition~2.1 in \citealp{HalPai2006}). Lemma~\ref{Lemeigenvalues} then identifies four regimes that will be present throughout our double asymptotic investigation below, namely \emph{away from contiguity} ($r_n\equiv 1$, case (i)), \emph{above contiguity} ($r_n$ is $o(1)$ with $\sqrt{n} r_n \to \infty$, case (ii)), \emph{under contiguity} ($r_n=1/\sqrt{n}$, case (iii)), and \emph{under strict contiguity} ($r_n=o(1/\sqrt{n})$, case (iv)).

In the high-dimensional setup where~$p=p_n$ is such that~$p/n\to \gamma^{-2}\in(0,1]$, a related phase transition phenomenon has been identified in \cite{Baietal2005}, in the case of complex-valued Gaussian observations. More precisely, in the single-spiked case
considered in the present paper, Theorem~1.1 of \cite{Baietal2005} proves that the asymptotic behavior of~$\hat{\lambda}_{n1}$ crucially depends on the ratio~$\rho$ of~$\lambda_{n1}$ to the common value of~$\lambda_{nj}$, $j=2,\ldots,p$; there, $\rho$ is essentially of the form~$\rho=1+C\sqrt{p/n}(\to 1+C\gamma)$, for some constant~$C\geq 1$ whose value is showed to strongly impact the weak limit and consistency rate of~$\hat{\lambda}_{n1}$. Note that, in contrast, several rates are considered for $\rho=\rho_n$ in Lemma~\ref{Lemeigenvalues} above, and that~$\hat{\lambda}_{n1}$ exhibits the same consistency rate in each case. 

While Lemmas~\ref{LemLL}--\ref{Lemeigenvalues} will be sufficient to study the asymptotic behavior of~$\phi_{\rm HPV}$, 
\vspace{-.3mm} 
the test~$\phi_{\rm A}$, as hinted by the expression in~(\ref{Andstat}), will further require investigating the joint asymptotic behavior of~$\hat{\thetab}_{n2}\pr \thetab_1^0, \ldots, \hat{\thetab}_{np}\pr \thetab_1^0$. To do so, fix arbitrary $p$-vectors~$\thetab_2, \ldots, \thetab_p$ such that~${\Gamb}:=(\thetab_1^0, \thetab_2, \ldots, \thetab_p)$ is orthogonal. Let further~$\hat{\Gamb}_n:=(\hat{\thetab}_{n1}, \ldots, \hat{\thetab}_{np})$,
where the ``signs" of~$\hat{\thetab}_{nj}$, $j=1,\ldots,p$, are fixed by the constraint that, with probability one, all entries in the first column of
\begin{eqnarray} 
\label{defE} 
{\bf E}_n
:= 
\hat{\Gamb}_n\pr\Gamb =\left( \begin{array}{cc} E_{n,11} & {\bf E}_{n,12} \\ 
{\bf E}_{n,21} & {\bf E}_{n,22} \end{array} \right)
\end{eqnarray}
are positive (note that all 
\vspace{-.5mm}
entries of~${\bf E}_n$ are almost surely non-zero). 
With this notation, ${\bf E}_{n,21}$ collects the random variables~$\hat{\thetab}_{n2}\pr \thetab_1^0, \ldots, \hat{\thetab}_{np}\pr \thetab_1^0$ of interest above. We then have the following result. 

%
%
%
%
%
%

\begin{Lem} 
\label{Lemeigenvectors}
Fix a unit $p$-vector~$\thetab_1$, $v>0$ and a bounded positive real sequence~$(r_n)$. 
Let~$\Zb$ be a $p\times p$ random matrix such that~${\rm vec}({\bf Z})\sim {\cal N}({\bf 0}, {\bf I}_{p^2}+ {\bf K}_p)$. Let ${\bf E}(v):=(\wb_1(v),\ldots,\wb_p(v))'$, where~$\wb_j(v)=(w_{j1}(v),\ldots,w_{jp}(v))'$ is the unit eigenvector associated with the $j$th largest eigenvalue of~${\bf Z}+{\rm diag}(v,0,\ldots,0)$ and such that~$w_{j1}(v)>0$ almost surely. Extending the definitions to the case~$v=0$, write~${\bf E}:={\bf E}(0)$. Then, we have the following 
under~${\rm P}_{\thetab_1,r_n,v}$:
\begin{itemize} 
\item[(i)]  
if $r_n\equiv 1$, then~$E_{n,11}=1+o_{\rm P}(1)$, ${\bf E}_{n,22}{\bf E}_{n,22}\pr=\mathbf{I}_{p-1}+o_{\rm P}(1)$, $\sqrt{n} {\bf E}_{n,21}=O_{\rm P}(1)$, and both~$\sqrt{n}{\bf E}_{n,22}\pr{\bf E}_{n,21}$ and~$\sqrt{n}{\bf E}_{n,12}'$ are asymptotically normal with mean zero and covariance matrix~$v^{-2}(1+v) {\bf I}_{p-1};$ 
\vspace{1mm} 
\item[(ii)]   
if~$r_n$ is $o(1)$ with $\sqrt{n} r_n \to \infty$, then~$E_{n,11}=1+o_{\rm P}(1)$, ${\bf E}_{n,22}{\bf E}_{n,22}\pr=\mathbf{I}_{p-1}+o_{\rm P}(1)$, $\sqrt{n} r_n {\bf E}_{n,21}=O_{\rm P}(1)$, and both~$\sqrt{n}r_n{\bf E}_{n,22}\pr{\bf E}_{n,21}$ and \linebreak
$\sqrt{n}r_n{\bf E}_{n,12}'$ are asymptotically normal with mean zero and covariance matrix~$v^{-2}{\bf I}_{p-1};$
\vspace{1mm} 
\item[(iii)]  
if~$r_n=1/\sqrt{n}$, then~${\bf E}_n$ converges weakly to~${\bf E}(v);$
\vspace{1mm} 
\item[(iv)] if~$r_n=o(1/\sqrt{n})$, then ${\bf E}_n$ converges weakly to~${\bf E}$.  
\end{itemize}
\end{Lem}
 
This result shows that the asymptotic behavior of~${\bf E}_{n,21}$, which, as mentioned above, is the only part of~${\bf E}_n$ involved in the Anderson test statistic~$Q_{\rm A}$, depends on the regimes identified in Lemma~\ref{Lemeigenvalues}. (i) Away from contiguity, ${\bf E}_{n,21}$ converges to the zero vector in probability at the standard root-$n$ rate. (ii) Above contiguity, ${\bf E}_{n,21}$ is still $o_{\rm P}(1)$, but the rate of convergence is slower. (iii) Under contiguity, consistency is lost and~${\bf E}_{n,21}$ converges weakly to a distribution that still depends on~$v$. (iv) Under strict contiguity, on the contrary, the asymptotic distribution of~${\bf E}_n$ does not depend on~$v$ and inspection of the proof of Lemma~\ref{Lemeigenvectors} reveals that this asymptotic distribution is the same as the one we would obtain for~$v=0$. In other words, the asymptotic distribution of~${\bf E}_n$ is then the same as in the spherical Gaussian case, so that, up to the fact that~$\bf E$ has almost surely positive entries in its first column (a constraint inherited from the corresponding one on~${\bf E}_n$), this asymptotic distribution is the invariant Haar distribution on the group of~$p\times p$ orthogonal matrices; see \cite{And63}, page~126.


%
%
%
%
%
%

\section{Null results}
\label{sec:Null}

\begin{figure}[htbp!]
\begin{center} 
\includegraphics[width=\textwidth]{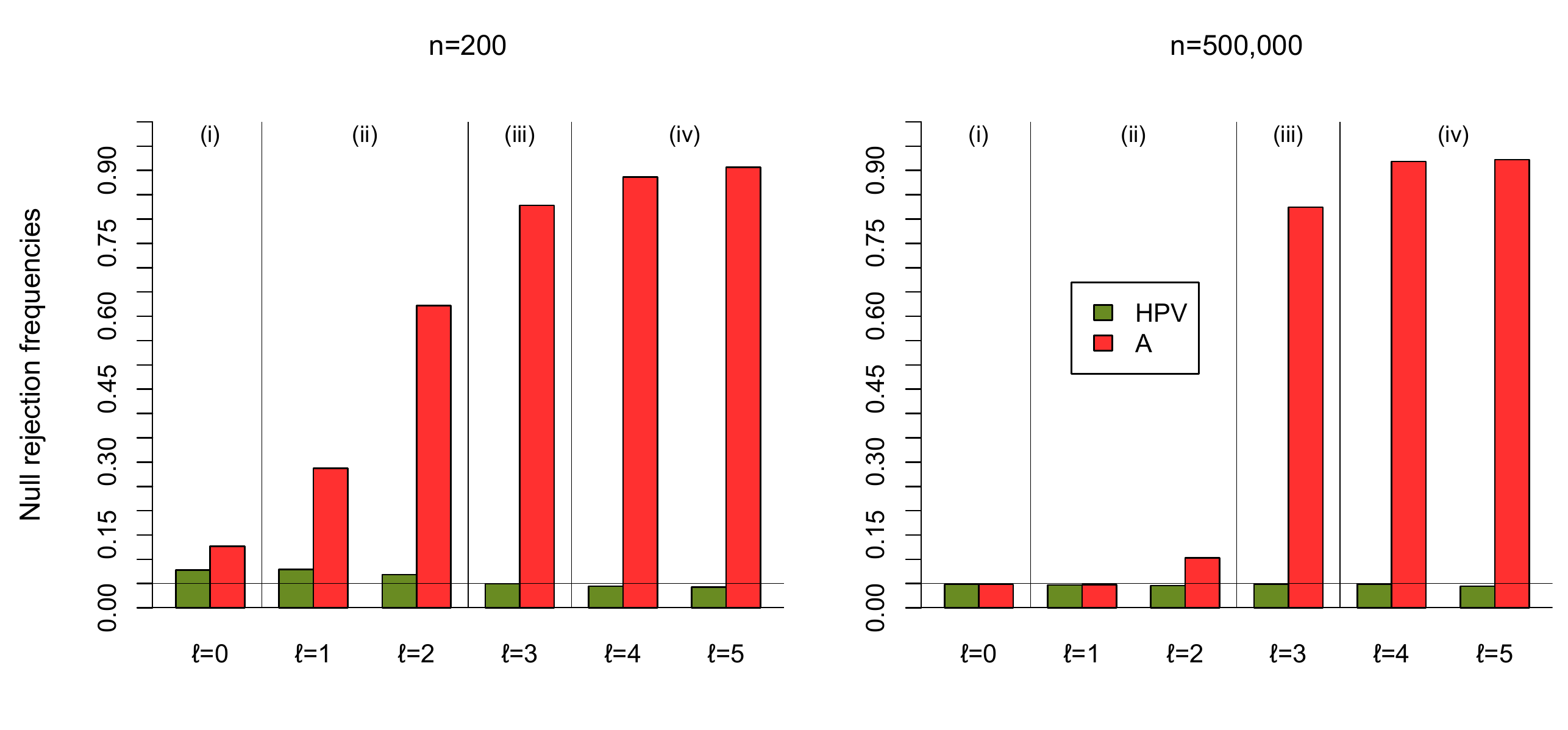}
\vspace{-5mm}
\caption{Empirical rejection frequencies, under the null hypothesis, of the tests~$\phi_{\rm HPV}$ and~$\phi_{\rm A}$ performed at nominal level~$5\%$. Results are based on $M=10,\!000$ independent ten-dimensional Gaussian random samples of size~$n=200$ (left) and size~$n=500,\!000$ (right). Increasing values of~$\ell$ bring the underlying spiked covariance matrix closer and closer to a multiple of the identity matrix; see Section~\ref{sec:Null} for details. The link between the values of~$\ell$ and the asymptotic regimes (i)--(iv) from Section~\ref{sec:preliminaries} is provided in each barplot.}
\label{Fig1}
\end{center}
\end{figure}

In this section, we will study the null asymptotic behaviors of~$\phi_{\rm A}$ and~$\phi_{\rm HPV}$ under weak identifiability, that is, we do so under the sequences of (null) hypotheses~${\rm P}_{\thetab_1^0,r_n,v}$ introduced in the previous section. Before doing so theoretically, we consider the following Monte Carlo exercise. For any~$\ell=0,1,\ldots,5$, we 
\vspace{-.7mm}
 generated $M=10,\!000$ mutually independent random samples~$\Xb_i^{(\ell)}$, $i=1, \ldots, n$, from the $(p=10)$-variate normal distribution with mean zero and covariance matrix
\begin{equation}
\Sigb_n^{(\ell)}:= {\bf I}_p+n^{-\ell/6} \thetab_1^0 \thetab_1^{0\prime}
,
\end{equation}
where~$\thetab_1^0$ is the first vector of the canonical basis of~$\R^p$. 
\vspace{-.5mm}
For each sample, we performed the tests $\phi_{\rm HPV}$ and~$\phi_{\rm A}$ for ${\cal H}_0\n: \thetab_1= \thetab_1^0$ at nominal level~$5\%$. 
The value of~$\ell$ allows to consider the various regimes above, namely
(i) away from contiguity ($\ell=0$), 
(ii) beyond contiguity ($\ell=1,2$),
(iii) under contiguity ($\ell=3$),
and 
(iv) under strict contiguity ($\ell=4,5$). Increasing values of~$\ell$ therefore provide harder and harder inference problems. Figure~\ref{Fig1}, that reports the resulting rejection frequencies for~$n=200$ and~$n=500,\!000$, suggests that~$\phi_{\rm HPV}$ asymptotically shows the target Type~1 risk in all regimes, hence is \emph{validity-robust} to weak identifiability. In sharp contrast, $\phi_{\rm A}$ seems to exhibit the right asymptotic null size in regimes~(i)--(ii) only, as it dramatically overrejects the null hypothesis (even asymptotically) in regimes~(iii)--(iv).

We now turn to the theoretical investigation of the asymptotic behaviors of~$\phi_{\rm A}$ and~$\phi_{\rm HPV}$ under weak identifiability. Obviously, this will heavily rely on Lemmas~\ref{Lemeigenvalues}--\ref{Lemeigenvectors}. For~$\phi_{\rm HPV}$, we have the following result.

\begin{Theor} 
\label{TheorHPVNull}
Fix a unit $p$-vector~$\thetab_1^0$, $v>0$ and a bounded positive real sequence~$(r_n)$. Then, under~${\rm P}_{\thetab_1^0,r_n,v}$, 
$$
Q_{\rm HPV} \stackrel{\mathcal{D}}{\to} \chi^2_{p-1}
,
$$
so that, in all regimes~(i)--(iv) from the previous section, the test~$\phi_{\rm HPV}$ has asymptotic size~$\alpha$ under the null hypothesis.
\end{Theor}

This result confirms that the test~$\phi_{\rm HPV}$ is validity-robust to weak identifiability in the sense that it will asymptotically meet the nominal level constraint in scenarios that are arbitrarily close to the spherical case. As hinted by the above Monte Carlo exercise, the situation is more complex for the Anderson test~$\phi_{\rm A}$. We have the following result. 

\begin{Theor} 
\label{TheorAndNull}
Fix a unit $p$-vector~$\thetab_1^0$, $v>0$ and a bounded positive real sequence~$(r_n)$. Let~$\Zb$ be a $p\times p$ random matrix such that~${\rm vec}({\bf Z})\sim {\cal N}({\bf 0}, {\bf I}_{p^2}+ {\bf K}_p)$. Then, we have the following 
under~${\rm P}_{\thetab_1^0,r_n,v}$:
\begin{itemize} 
\item[(i)--(ii)]  
if $r_n\equiv 1$ or if~$r_n$ is $o(1)$ with $\sqrt{n} r_n \to \infty$, then
$$
Q_{\rm A} \stackrel{\mathcal{D}}{\to} \chi^2_{p-1}
,
$$
so that the test~$\phi_{\rm A}$ has asymptotic size~$\alpha$ under the null hypothesis;
\vspace{1mm} 
\item[(iii)]  
if~$r_n=1/\sqrt{n}$, then
$$
Q_{\rm A} \stackrel{\mathcal{D}}{\to} 
\sum_{j=2}^p \, (\ell_1(v)-\ell_j(v))^2 (w_{j1}(v))^2
,
$$
where~$\ell_1(v)\geq \ldots\geq \ell_p(v)$ are the eigenvalues of~${\bf Z}+{\rm diag}(v,0,\ldots,0)$ and~$\wb_j(v)=(w_{j1}(v),\ldots,w_{jp}(v))'$ is an arbitrary unit eigenvector associated with~$\ell_j(v)$ 
$($with probability one, the only freedom in the choice of~$\wb_j(v)$ is in the sign of~$w_{j1}(v)$, that is clearly irrelevant here$)$;
\vspace{1mm} 
\item[(iv)] if~$r_n=o(1/\sqrt{n})$, then   
$$
Q_{\rm A} \stackrel{\mathcal{D}}{\to} 
\sum_{j=2}^p \, (\ell_1-\ell_j)^2 w_{j1}^2
,
$$
where~$\ell_1\geq \ldots\geq \ell_p$ are the eigenvalues of~${\bf Z}$ and~$\wb_j=(w_{j1},\ldots,w_{jp})'$ is an arbitrary unit eigenvector associated with~$\ell_j$.
\end{itemize}
\end{Theor}

\begin{figure}[htbp!]  
\begin{center} 
\includegraphics[width=\textwidth]{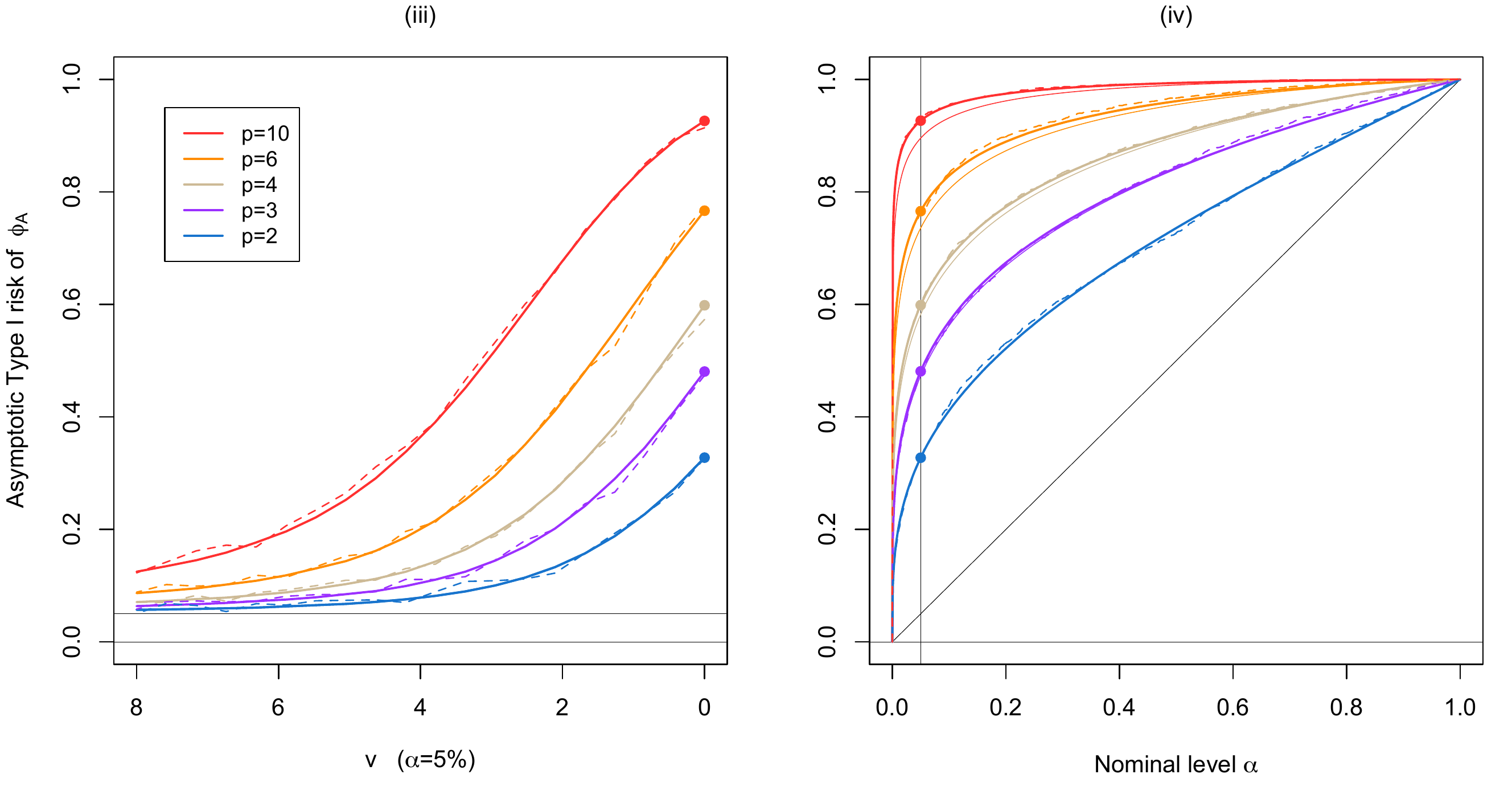}
\vspace{-5mm}
\caption{(Left:) Plots, for various dimensions~$p$, of the approximate asymptotic Type~I risk~$r_{p,.95}^{\rm (iii)}(v)$ (see~(\ref{approxtype1a})) 
of the $5\%$-level Anderson test for~$\mathcal{H}_0:\thetab_1=\thetab_1^0$ under~${\rm P}_{\thetab_1^0,r_n,v}$, with~$r_n=1/\sqrt{n}$ (regime~(iii)). The dashed curves report the corresponding rejection frequencies for sample size~$n=10,\!000$ (a regular grid of~20 $v$-values in~$[0,8]$ was considered and rejection frequencies were computed from $2,\!500$ independent replications in each case). 
(Right:) Plots, for the same dimensions~$p$, of the approximate asymptotic Type~I risk~$r_{p,\alpha}^{\rm (iv)}$ (see~(\ref{approxtype1b}))
of the level-$\alpha$ Anderson test for~$\mathcal{H}_0:\thetab_1=\thetab_1^0$ under~${\rm P}_{\thetab_1^0,r_n,v}$, with~$r_n=o(1/\sqrt{n})$ (regime~(iv)). The dashed curves report the corresponding rejection frequencies computed from $2,\!500$ independent standard normal samples of size~$n=10,\!000$.  The thin curves represent what the asymptotic Type~I risk of the level-$\alpha$ Anderson test would be if the null asymptotic distribution
of $Q_A$ in regime~(iv) would be~$4\chi^2_{p-1}$; see the discussion below Corollary~\ref{CorolAnd}.}  
\label{Fig2}
\end{center} 
\end{figure}

This result ensures that the Anderson test asymptotically meets the nominal level constraint in regimes~(i)--(ii). To see whether or not this extends to regimes~(iii)--(iv), we need to investigate the asymptotic distributions in Theorem~\ref{TheorAndNull}(iii)--(iv). We consider first the asymptotic distribution of $Q_{\rm A}$ under~${\rm P}_{\thetab_1^0, 1/\sqrt{n}, v}$, that is,
\vspace{-.2mm}
  in the contiguity regime (iii). To do so, we generated, for various dimensions~$p$ and for each $v=8(\ell-1)/19$, $\ell=1, \ldots, 20$, in a regular grid of 20~$v$-values in~$[0,8]$, a collection of $M=10^6$ independent values $Z_1(v), \ldots, Z_M(v)$ from the asymptotic distribution in Theorem~\ref{TheorAndNull}(iii). For each~$p$ and~$v$, we then recorded 
\begin{equation}
	\label{approxtype1a}
r_{p,.95}^{\rm (iii)}(v)
:=
\frac{1}{M} \sum_{m=1}^M \mathbb{I}\big[Z_m(v)>\chi^2_{p-1,.95}\big] 
,
\end{equation}
which is an excellent approximation of the asymptotic null size of the $5\%$-level Anderson test under~${\rm P}_{\thetab_1^0,1/\sqrt{n},v}$ (a $99\%$-confidence interval for the true asymptotic size has a length smaller than~$.0026$). 
\vspace{-.5mm}
 The left panel of Figure~\ref{Fig2} plots~$r_{p,.95}^{\rm (iii)}(v)$ as a function of~$v$ for several dimensions~$p$. Clearly, the Anderson test is, irrespective of~$p$ and~$v$, asymptotically overrejecting the null hypothesis. The asymptotic Type I risk increases with~$p$ and decreases with~$v$ (letting~$v$ go to infinity essentially provides regime~(ii), which explains that the Type I risk then converges to the nominal level). In the right panel of Figure~\ref{Fig2}, we generated, still for various values of~$p$, $M=10^6$ independent values~$Z_1,\ldots,Z_M$ from the asymptotic distribution in Theorem~\ref{TheorAndNull}(iv) and plotted the function mapping~$\alpha$ to
\begin{equation}
	\label{approxtype1b}
r_{p,\alpha}^{\rm (iv)}
:=
\frac{1}{M} \sum_{m=1}^M \mathbb{I}\big[Z_m>\chi^2_{p-1,1-\alpha}\big]
,
\end{equation}
which accurately approximates the asymptotic Type I risk of the level-$\alpha$ Anderson test in regime~(iv). Irrespective of~$\alpha$, the Anderson test is still overrejecting the null hypothesis asymptotically and the asymptotic Type I risk increases with~$p$. Overrejection is dramatic: for instance, in dimension~10, the asymptotic Type I risk of the $5\%$-level Anderson test exceeds~$92\%$. Empirical rejection frequencies of the Anderson test, which are also showed in Figure~\ref{Fig2}, clearly support the asymptotic results in Theorem~\ref{TheorAndNull}(iii)--(iv).

The asymptotic distributions in Theorem~\ref{TheorAndNull}(iii)--(iv) are explicitly described yet are quite complicated. Remarkably, for regime~(iv), the asymptotic distribution is a classical one in the bivariate case~$p=2$. More precisely, we have the following result. 

\begin{Corol} 
\label{CorolAnd}
Fix~$p=2$, a unit $p$-vector~$\thetab_1^0$, $v>0$ and a positive real sequence~$(r_n)$ such that~$\sqrt{n}r_n\to 0$. Then, under~${\rm P}_{\thetab_1^0,r_n,v}$, 
$$
Q_{\rm A} \stackrel{\mathcal{D}}{\to} 4\chi^2_{1}
,
$$
so that, irrespective of~$\alpha\in(0,1)$, the asymptotic size of~$\phi_{\rm A}$ under the null hypothesis is strictly larger than~$\alpha$.
\end{Corol}

This result shows in a striking way the impact weak identifiability may have, in the bivariate case, on the null asymptotic distribution of the Anderson test statistic~$Q_{\rm A}$: away from contiguity and beyond contiguity (regimes (i)--(ii)), $Q_{\rm A}$ is asymptotically~$\chi^2_{1}$ under the null hypothesis (Theorem~\ref{TheorAndNull}), whereas under strict contiguity (regime~(iv)), this statistic is asymptotically~$4\chi^2_1$ under the null hypothesis. The result also allows us to quantify, for any nominal level~$\alpha$, how much the bivariate Anderson test will asymptotically overreject the null hypothesis in regime~(iv). More precisely, the Type~1 risk, in this regime, of the \mbox{level-$\alpha$}  bivariate Anderson test converges to~${\rm P}[4Z^2 > \chi^2_{1,1-\alpha}]$, where~$Z$ is standard normal. For~$\alpha=0.1\%$, $1\%$ and $5\%$, this provides in regime~(iv) an asymptotic Type~1 risk of about~$10\%$, $19.8\%$ and~$32.7\%$, respectively (which exceeds the nominal level by about a factor~100, 20 and 6.5, respectively!)
In dimensions~$p\geq 3$, the null asymptotic distribution of~$Q_{\rm A}$ in regime~(iv), as showed in the right panel of Figure~\ref{Fig2}, is very close to~$4\chi^2_{p-1}$, particularly so for~$p=3$. Yet the distribution is not~$4\chi^2_{p-1}$. For instance, in dimension~$p=3$, it can be showed that the null asymptotic  distribution of~$Q_{\rm A}$ in regime~(iv) has mean~$49/6$, whereas the distribution~$4\chi^2_{2}$ has mean~$8=48/6$ (also, computing the variance of the null asymptotic  distribution of~$Q_{\rm A}$ shows that this distribution is not of the form~$\lambda \chi^2_{2}$ for any~$\lambda>0$).  
 
\begin{figure}[htbp!] 
\begin{center} 
\includegraphics[width=\textwidth]{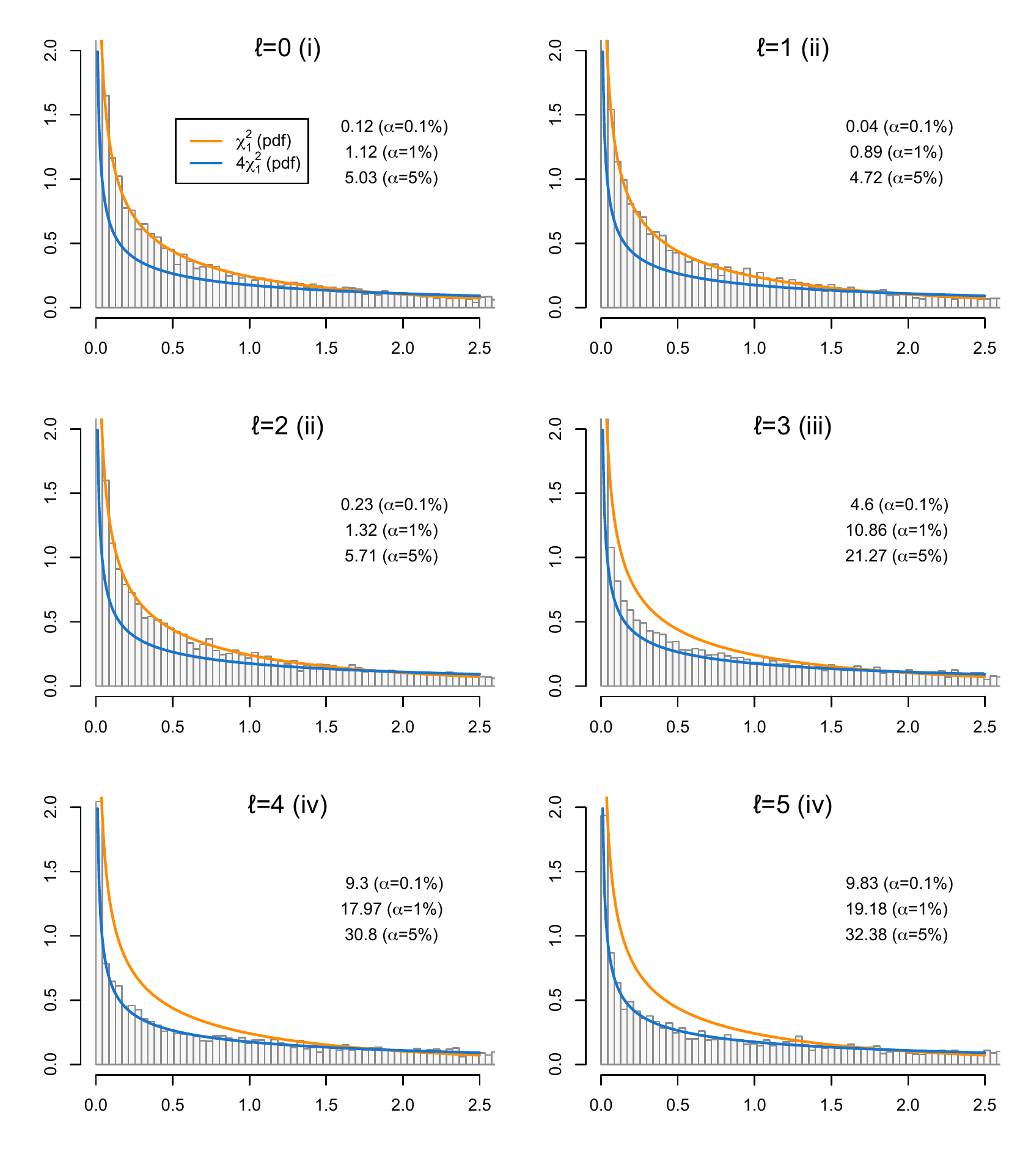}
\vspace{-5mm}
\caption{For each~$\ell=0,1,\ldots,5$, histograms of values of~$Q_{\rm A}$ from $M=10,\!000$ independent (null)  Gaussian random samples of size~$n=500,\!000$ and dimension~$p=2$. Increasing values of~$\ell$ bring the underlying spiked covariance matrix closer and closer to a multiple of the identity matrix; see Section~\ref{sec:Null} for details. The links between the values of~$\ell$ and the asymptotic regimes~(i)--(iv) from Section~\ref{sec:preliminaries} are provided in each case. For any value of~$\ell$, the density of the null asymptotic distribution of~$Q_{\rm A}$ in regimes~(i)--(ii) (resp., in regime~(iv)) is plotted in orange (resp., in blue) and the empirical Type~1 risk of the \mbox{level-$\alpha$} Anderson test is provided for~$\alpha=0.1\%$, $1\%$ and~$5\%$.} 
\label{Fig3}
\end{center}
\end{figure}
   
We close this section with a last simulation illustrating Theorem~\ref{TheorAndNull} and its consequences. To do so, we generated, for any~$\ell=0,1,\ldots,5$, a collection of~$M=10,\!000$ mutually independent random samples of size~$n=500,\!000$ from the bivariate normal distribution with mean zero and covariance matrix~$\Sigb_n^{(\ell)}:= {\bf I}_2+n^{-\ell/6} \thetab_1^0 \thetab_1^{0\prime}$, with~$\thetab_1^0=(1,0)'$. This is therefore essentially the bivariate version of the ten-dimensional numerical exercise leading to Figure~\ref{Fig1}. For each value of~$\ell$, Figure~\ref{Fig3} provides histograms of the resulting~$M$ values of the Anderson test statistic~$Q_{\rm A}$, along with plots of the densities of the~$\chi^2_1$ and~$4\chi^2_1$ distributions, that is, densities of the null asymptotic distribution of~$Q_{\rm A}$ in regimes~(i)--(ii) and in regime~(iv), respectively. In these three regimes, the histograms are perfectly fitted by the corresponding density. The figure also provides the empirical Type~1 risks of the \mbox{level-$\alpha$} Anderson test for~$\alpha=0.1\%$, $1\%$ and~$5\%$. Clearly, these Type~1 risks are close to the theoretical asymptotic ones both in regimes~(i)--(ii) (namely, $\alpha$) and in regime~(iv) (namely, the Type~1 risks provided in the previous paragraph). 
\vspace{4mm}

 
 %
 %
 %
 %
 %
 %

\section{Non-null and optimality results}
\label{sec:LAN}

The previous section shows that, unlike~$\phi_{\rm A}$, the test~$\phi_{\rm HPV}$ is validity-robust to weak identifiability. However, the trivial \mbox{level-$\alpha$} test, that randomly rejects the null hypothesis with probability~$\alpha$, of course enjoys the same robustness property. This motivates investigating whether or not the validity-robustness of~$\phi_{\rm HPV}$ is obtained at the expense of power. In this section, we therefore study the asymptotic non-null behavior of~$\phi_{\rm HPV}$ and show that this test actually still enjoys strong optimality properties under weak identifiability.

Throughout, optimality is to be understood in the Le Cam sense. In the present hypothesis testing context, Le Cam optimality requires studying local log-likelihood ratios of the form 
\begin{equation}
	\label{deflikeratio}
\Lambda_n
:= 
\log \frac{d{\rm P}_{\thetab_1^0+ \nu_n \taub_n,r_n,v}}{d{\rm P}_{\thetab_1^0,r_n,v}}
,
\end{equation}
where the bounded sequence~$(\taub_n)$ in $\R^p$ and the positive real sequence~$(\nu_n)$ are such that, for any~$n$, $\thetab_1^0+ \nu_n \taub_n$ is a unit $p$-vector, hence, is an admissible perturbation of~$\thetab_1^0$. This imposes that~$(\taub_n)$ satisfies 
\begin{equation} 
\label{tau}
\thetab_1^{0\prime} \taub_n = -\frac{\nu_n}{2} \| \taub_n\|^2
\end{equation}
for any~$n$. The following result then describes, in each of the four regimes~(i)--(iv) considered in the previous sections, the asymptotic behavior of the log-likelihood ratio~$\Lambda_n$.

\begin{Theor}
\label{TheorLAN}
Fix a unit $p$-vector~$\thetab_1^0$, $v>0$ and a bounded positive real sequence~$(r_n)$. 
Then, we have the following under~${\rm P}_{\thetab_1^0,r_n,v}$:
\begin{itemize}
\item[(i)] if $r_n\equiv 1$, then, with~$\nu_n=1/\sqrt{n}$,
$$
\Deltab_n
=
\frac{v}{1+v}\sqrt{n}({\bf I}_p-\thetab_1^0\thetab_1^{0\prime})({\bf S}_n-\Sigb_n) \thetab_1^0
\ \textrm{ and }\
\Gamb
=
\frac{v^2}{1+v} ({\bf I}_p- \thetab_1^0\thetab_1^{0\prime})
,
$$
we have that
$\Lambda_n=\taub_n\pr  \, \Deltab_n -\frac{1}{2}\, \taub_n\pr \Gamb \taub_n + o_{\rm P}(1)$ and that~$\Deltab_n$ is asymptotically normal with mean zero and covariance matrix~$\Gamb$;
\vspace{1mm}
\item[(ii)] if $r_n$ is $o(1)$ with $\sqrt{n}r_n \rightarrow \infty$, then, with~$\nu_n=1/(\sqrt{n}r_n)$,
$$
\Deltab_n
=
v\sqrt{n}({\bf I}_p-\thetab_1^0\thetab_1^{0\prime})({\bf S}_n-\Sigb_n) \thetab_1^0
\ \textrm{ and }\
\Gamb
=
v^2 ({\bf I}_p- \thetab_1^0\thetab_1^{0\prime})
,
$$
we similarly have that
$\Lambda_n=\taub_n\pr  \, \Deltab_n -\frac{1}{2}\, \taub_n\pr \Gamb \taub_n + o_{\rm P}(1)$ and that~$\Deltab_n$ is still asymptotically normal with mean zero and covariance matrix~$\Gamb$;
\vspace{1mm}
\item[(iii)] if $r_n=1/\sqrt{n}$, then, letting~$\nu_n \equiv 1$,
\begin{eqnarray} 
\label{LAQvic}
\lefteqn{
\Lambda_n
=
 \taub_n\pr  \Big[ v\sqrt{n} ({\bf S}_n-\Sigb_n) \Big(\thetab_1^0 + \frac{1}{2} \taub_n\Big)\Big]
 }
 \\[2mm]
 & & 
 \hspace{28mm} 
 -\frac{v^2}{2} \| \taub_n \|^2+ \frac{v^2}{8} \|\taub_n\|^4 +o_{\rm P}(1),
\nonumber
\end{eqnarray}
where, if~$(\taub_n)\to \taub$, then $\taub_n\pr \sqrt{n} ({\bf S}_n-\Sigb_n) (\thetab_1^0 + \frac{1}{2} \taub_n)$ is asymptotically normal with mean zero and variance~$\| \taub \|^2- \frac{1}{4} \| \taub \|^4$;
\vspace{1mm}
\item[(iv)] if $r_n=o(1/\sqrt{n})$, then, even with~$\nu_n \equiv 1$, we have that $\Lambda_n$ is $o_{\rm P}(1)$.
\end{itemize}
\end{Theor}

This result shows that, for any fixed~$v>0$ and for any fixed sequence~$(r_n)$ associated with regime~(i) (away from contiguity) or regime~(ii) (beyond contiguity), the sequence of models~$\{{\rm P}_{\thetab_1,r_n,v}:\thetab_1\in\mathcal{S}^{p-1}\}$ is LAN (locally asymptotically normal), with central sequence
\begin{equation}
\label{centralvv}	
\Deltab_{n,\delta}
:= 
\frac{\sqrt{n}v}{1+\delta v} 
\big({\bf I}_p-\thetab_1^0\thetab_1^{0\prime}\big)
({\bf S}_n-\Sigb_n) \thetab_1^0
\end{equation}
and Fisher information matrix
\begin{equation}
\label{infovv}	
\Gamb_{\delta}
:=
\frac{v^2}{1+\delta v}  
\big({\bf I}_p- \thetab_1^0\thetab_1^{0\prime}\big)
,
\end{equation}
where~$\delta:=1$ if regime~(i) is considered and~$\delta:=0$ otherwise. Denoting as~${\bf A}^{-}$ the Moore-Penrose inverse
\vspace{-.5mm}
  of~${\bf A}$, it follows that the locally asymptotically maximin test for~${\cal H}_{0}\n: \thetab_1= \thetab_1^0$ against~${\cal H}_{1}\n: \thetab_1 \neq \thetab_1^0$ rejects the null hypothesis at asymptotic level~$\alpha$ when
\begin{equation} 
\label{optbeyond}
Q_\delta
=
\Deltab'_{n,\delta} \Gamb_{\delta}^{-}\Deltab_{n,\delta}
=
\frac{n}{1+\delta v} 
\,
\thetab_1^{0\prime}
 {\bf S}_n 
\big({\bf I}_p-\thetab_1^0\thetab_1^{0\prime}\big) 
{\bf S}_n\thetab_1^0 
>
 \chi^2_{p-1, 1-\alpha}
.
\end{equation}
In view of~(\ref{asympeqproofpre}) in the proof of Theorem~\ref{TheorHPVNull}, we have that~$Q_{\rm HPV}= Q_\delta+o_{\rm P}(1)$ under~${\rm P}_{\thetab_1^0,r_n,v}$, for any~$v>0$ and any bounded positive sequence~$(r_n)$, hence also, from contiguity, under local alternatives of the form~${\rm P}_{\thetab_1^0+ \taub_n/(\sqrt{n} r_n),r_n,v}$. We may therefore conclude that, away from contiguity and beyond contiguity, the test~$\phi_{\rm HPV}$ is Le Cam optimal for the problem at hand. We have the following result. 
 
 \begin{Theor}
 \label{HPVnonull12}
 	Fix a unit $p$-vector~$\thetab_1^0$, $v>0$ and a positive real sequence~$(r_n)$ 
\vspace{-.3mm}
satisfying~(i) $r_n\equiv 1$ 
or (ii) $r_n=o(1)$ with~$\sqrt{n}r_n\to \infty$. Then, the test~$\phi_{\rm HPV}$ is locally asymptotically
\vspace{-.8mm}
 maximin at level~$\alpha$ when testing~${\cal H}_{0}\n: \thetab_1= \thetab_1^0$ against~${\cal H}_{1}\n: \thetab_1 \neq \thetab_1^0$. Moreover, under~${\rm P}_{\thetab_1^0+ \taub_n/(\sqrt{n} r_n),r_n,v}$, with~$(\taub_n)
 \linebreak
 \to \taub$, the statistic~$Q_{\rm HPV}$ is asymptotically non-central chi-square with~\mbox{$p-1$} degrees of freedom and with non-centrality parameter~$(v^2/(1+\delta v))\|\taub\|^2$. 
 \end{Theor}

Under strict contiguity (Theorem~\ref{TheorLAN}(iv)), no asymptotic \mbox{level-$\alpha$} test can show non-trivial asymptotic powers against the most severe alternatives of the form $\thetab_1^0+ \taub$. Therefore, the test~$\phi_{\rm HPV}$ is also optimal in regime~(iv), even though this optimality is degenerate since the trivial \mbox{level-$\alpha$} test is also optimal in this regime. We then turn to Theorem~\ref{TheorLAN}(iii), where the situation is much less standard, as the sequence of experiments~$\{{\rm P}_{\thetab_1,1/\sqrt{n},v}:\thetab_1\in\mathcal{S}^{p-1}\}$ there is neither LAN, nor LAMN (locally asymptotically mixed normal), nor LAQ (locally asymptotically quadratic); see \cite{Je95} and \cite{RO11}. While, to the best of our knowledge, the form of optimal tests in such non-standard limiting experiments remains unknown, we will still be able below to draw conclusions about optimality for small values of $\taub$. Before doing so, note that, by using the Le Cam first lemma (see, e.g., Lemma~6.4 in \citealp{van1998}), Theorem~\ref{TheorLAN}(iii) readily entails that, for any~$v>0$, the sequences of hypotheses~${\rm P}_{\thetab_1^0,1/\sqrt{n},v}$ and~${\rm P}_{\thetab_1^0+\taub_n,1/\sqrt{n},v}$, with~$(\taub_n)\to\taub$, are mutually contiguous. Consequently, the asymptotic non-null distribution of the test statistic~$Q_{\rm HPV}$ under contiguous alternatives may still be obtained from the Le Cam third lemma. We have the following result.
 
 \begin{Theor}
 \label{HPVnonull34}
 	Fix a unit $p$-vector~$\thetab_1^0$, $v>0$ and a positive real sequence~$(r_n)$ satisfying~(iii) $r_n=1/\sqrt{n}$ or (iv) $r_n=o(1/\sqrt{n})$. Then, in case~(iii), the test statistic~$Q_{\rm HPV}$, under~${\rm P}_{\thetab_1^0+ \taub_n,r_n,v}$, with~$(\taub_n)\to \taub$, is asymptotically non-central chi-square with~$p-1$ degrees of freedom and with non-centrality parameter
\begin{equation} 
\label{ncpHPV}
\frac{v^2}{16}
\| \taub \|^2
\big(4- \| \taub \|^2\big)
\big( 2-\| \taub\|^2\big)^2
,
\end{equation}
so that the test~$\phi_{\rm HPV}$ is rate-consistent.
\vspace{-.5mm}
  In case~(iv), this test is locally asymptotically maximin at level~$\alpha$ when testing~${\cal H}_{0}\n: \thetab_1= \thetab_1^0$ against~${\cal H}_{1}\n: \thetab_1 \neq \thetab_1^0$, but trivially so since its asymptotic power against any sequence of alternatives~${\rm P}_{\thetab_1^0+ \taub_n,r_n,v}$, with~$(\taub_n)\to \taub$, is then equal to~$\alpha$. 
 \end{Theor}

\begin{figure}[htbp!] 
\begin{center} 
\includegraphics[width=\textwidth]{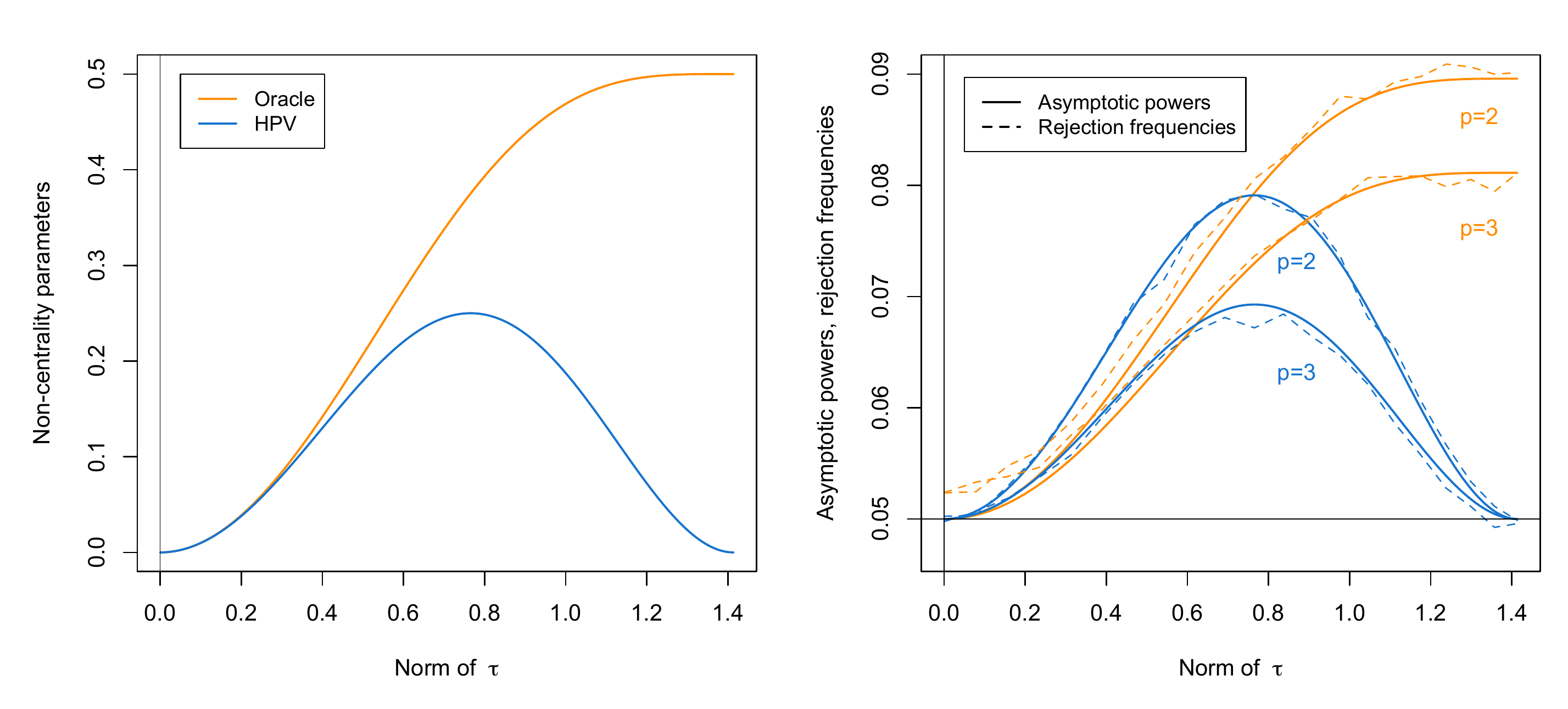}
\vspace{-5mm}
\caption{(Left:) Non-centrality parameters~(\ref{ncpHPV}) and~(\ref{ncporacle}), as a function of~$\|\taub\|(\in[0,\sqrt{2}])$, in the asymptotic non-central chi-square distributions of the test statistics of~$\phi_{\rm HPV}$ and~$\phi_{\rm oracle}$, respectively, under alternatives of the form~${\rm P}_{\thetab_1^0+ \taub,1/\sqrt{n},1}$. (Right:) The corresponding asymptotic power curves in dimensions~$p=2$ and~$p=3$, as well as the empirical power curves resulting from the Monte Carlo exercise described at the end of Section~\ref{sec:LAN}.}
\label{Fig4}
\end{center} 
\end{figure}

Figure~\ref{Fig4} plots the non-centrality parameter in~(\ref{ncpHPV}) as a function of~$\taub\in[0,\sqrt{2}]$ (since~$\thetab_1$ is defined up to a sign, one may restrict to alternatives~$\thetab_1^0+\taub$ in the hemisphere centered at~$\thetab_1^0$), as well as the resulting asymptotic powers of the test~$\phi_{\rm HPV}$ in dimensions~$p=2,3$. This test shows no asymptotic power when~$\|\taub\|=\sqrt{2}$ (that is, when~$\thetab_1$ is orthogonal to~$\thetab_1^0$), hence clearly does not enjoy \emph{global-in-$\taub$} optimality properties in regime~(iii). 
As we now explain, however,~$\phi_{\rm HPV}$ exhibits excellent \mbox{\emph{local-in-$\taub$}} optimality properties in this regime. In order to see this, note that decomposing~${\bf I}_p$ into~$\big({\bf I}_p-\thetab_1^0\thetab_1^{0\prime}\big)+\thetab_1^0\thetab_1^{0\prime}$ and using repeatedly~(\ref{tau}) allows to rewrite~(\ref{LAQvic}) as
\begin{eqnarray}
\lefteqn{
\hspace{3mm} 
\Lambda_n
=
 \taub_n\pr \Deltab_{n,0} 
 - \frac{1}{2} \taub_n'\Gamb_0\taub_n
}
\label{LAQvic2}
\\[2mm]
 & & 
 \hspace{3mm} 
 -\frac{\sqrt{n}v}{2}  \|\taub_n\|^2\thetab_1^{0\prime} ({\bf S}_n-\Sigb_n) \thetab_1^0
+
 \frac{\sqrt{n}v}{2} \taub_n\pr ({\bf S}_n-\Sigb_n) \taub_n
+o_{\rm P}(1)
,
\nonumber
\end{eqnarray}
where~$\Deltab_{n,0}$ and~$\Gamb_0$ were defined in~(\ref{centralvv}) and~(\ref{infovv}), respectively. For small perturbations~$\taub_n$,  the righthand side of~(\ref{LAQvic2}), after neglecting the second-order random terms in~$\taub_n$, becomes 
$$
 \taub_n\pr \Deltab_{n,0} 
- \frac{1}{2} \taub_n'\Gamb_0\taub_n
+o_{\rm P}(1)
,
$$
so that the sequence of experiments is then LAN again, with central sequence~$\Deltab_{n,0}$ and Fisher information matrix~$\Gamb_0$. This implies that the test in~(\ref{optbeyond}) and (in view of the asymptotic equivalence stated below~(\ref{optbeyond})) the test~$\phi_{\rm HPV}$ are locally(-in-$\taub$) asymptotically maximin.  

Now, if the objective is to construct a test that will perform well also for large perturbations~$\taub_n$ in regime~(iii), it may be tempting to consider as a test statistic the linear-in-$\taub$ part of the random term in~(\ref{LAQvic}), namely~$\tilde{\Deltab}_n:=v\sqrt{n} ({\bf S}_n-\Sigb_n) \thetab_1^0$. Since~$\tilde{\Deltab}_n$, under~${\rm P}_{\thetab_1^0,1/\sqrt{n},v}$, is 
\vspace{-.5mm}
 asymptotically normal with mean zero and covariance~$\tilde{\Gamb}:=v^2 ({\bf I}_p+ \thetab_1^0\thetab_1^{0\prime})$, the resulting test,~$\phi_{\rm oracle}$ say, rejects the null hypothesis at asymptotic level $\alpha$ when
\begin{eqnarray}
\lefteqn{
\hspace{10mm} 
\tilde{Q}  
:=
\tilde{\Deltab}_n' 
\tilde{\Gamb}^{-1}
\tilde{\Deltab}_n
=
n \thetab_1^{0\prime} ({\bf S}_n- \Sigb_n) 
\big({\bf I}_p+\thetab_1^0\thetab_1^{0\prime}\big)^{-1}
({\bf S}_n- \Sigb_n)\thetab_1^0 
}
\label{oracletest}
\\[2mm]
& & 
\hspace{-6mm} 
=
n 
\big({\bf S}_n\thetab_1^0 - \big({\textstyle{1+\frac{v}{\sqrt{n}}}}\big)\thetab_1^0 \big)'
\big({\bf I}_p- {\textstyle{\frac{1}{2}}} \thetab_1^0\thetab_1^{0\prime}\big)
\big({\bf S}_n\thetab_1^0 - \big({\textstyle{1+\frac{v}{\sqrt{n}}}}\big)\thetab_1^0 \big)
>
\chi^2_{p,1-\alpha}
;
\nonumber
\end{eqnarray}
the terminology ``oracle" stresses that this test requires knowing the true value of~$v$. The Le Cam third lemma 
\vspace{-.3mm}
entails that~$\tilde{\Deltab}_n$, under~${\rm P}_{\thetab_1^0+ \taub_n,1/\sqrt{n},v}$, with~$(\taub_n)\to \taub$, is asymptotically normal with mean $v^2(1-\frac{1}{2} \| \taub\|^2) \taub-(v^2/2) \| \taub\|^2 \thetab_1^0$ and covariance matrix~$\tilde{\Gamb}$. Therefore, under the same sequence of hypotheses, $\tilde{Q}$ is asymptotically non-central chi-square with~$p$ degrees of freedom and with non-centrality parameter 
\begin{equation} 
\label{ncporacle}
\frac{v^2}{16}
\| \taub \|^2
\big(4- \| \taub \|^2\big) 
\Big( 4-2\| \taub\|^2+\frac{1}{2}\| \taub\|^4\Big)
.
\end{equation}
Note that the difference between this non-centrality parameter and the one in~(\ref{ncpHPV}) is~$O(\|\taub\|^4)$ as~$\|\taub\|$ goes to zero. Since~$\phi_{\rm HPV}$ is based on a smaller number of degrees of freedom ($p-1$, versus~$p$ for the oracle test), it will therefore exhibit larger asymptotic powers than the oracle test for small values of~$\taub$, which reflects the aforementioned local-in-$\taub$ optimality of~$\phi_{\rm HPV}$. 

Figure~\ref{Fig4} also plots the non-centrality parameter in~(\ref{ncporacle}) as a function of~$\|\taub\|$, as well as the asymptotic powers of the oracle test in dimensions~$p=2,3$. As predicted above,~$\phi_{\rm HPV}$ dominates~$\phi_{\rm oracle}$ for small values of~$\|\taub\|$, that is, for small perturbations. The opposite happens for large values of the perturbation and it is seen that~$\phi_{\rm oracle}$ overall is quite efficient. It is important to recall, however, that this test cannot be used in practice since it requires knowing the value of~$v$.
The figure further reports the results of a Monte Carlo exercise we conducted to check correctness of the highly non-standard asymptotic results obtained 
\vspace{-.7mm}
 in the present regime~(iii). In this simulation, we generated $M=200,\!000$ mutually independent random samples~$\Xb_i^{(k)}$, $i=1, \ldots, n=10,\!000$, $k=0,1,\ldots,20$, of $p$-variate ($p=2,3$) Gaussian random vectors with mean zero and covariance matrix
\begin{equation}
\Sigb_n^{(k)}
:= 
{\bf I}_p+n^{-1/2} (\thetab_1^0+ \taub_{k}) (\thetab_1^0+ \taub_k)\pr
 \end{equation}
where $\thetab_1^0=(1,0,\ldots,0)'\in\R^p$ and $\thetab_1^0+\taub_{k}=(\cos(k \pi/40), \sin(k \pi/40),0,\ldots,0)\pr
\linebreak
\in\R^p$. In each replication, we performed the tests $\phi_{\rm HPV}$ and $\phi_{\rm oracle}$ for ${\cal H}_0\n: \thetab_1= \thetab_1^0$ at nominal level~$5\%$. The value~$k=0$ is associated with the null hypothesis, whereas the values~$k=1, \ldots,20$ provide increasingly severe alternatives, the most severe of which involves a first eigenvector that is orthogonal to~$\thetab_1^0$. The resulting rejection frequencies are plotted in the right panel of Figure~\ref{Fig4}. Clearly, they are in perfect agreement with the asymptotic powers, which supports our theoretical results.
 
We conclude this section by stressing that, as announced in the introduction, the contiguity rate in Theorem~\ref{TheorLAN} depends on the regime considered. Clearly, the weaker the identifiability (that is, the closer the underlying distribution to the spherical Gaussian one), the slower the contiguity rate, that is, the hardest the inference problem on~$\thetab_1$.


\section{Extension to the elliptical case}
\label{sec:ellipt} 

Since we focused so far on multinormal distributions, a natural question is whether or not our results extend away from the Gaussian case. In this section, we discuss this in the framework of the most classical extension of multinormal distributions, namely in the class of elliptical distributions. More specifically, we will consider triangular arrays of $p$-variate observations~$\Xb_{ni}$, $i=1,\ldots,n$,  $n=1,2,\ldots$, where~$\Xb_{n1},\ldots,\Xb_{nn}$ form a random sample from the $p$-variate elliptical distribution with location~${\pmb\mu}_n$, covariance matrix~$\Sigb_n=\sigma_n({\bf I}_p+ r_n v\, \thetab_1\thetab_1\pr)$ (as in \eqref{localspiked}) and \emph{radial density}~$f$. That is, we assume that~$\Xb_{ni}$ admits the probability density function (with respect to the Lebesgue measure on~$\R^p$)
\begin{equation} 
\label{ellidens}
\xb 
\mapsto 
\frac{c_{p,f}}{(\det \Sigb_n)^{1/2}}
\, 
f\Big( \sqrt{(\xb-{\pmb\mu}_n)\pr \Sigb_n^{-1} (\xb-{\pmb\mu}_n)} \, \Big)
,
\end{equation}
where~$c_{p,f}>0$ is a normalization factor and where the radial density~$f:\R^+\to\R^+$ is such that the covariance matrix of~$\Xb_{ni}$ exists and is equal to~$\Sigb_n$;~$f$ is not a genuine density (as it does not integrate to one), but it determines the density of the Mahalanobis distance
$$
d_{ni}
:=
\sqrt{({\Xb}_{ni}-{\pmb\mu}_n)\pr \Sigb_n^{-1} ({\Xb}_{ni}-{\pmb\mu}_n)}
,
$$
which is given by~$r\mapsto (\mu_{p-1,f})^{-1} r^{p-1}f(r)\mathbb{I}[r\geq 0]$, with~$\mu_{\ell,f}:=\int_0^\infty r^{\ell}f(r)\,dr$.
In this section, we will assume that~$\Xb_{ni}$, or equivalently~$d_{ni}$, has finite fourth-order moments, that is, we will assume that~$f$ belongs to the collection~$\mathcal{F}$ of radial densities~$f$ above that further satisfy~$\mu_{p+3,f}<\infty$. This guarantees finiteness of the elliptical kurtosis coefficient
\begin{equation} 
\label{truekur}
\kappa_{p}(f)
:= 
\frac{p {\rm E}[ d_{ni}^4 ]}{(p+2) ({\rm E}[d_{ni}^2])^2}-1
\quad
\bigg(
\!
= \frac{p \mu_{p-1,f}\mu_{p+3,f}}{(p+2) \mu_{p+1,f}^2}-1
\bigg)
;
\end{equation} 
see, e.g., page~54 of \cite{And2003}. 
Classical radial densities in~$\mathcal{F}$ include the Gaussian one~$\phi(r)=\exp(-r^2/2)$ or the Student $t_\nu$ one~$f_\nu(r)=(1+r^2/(\nu-2))^{-(p+\nu)/2}$, with~$\nu>4$. The sequence of hypotheses associated with the triangular arrays of observations above will be denoted as~${\rm P}_{{\pmb \mu}_n,\sigma_n,\thetab_1,r_n,v,f}$.  

When it comes to testing the null hypothesis~$\mathcal{H}_0^{(n)}:\thetab_1=\thetab_1^0$, it is well-known that, even in the standard regime (i) ($r_n\equiv 1$), the Anderson test statistic~$Q_{\rm A}$ in~(\ref{Andstat}) and the HPV test statistic~$Q_{\rm HPV}$ in~(\ref{HPV2010}) are asymptotically $\chi^2_{p-1}$ under the sequence of null hypotheses~${\rm P}_{{\pmb \mu}_n,\sigma_n,\thetab_1^0,r_n,v,f}$ if and only if~$\kappa_{p}(f)$ takes the same value~$\kappa_p(\phi)=0$ as in the Gaussian case;
see, e.g., \cite{HPV10}. Consequently, there is no guarantee, even in regime~(i), that the corresponding tests~$\phi_{\rm A}$ and~$\phi_{\rm HPV}$ meet the asymptotic nominal level constraint under ellipticity, and it therefore makes little sense,  in the elliptical case, to investigate the robustness of these tests to weak identifiability. This explains why we will rather focus on their robustified versions~$\phi_{\rm A}^\dagger$ and~$\phi_{\rm HPV}^\dagger$, that reject the null hypothesis at asymptotic level~$\alpha$ whenever 
\begin{equation}  
\label{pseudos}
Q_{{\rm A}}^{(n)\dagger}
:=
\frac{Q_{\rm A}\n}{1+ \hat{\kappa}_p\n}
> 
\chi^2_{p-1,1-\alpha}
\quad \textrm{and} \quad
Q_{{\rm HPV}}^{(n)\dagger}
:=
\frac{Q_{\rm HPV}\n}{1+ \hat{\kappa}_p\n}
> 
\chi^2_{p-1,1-\alpha}
,
\end{equation}
respectively, where 
\begin{eqnarray} 
\hat{\kappa}_p\n
&\!\!\!:=\!\!\!&
\frac{p\big\{\frac{1}{n} \sum_{i=1}^{n} ( ({\Xb}_{ni}-\bar{\Xb}_n)'\Sb_n^{-1}({\Xb}_{ni}-\bar{\Xb}_n) )^2\big\}}{(p+2) \big\{\frac{1}{n} \sum_{i=1}^{n}  ({\Xb}_{ni}-\bar{\Xb}_n)'\Sb_n^{-1}({\Xb}_{ni}-\bar{\Xb}_n)\big\}^2}-1
\nonumber
\\[2mm]
&\!\!\!=\!\!\!&
\frac{1}{np(p+2)} 
\sum_{i=1}^{n}  
\,
( ({\Xb}_{ni}-\bar{\Xb}_n)'\Sb_n^{-1}({\Xb}_{ni}-\bar{\Xb}_n) )^2
-1
\label{hatkur}
\end{eqnarray}
is the natural estimator of the kurtosis coefficient~$\kappa_{p}(f)$. 
In the standard regime~(i), \cite{Tyl81,Tyl83} showed that~$\phi_{\rm A}^\dagger$ has asymptotic size~$\alpha$ under~${\rm P}_{{\pmb \mu}_n,\sigma_n,\thetab_1^0,r_n,v,f}$, whereas
\vspace{-.9mm}
  \cite{HPV10} proved the same result for~$\phi_{\rm HPV}^\dagger$ and established the asymptotic equivalence of both tests in probability. In the Gaussian case, these tests are, still in regime~(i), asymptotically equivalent to their original versions~$\phi_{\rm A}$ and~$\phi_{\rm HPV}$, hence inherit the optimality properties of the latter. The tests~$\phi_{\rm A}^\dagger$ and~$\phi_{\rm HPV}^\dagger$ may therefore be considered \emph{pseudo-Gaussian} versions of their antecedents, since they extend their validity to the class of elliptical distributions with finite fourth-order moments without sacrificing optimality in the Gaussian case.    

The above considerations make it natural to investigate the robustness of these pseudo-Gaussian tests to weak identifiability. Since these tests are invariant under translations and scale transformations, we will still assume, without loss of generality, that~${\pmb \mu}_n \equiv {\bf 0}$ and~$\sigma_n \equiv 1$ (see the discussion below Lemma \ref{LemLL}), and we will write accordingly~${\rm P}_{\thetab_1,r_n,v,f}:={\rm P}_{{\pmb 0},1,\thetab_1,r_n,v,f}$. Note that the Gaussian hypotheses~${\rm P}_{\thetab_1,r_n,v}={\rm P}_{\thetab_1,r_n,v,\phi}$ are those we considered in the previous sections of the paper. 
Our results will build on the following elliptical extension of Lemma~\ref{LemLL}.

\begin{Lem}
\label{LemLLellipt}
Fix a unit $p$-vector~$\thetab_1^0$, $v>0$, a bounded positive real sequence $(r_n)$, and~$f\in\mathcal{F}$. Then, under ${\rm P}_{\thetab_1,r_n,v,f}$, 
	$\sqrt{n}(\Sigb_n^{-1/2})^{\otimes 2} {\rm vec}
	({\bf S}_n- \Sigb_n)$ is asymptotically normal with mean zero and covariance matrix~$(1+ \kappa_p(f))({\bf I}_{p^2}+ {\bf K}_p)+\kappa_p(f) ({\rm vec}\,{\bf I}_p) ({\rm vec}\,{\bf I}_p)\pr$. In particular, (i) if $r_n\equiv 1$, then $\sqrt{n} \, {\rm vec}
	({\bf S}_n- \Sigb_n)$ is asymptotically normal with mean zero and covariance matrix~$(1+ \kappa_p(f))({\bf I}_{p^2}+ {\bf K}_p)(\Sigb(v))^{\otimes 2}+\kappa_p(f) ({\rm vec}\,\Sigb(v)) ({\rm vec}\,\Sigb(v))\pr$, still with $\Sigb(v)
	%
	%
	:={\bf I}_p+ v\, \thetab_1\thetab_1\pr$; 
	(ii) if~$r_n$ is~$o(1)$, then~$\sqrt{n} \, {\rm vec}
	({\bf S}_n- \Sigb_n)$ is asymptotically normal with mean zero and covariance matrix~$(1+ \kappa_p(f))({\bf I}_{p^2}+ {\bf K}_p)+\kappa_p(f) ({\rm vec}\,{\bf I}_p) ({\rm vec}\,{\bf I}_p)\pr$.
\end{Lem}

The main result of this section is the following theorem, that in particular extends Theorem~\ref{TheorHPVNull} to the elliptical setup.

\begin{Theor} 
\label{TheorHPVNullellipt}
Fix a unit $p$-vector~$\thetab_1^0$, $v>0$, a bounded positive real sequence~$(r_n)$, and~$f\in\mathcal{F}$. Then, under~${\rm P}_{\thetab_1^0,r_n,v, f}$, 
$$
Q_{\rm HPV}^{(n)\dagger} \stackrel{\mathcal{D}}{\to} \chi^2_{p-1},
$$
so that, in all regimes~(i)--(iv), the test~$\phi_{\rm HPV}^\dagger$, irrespective of the radial density~$f\in\mathcal{F}$, has asymptotic size~$\alpha$ under the null hypothesis. Moreover, under~${\rm P}_{\thetab_1^0,r_n,v,\phi}$,  
\begin{equation}
	\label{asymequipse}
Q_{\rm HPV}^{(n)\dagger}=Q_{\rm HPV}^{(n)}+o_{\rm P}(1)
\end{equation}
as~$n\to\infty$.   
\end{Theor}

This result shows that the pseudo-Gaussian version~$\phi_{\rm HPV}^\dagger$ of~$\phi_{\rm HPV}$ is robust to weak identifiability under any elliptical distribution with finite fourth-order moments. Since the asymptotic equivalence in~(\ref{asymequipse}) extends, from contiguity, to the (Gaussian) local
\vspace{-.5mm}
  alternatives identified in Theorem~\ref{TheorLAN}, it also directly follows from Theorem~\ref{TheorHPVNullellipt} that~$\phi_{\rm HPV}^\dagger$ inherits the optimality properties of~$\phi_{\rm HPV}$ in the multinormal case. For the sake of completeness, we mention that, by using elliptical extensions of Lemmas~\ref{Lemeigenvalues}--\ref{Lemeigenvectors} (see Lemmas~\ref{Lemeigenvaluesellipt}--\ref{Lemeigenvectorsellipt} in the appendix), it can be showed that, irrespective of the elliptical distribution considered, the pseudo-Gaussian test~$\phi_{\rm A}^\dagger$ asymptotically meets the nominal level constraint in regimes~(i)--(ii) only, hence is not robust to weak identifiability. Remarkably, by using the same results, it can also be showed that, under any bivariate elliptical distribution with finite fourth-order moments, the null asymptotic distribution of~$Q_{\rm A}^{(n)\dagger}$ is still~$4\chi^2_1$ in regime~(iv), which extends Corollary~\ref{CorolAnd} to the elliptical setup.

\begin{figure}[htbp!]
\begin{center} 
\includegraphics[width=\textwidth]{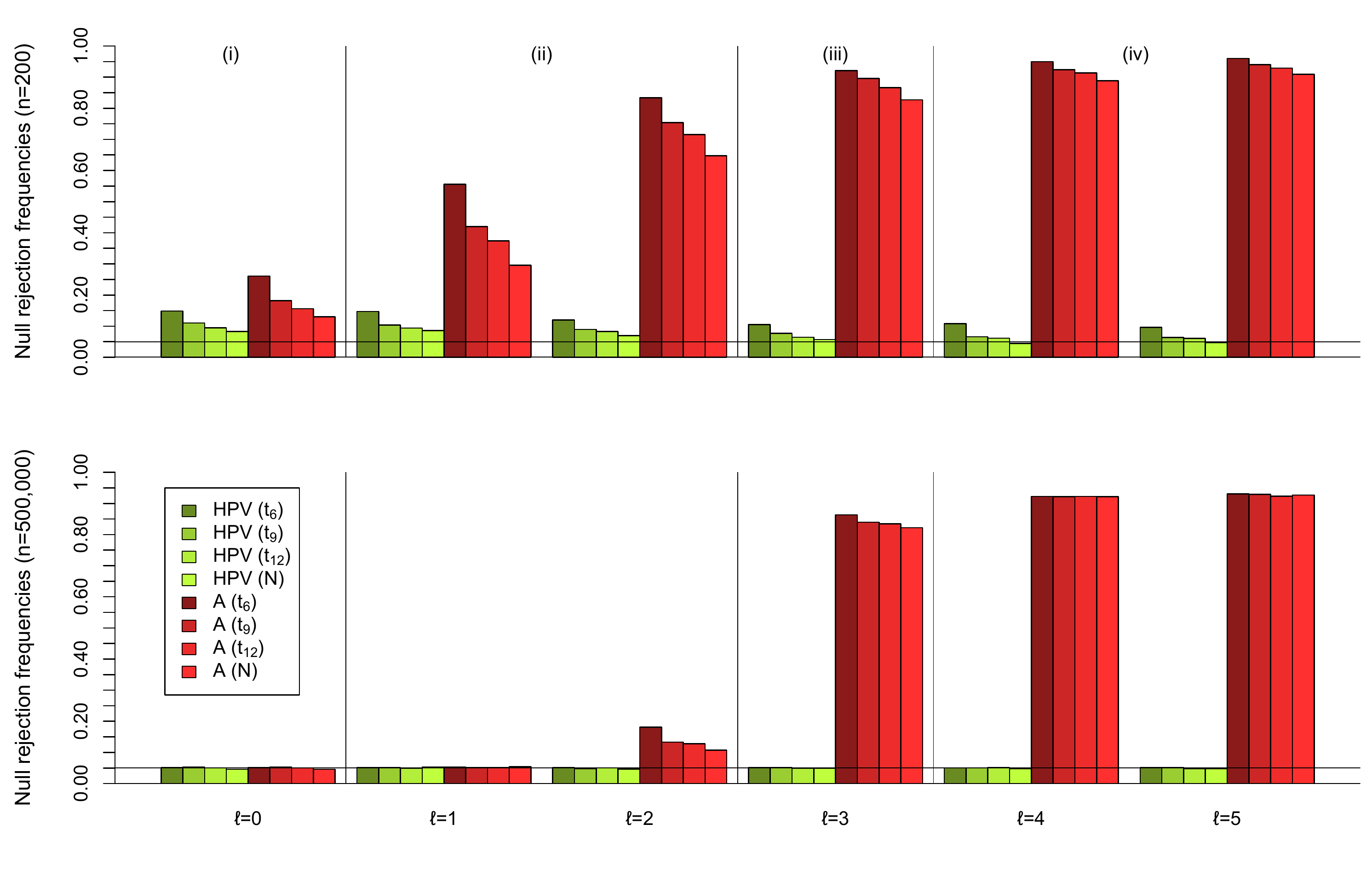}
\vspace{-5mm}
\caption{Empirical rejection frequencies, under the null hypothesis, of the tests~$\phi^\dagger_{\rm HPV}$ and~$\phi^\dagger_{\rm A}$ performed at nominal level~$5\%$. Results are based on $M=10,\!000$ independent ten-dimensional random samples of size~$n=200$ and size~$n=500,\!000$, drawn from $t_6$, $t_9$, $t_{12}$ and Gaussian distributions. Increasing values of~$\ell$ bring the underlying spiked covariance matrix closer and closer to a multiple of the identity matrix; see Section~\ref{sec:ellipt} for details.}
\label{Fig5}
\end{center}
\end{figure}

We now quickly illustrate these results through a Monte Carlo exercise that extends to the elliptical setup the one conducted in Figure~\ref{Fig1}. To do so, for any~$\ell=0,1,\ldots,5$, we 
\vspace{-.4mm}
 generated $M=10,\!000$ mutually independent random samples~$\Xb_i^{(\ell,s)}\!$, $i=1, \ldots, n$, from the $(p=10)$-variate $t_6$ ($s=1$), $t_9$ ($s=2$), $t_{12}$ ($s=3$), and normal ($s=4$) distributions with mean zero and covariance matrix~$\Sigb_n^{(\ell)}:= {\bf I}_p+n^{-\ell/6} \thetab_1^0 \thetab_1^{0\prime}$, where~$\thetab_1^0$ is still the first vector of the canonical basis of~$\R^p$. As in Figure~\ref{Fig1}, this covers regimes~(i) ($\ell=0$), (ii) ($\ell=1,2$), (iii) ($\ell=3$), and (iv) ($\ell=4,5$). Figure~\ref{Fig5} reports, for~$n=200$ and~$n=500,\!000$, the resulting rejection frequencies
\vspace{-.5mm}
of the pseudo-Gaussian tests~$\phi_{\rm HPV}^\dagger$ and~$\phi_{\rm A}^\dagger$ for ${\cal H}_0\n: \thetab_1= \thetab_1^0$ at nominal level~$5\%$.  Clearly, the results confirm that, irrespective of the underlying elliptical distribution, the pseudo-Gaussian HPV test is robust to weak identifiability, while the pseudo-Gaussian Anderson test meets the asymptotic level constraint only in regimes (i)--(ii) (this test strongly overrejects the null hypothesis in other regimes).

\section{Real data example}  
\label{sec:real}

We now provide a real data illustration on the celebrated Swiss banknote dataset, which has been considered in numerous multivariate statistics monographs, such as  \cite{FluRie1988}, \cite{Atketal2004}, \cite{HarSim2007} and \cite{Koc2013}, but also in many research papers; see, e.g., \cite{Saletal2006} or \cite{BurPol2009}. The dataset, that is available in the R package \verb|uskewfactors| (\citealp{uskew}), offers six measurements on 100 genuine and 100 counterfeit old Swiss 1000-franc banknotes. This dataset was often used to illustrate various multivariate statistics procedures such as, e.g., linear discriminant analysis (\citealp{FluRie1988}), principal component analysis (\citealp{Flu88}), or independent component analysis (\citealp{Gir1999}); we also refer to \cite{Shi2016} for a recent account on discriminant analysis for this dataset.

Here, we aim to complement the PCA analysis conducted in \cite{Flu88} (see pp.~41--43), hence use the exact same subset of the Swiss banknote data as the one considered there. More precisely, (i) we focus on four of the six available measurements, namely the width~$L$ of the left side of the banknote, the width~$R$ on its right side, the width~$B$ of the bottom margin and the width~$T$ of the top margin, all measured in~mm$\times 10^{-1}$ (rather than in the original mm); (ii) we also restrict to $n=85$ counterfeit bills made by the same forger (it is well-known that the $100$ counterfeit bills were made by two different forgers; see, e.g., \cite{FluRie1988}, page~250, or \cite{Frietal2012}, page~22).
Letting~$c_n=(n-1)/n\approx .99$, the resulting sample covariance matrix is 
$$
{\bf S}
= 
c_n
\left(
\begin{array}{cccc} 
6.41 & 4.89 & 2.89 &  -1.30 \\ 
4.89 & 9.40 & -1.09 & 0.71 \\ 
2.89 & -1.09 & 72.42 & -43.30 \\ 
-1.30 & 0.71& -43.30 & 40.39 
\end{array}
\right)
,
$$
with eigenvalues of $\hat{\lambda}_1=102.69c_n$, $\hat{\lambda}_2=13.05c_n$, $\hat{\lambda}_3=10.23c_n$ and $\hat{\lambda}_4=2.66c_n$, and corresponding eigenvectors
$$
\label{eigenvectorsrealdata}
\hat{\thetab}_1 
=
\left(
\!
\begin{array}{c} 
.032 \\
 -.012 \\ 
 .820 \\ 
 -.571 
 \end{array}
\!
 \right)
\!
,
\
\hat{\thetab}_2
=
\left(
\!
\begin{array}{cccc} 
 .593  \\
 .797 \\ 
 .057 \\ 
 .097 
 \end{array}
\!
 \right)
\!
,
\
\hat{\thetab}_3
=
\left(
\!
\begin{array}{cccc} 
-.015 \\
-.129 \\ 
.566  \\ 
.814 
 \end{array}
\!
 \right)
\!
,
\textrm{ and }
\hat{\thetab}_4
=
\left(
\!
\begin{array}{cccc} 
 .804 \\
 -.590 \\ 
-.064 \\ 
 -.035 
 \end{array}
\!
 \right)
\!
;
$$
the unimportant factor~$c_n$ is used here to ease the comparison with \cite{Flu88}, where the unbiased version of the sample covariance matrix was adopted throughout. From these estimates, Flury concludes that the first principal component is a contrast between~$B$ and~$T$, hence can be interpreted as the vertical position of the print image on the bill. It is tempting to interpret the second principal component as an aggregate of~$L$ and~$R$, that is, essentially as the vertical size of the bill. Flury, however, explicitly writes \emph{``beware: the second and third roots are quite close and so the computation of standard errors for the coefficients of~$\hat{\thetab}_2$ and~$\hat{\thetab}_3$ may be hazardous"}. He reports that these eigenvectors should be considered spherical and that the corresponding standard errors should be ignored. In other words, Flury, due to the structure of the spectrum, refrains from drawing any conclusion about the second component. 
 
The considerations above make it natural to test that~$L$ and~$R$ contribute equally to the second principal component and that they are the only variables to contribute to it. In other words, it is natural to test the null hypothesis~${\cal H}_0: \thetab_2=\thetab_2^0$, with~$\thetab_2^0:=(1,1,0,0)\pr/\sqrt{2}$. While the tests discussed in the present paper address testing problems on the first eigenvector~$\thetab_1$, obvious modifications of these tests allow performing inference on any eigenvector~$\thetab_j$, $j=2, \ldots, p$. In particular, 
\vspace{-.6mm}
 the Anderson test~$\phi\n_{\rm A}$ and HPV test~$\phi\n_{\rm HPV}$ for~${\cal H}_0: \thetab_2= \thetab_2^0$ against~$ {\cal H}_1: \thetab_2 \neq \thetab_2^0$ reject the null hypothesis at asymptotic level~$\alpha$ whenever
\begin{equation} 
\label{Andsec}
n \big( \hat{\lambda}_{2} \thetab_{2}^{0\prime} {\bf S}^{-1} \thetab_{2}^0 +  \hat{\lambda}_{2}^{-1} \thetab_{2}^{0\prime} {\bf S}\, \thetab_{2}^0-2\big) > \chi_{p-1, 1-\alpha}^2
\end{equation}
and
\begin{equation} 
\label{HPVsec}
 \frac{n}{\hat{\lambda}_{2}}  \sum_{j=1,j\neq 2}^p \hat{\lambda}_{j}^{-1} 
 \big(\tilde{\thetab}_j\pr {\bf S} \thetab_{2}^0\big)^2 > \chi_{p-1, 1-\alpha}^2
 ,
\end{equation}
respectively, where, parallel to~(\ref{GS}), $\tilde{\thetab}_1, \thetab_{2}^0, \tilde{\thetab}_3, \ldots, \tilde{\thetab}_p$ results from a Gram-Schmidt orthogonalization of~$\hat{\thetab}_1, \thetab_{2}^0, \hat{\thetab}_3, \ldots, \hat{\thetab}_p$. 
When applied with~$\thetab_2^0:=(1,1,0,0)\pr/\sqrt{2}$, this HPV test provides a $p$-value equal to~$.177$, hence does not lead to rejection of the null hypothesis at any usual nominal level. In contrast, the $p$-value of the Anderson test in~(\ref{Andsec}) is~$.099$, so that this test rejects the null hypothesis at the level $10\%$. Since the results of this paper show that the Anderson test tends to strongly overreject the null hypothesis when eigenvalues are close, practitioners should here be confident that the HPV test provides the right decision. 

\begin{figure}[htbp!]
\begin{center} 
\includegraphics[width=80mm]{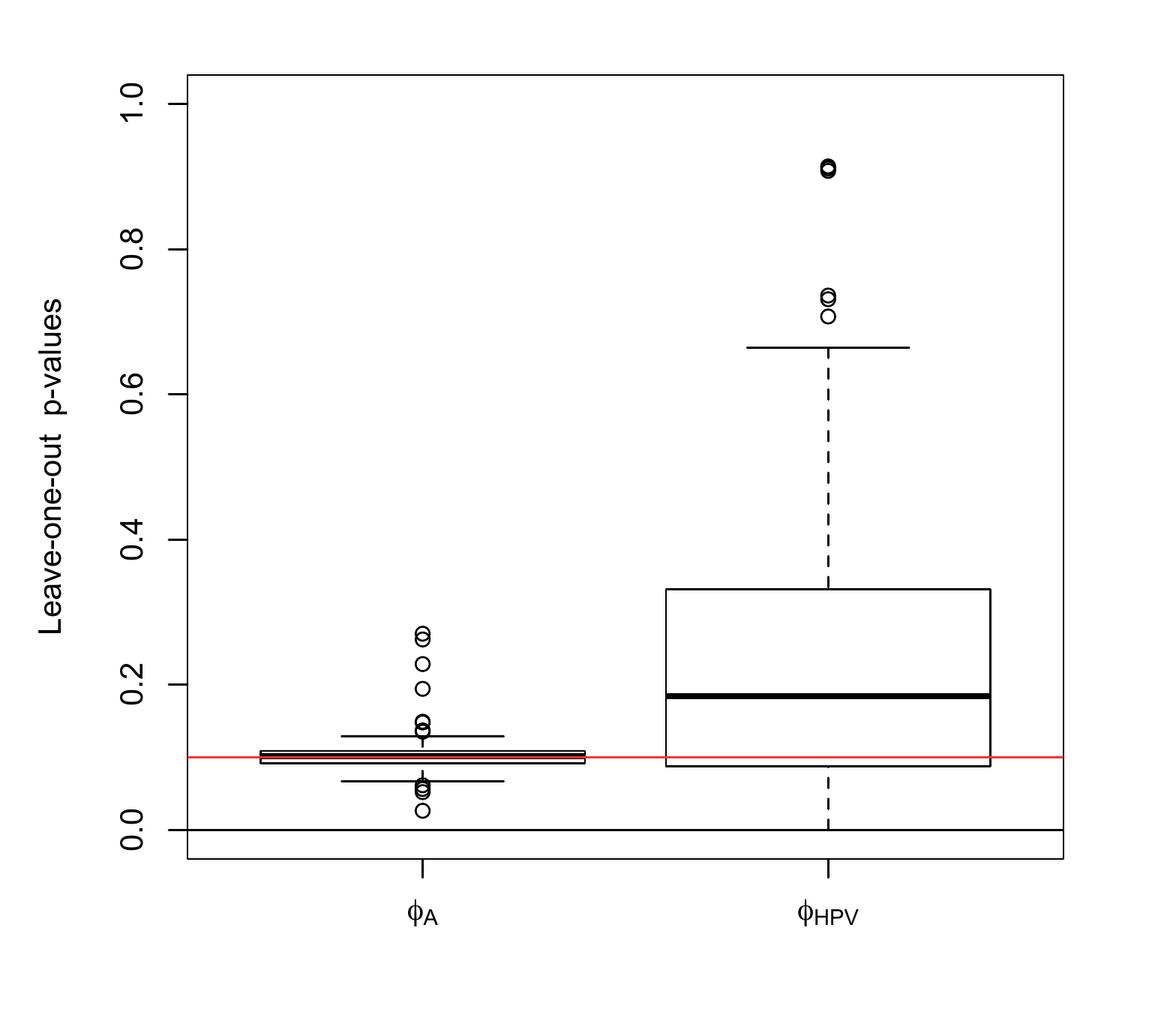}
\vspace{-5mm}
\caption{Boxplots of the 85 ``leave-one-out" $p$-values of the Anderson test in~(\ref{Andsec}) (left) and HPV test in~(\ref{HPVsec}) (right) when testing the null hypothesis~$\mathcal{H}_0:\thetab_2:=(1,1,0,0)\pr/\sqrt{2}$. More precisely, these $p$-values are those obtained when applying the corresponding tests to the 85 subsample of size~$84$ obtained by removing one observation in the real data set considered in the PCA analysis of~\cite{Flu88}, pp.~41--43.}
\label{Fig6}
\end{center}
\end{figure} 

To somewhat assess the robustness of this result, we performed the same HPV and Anderson tests on the~85 subsamples obtained by removing one observation from the sample considered above. For each test, a boxplot of the resulting 85 ``leave-one-out" $p$-values is provided in Figure~\ref{Fig6}. Clearly, these boxplots reveal that the Anderson test rejects the null hypothesis much more often than the HPV test. Again, the results of the paper provide a strong motivation to rely on the outcome of the HPV  test in the present context.


\section{Wrap up and perspectives}
\label{sec:wrapup} 

In this paper, 
\vspace{-.3mm}
we tackled the problem of testing the null
\vspace{-.5mm}
  hypothesis~${\cal H}_0\n: \thetab_1= \thetab_1^0$ against the alternative~${\cal H}_1\n: \thetab_1 \neq \thetab_1^0$, where~$\thetab_1$ is the eigenvector associated with the largest eigenvalue of the underlying covariance matrix and where~$\thetab_1^0$ is some fixed unit vector. We analyzed the asymptotic behavior of the classical \cite{And63} test~$\phi_{\rm A}$ and of the \cite{HPV10} test~$\phi_{\rm HPV}$ under sequences of $p$-variate Gaussian models with spiked covariance matrices of the form~$\Sigb_n=\sigma_n^2 ({\bf I}_p+ r_n v \thetab_1\thetab_1\pr)$, where~$(\sigma_n)$ is a positive sequence,~$v>0$ is fixed, and~$(r_n)$ is a positive sequence that converges to zero. We showed that in these situations where~$\thetab_1$ is closer and closer to being unidentified, $\phi_{\rm HPV}$ performs better than $\phi_{\rm A}$: (i) $\phi_{\rm HPV}$, unlike~$\phi_{\rm A}$, meets asymptotically the nominal level constraint without any condition on the rate at which~$r_n$ converges to zero, and (ii) $\phi_{\rm HPV}$ remains locally asymptotically maximin in all regimes, but in the contiguity regime~$r_n=1/\sqrt{n}$ where~$\phi_{\rm HPV}$ still enjoys the same optimality property locally in $\taub$. These considerations, along with the asymptotic equivalence of~$\phi_{\rm HPV}$ and~$\phi_{\rm A}$ in the standard case~$r_n\equiv 1$, clearly imply that the test~$\phi_{\rm HPV}$, for all practical purposes, should be favored over~$\phi_{\rm A}$, all the more so that the results above extend to elliptical distributions if the Anderson and HPV tests are replaced with their pseudo-Gaussian versions~$\phi_{\rm A}^\dagger$ and~$\phi_{\rm HPV}^\dagger$.
  
To conclude, we discuss some research perspectives. Throughout the paper, we assumed that the dimension $p$ is fixed. Since PCA is often used for dimension reduction, it would be of interest to consider tests that can cope with high-dimensional situations where~$p$ is as large as~$n$ or even larger than~$n$, and to investigate the robustness of these tests to weak identifiability. The tests considered in the present paper, however, are not suitable in high dimensions. This is clear for the Anderson test~$\phi_{\rm A}$ since this test requires inverting the sample covariance matrix~$\mathbf{S}_n$, that fails to be invertible for~$p\geq n$. As for the HPV test~$\phi_{\rm HPV}$, our investigation of the asymptotic behavior of this test in the fixed-$p$ case crucially relied on the consistency of the eigenvalues~$\hat{\lambda}_{nj}$ of~$\mathbf{S}$; in high-dimensional regimes where~$p=p_n \to \infty$ so that~$p_n/n \to c$, however, these sample eigenvalues are no longer consistent (see, e.g., \citealp{Baietal2005}), which suggests that~$\phi_{\rm HPV}$ is not robust to high dimensionality. To explore this, we conducted the following Monte Carlo exercise: for~$n=200$ and each value of $p=cn$, with~$c\in\{0.5,0.75,1,1.5,2\}$, we generated $2,\!000$ mutually independent random samples~$\Xb_{1}, \ldots, \Xb_n$ from the $p$-variate normal distribution with mean zero and covariance matrix~$\Sigb= {\bf I}_p + \thetab_1^0 \thetab_1^{0 \prime}$, where~$\thetab_1^0$ is the first vector of the canonical basis of~$\R^p$. The resulting rejection frequencies of the test~$\phi_{\rm HPV}$ (resp., of the test~$\phi_{\rm A}$), conducted at asymptotic level~$5\%$, are 
$0.9255$ (resp.,~1) for~$c=0.5$,  
$0.9240$ (resp.,~1) for~$c=0.75$, 
$0.5000$ (resp., ---) for~$c=1$, 
$0.1985$ (resp., ---) for~$c=1.5$, 
and
$0.1715$ (resp., ---) for~$c=2$ (as indicated above, the Anderson test cannot be used for~$p\geq n$). This confirms that neither~${\phi}_{\rm HPV}$ nor~${\phi}_{\rm A}$ can cope with high dimensionality. As a result, the problem of providing a suitable test in the high-dimensional setup and of studying its robustness to weak identifiability is widely open and should be investigated in future research. 

\appendix
 

\section{Technical proofs for Section~2}
\label{SecSupProofLAN}


Throughout this appendix, we will write~$v_n:=r_n v$ and~$\delta$ will take value one if regime~(i) is considered and zero otherwise. All convergences are as~$n\to\infty$.  
\vspace{3mm}

{\bf Proof of Lemma~\ref{Lemeigenvalues}}. Fix arbitrarily~$\thetab_2, \ldots, \thetab_p$ such that the $p \times p$ matrix~$\Gamb:=(\thetab_1^0,\thetab_2, \ldots, \thetab_p)$ is orthogonal.
Letting~$\Lamb_n:={\rm diag}(1+v_n, 1, \ldots, 1)$, we have $\Sigb_n \Gamb=\Gamb\Lamb_n$, so that $\Gamb$ is an eigenvectors matrix for $\Sigb_n$.
Clearly, $\ell_{n1}:=\sqrt{n} (\hat{\lambda}_{n1}-(1+v_n))$ is the largest root of the polynomial~$P_{n1}(h):={\rm det}({\sqrt n}({\bf S}_n-(1+v_n) {\bf I}_p)- h {\bf I}_p) 
$ and $\ell_{n2}:=\sqrt{n}(\hat{\lambda}_{n2}-1), \ldots, \ell_{np}:=\sqrt{n}(\hat{\lambda}_{np}-1)$ are the $p-1$ smallest roots of the polynomial~$P_{n2}(h):={\rm det}({\sqrt n}({\bf S}_n- {\bf I}_p)- h {\bf I}_p)$. Letting
\begin{equation}
	\label{defZn}
{\bf Z}_n
:= 
\sqrt{n}(\Gamb\pr {\bf S}_n \Gamb- \Lamb_n)
=
 \Gamb\pr \sqrt{n}({\bf S}_n- \Sigb_n) \Gamb
 ,
\end{equation}
a key ingredient in this proof is to rewrite these polynomials as
\begin{eqnarray}
P_{n1}(h)
&\!\!\!=\!\!\!& {\rm det}(\sqrt{n}\Gamb\pr {\bf S}_n \Gamb- \sqrt{n}(1+v_n) {\bf I}_p - h {\bf I}_p) \nonumber  \\
&\!\!\!=\!\!\!& {\rm det}({\Zb}_n+\sqrt{n}(\Lamb_n-(1+v_n) {\bf I}_p)- h {\bf I}_p)
\label{firsteq}
\end{eqnarray}
and
\begin{eqnarray} 
P_{n2}(h)
&\!\!\!=\!\!\!& {\rm det}({\sqrt n}\Gamb\pr {\bf S}_n \Gamb- {\sqrt n}{\bf I}_p - h {\bf I}_p) \nonumber \\
&\!\!\!=\!\!\!& {\rm det}({\bf Z}_n+ {\sqrt n}(\Lamb_n- {\bf I}_p)- h {\bf I}_p)
.
\label{seceq}
\end{eqnarray}
Note that Lemma~\ref{LemLL} readily implies that~$\Zb_n$ converges weakly to~$\Zb(v)$ in case~(i) and to~$\Zb=\Zb(0)$ in cases~(ii)--(iv), where~$\Zb(v)$ is the random matrix defined in the statement of the theorem. In all cases, thus, $\Zb_n$ converges weakly to~$\Zb(\delta v)$, where~$\delta$ was introduced at the beginning of this appendix.  

We start with the proofs of~(i)--(ii). Partition~$\Zb_n$ and~$\Zb(\delta v)$ into 
\begin{equation}
\label{partitionZ}	
\Zb_n 
=
\bigg(
\begin{array}{cc}
Z_{n,11} & \Zb_{n,21}' \\[0mm]	
\Zb_{n,21} & \Zb_{n,22} 
\end{array}
\bigg)
\quad
\textrm{and}
\quad 
\Zb(\delta v)
=
\bigg(
\begin{array}{cc}
Z_{11}(\delta v) & \Zb_{21}'(\delta v) \\[0mm]	
\Zb_{21}(\delta v) & \Zb_{22}(\delta v) 
\end{array}
\bigg)
,
\end{equation}
where~$Z_{n,11}$ and~$Z_{11}(\delta v)$ are random variables, whereas~$\Zb_{n,22}$ and~$\Zb_{22}(\delta v)$ are $(p-1)\times (p-1)$ random matrices. It follows from the discussion above that $\ell_{n1}$ is the largest root of the polynomial
\begin{eqnarray*}
Q_{n1}(h)
&\!\!\!:=\!\!\!&
\frac{P_{n1}(h)}{(\sqrt{n} v_n)^{p-1}}
\\[2mm]
&\!\!\!=\!\!\!&
\frac{1}{(\sqrt{n} v_n)^{p-1}}
\,
{\rm det}\left(\begin{array}{cc} Z_{n,11}-h &  \Zb_{n,21}' \\ 
\Zb_{n,21} & \Zb_{n,22}- (\sqrt{n} v_n+h) {\bf I}_{p-1} \end{array}\right) 
\\[2mm]
&\!\!\!=\!\!\!&
{\rm det}\left(\begin{array}{cc} 
Z_{n,11}-h &  (\sqrt{n} v_n)^{-1} \Zb_{n,21}' \\
\Zb_{n,21} &(\sqrt{n} v_n)^{-1}( \Zb_{n,22}- (\sqrt{n} v_n+h) {\bf I}_{p-1}) \end{array}\right) 
.
\end{eqnarray*}
In cases~(i)--(ii), $(\sqrt{n} v_n)^{-1} \Zb_{n,21}$ is~$o_{\rm P}(1)$, so that the largest root~$\ell_{n1}$ converges weakly to the root of the weak limit of~$Z_{n,11}-h$, namely to~$Z_{11}(\delta v)$ (note that, still since $(\sqrt{n} v_n)^{-1} \Zb_{n,21}$ is~$o_{\rm P}(1)$, the~$p-1$ smallest roots of~$Q_{n1}(h)$ converge to~$-\infty$ in probability). Similarly, $\ell_{n2},\ldots, \ell_{np}$ are the $p-1$ smallest roots (in decreasing order) of the polynomial         
\begin{eqnarray*}
Q_{n2}(h)
&\!\!\!:=\!\!\!&
\frac{P_{n2}(h)}{\sqrt{n} v_n}
\\[2mm]
&\!\!\!=\!\!\!&
\frac{1}{\sqrt{n} v_n}
\,
{\rm det}\left(\begin{array}{cc} {Z}_{n,11}+ \sqrt{n} v_n-h &  {\bf Z}_{n,21}' \\  {\bf Z}_{n,21} & {\bf Z}_{n,22}- h {\bf I}_{p-1} \end{array}\right) 
\\[2mm]
&\!\!\!=\!\!\!&
{\rm det} \left(\begin{array}{cc} (\sqrt{n} v_n)^{-1}({Z}_{n,11}+ \sqrt{n} v_n-h) &  {\bf Z}_{n,12}' \\  
(\sqrt{n} v_n)^{-1}{\bf Z}_{n,21} & {\bf Z}_{n,22}- h {\bf I}_{p-1} \end{array}\right) 
.
\end{eqnarray*}
Using again the fact that $(\sqrt{n} v_n)^{-1} \Zb_{n,21}$ is~$o_{\rm P}(1)$ yields that~$(\ell_{n2},\ldots, \ell_{np})'$ converges weakly to the vector of roots (still in decreasing order) of the weak limit of the polynomial~$\det({\bf Z}_{n,22}(\delta v) - h {\bf I}_{p-1})=0$, namely to the vector of (ordered) eigenvalues of~${\bf Z}_{22}(\delta v)$; note that the largest root of~$Q_{n2}(h)$ converges to~$\infty$ in probability. The result then follows from the fact that, in both cases~(i)--(ii),~$Z_{11}(\delta v)\sim \mathcal{N}(0,{2(1+\delta v)^2})$ and~$\Zb_{22}(\delta v)$ are mutually independent.    

We turn to the proof of~(iii). In this case, we have
$
P_{n1}(h)
=
{\rm det}({\Zb}_n- {\rm diag}(0, v, \ldots, v)- h {\bf I}_p)
$
and 
$
P_{n2}(h)
=
{\rm det}({\Zb}_n + {\rm diag}(v,0, \ldots, 0)- h {\bf I}_p)
$. It readily follows that~$\ell_{n1}$ converges weakly to the largest root of the polynomial~${\rm det}({\Zb}- {\rm diag}(0, v, \ldots, v)- h {\bf I}_p)$, that is, to the largest eigenvalue of~${\Zb}- {\rm diag}(0, v, \ldots, v)$, while~$(\ell_{n2},\ldots, \ell_{np})'$ converges weakly to the vector of the $p-1$ smallest roots (in decreasing order) of the polynomial~${\rm det}({\Zb} + {\rm diag}(v,0, \ldots, 0)- h {\bf I}_p)$, namely to the vector of the $p-1$ smallest (ordered) eigenvalues of ${\Zb} + {\rm diag}(v,0, \ldots, 0)$. 
 
Finally, we prove the result in~(iv). In that case, we have 
$$
P_{n1}(h)
=
{\rm det}\big({\Zb}_n- \sqrt{n}v_n {\rm diag}(0, 1, \ldots,1) - h {\bf I}_p \big)
$$
and 
$$
P_{n2}(h)
=
{\rm det}\big({\Zb}_n + \sqrt{n}v_n {\rm diag}(1,0, \ldots, 0)- h {\bf I}_p\big)
,
$$
where~$\sqrt{n}v_n$ is~$o(1)$. Parallel as above, it follows that~$\ell_{n1}$ (resp.,~$(\ell_{n2},\ldots, \ell_{np})'$) converges weakly to the largest root (resp., the $p-1$ smallest roots, in decreasing order) of the polynomial~${\rm det}({\Zb} - h {\bf I}_p)$. In other words, ${\pmb \ell}=(\ell_1, \ldots,\ell_p)'$ is equal in distribution to the vector collecting the eigenvalues (in decreasing order) of~$\Zb$. It only remains to establish that the density of~$\ell$ is the one given in~(\ref{explidensl}). To do so, let~$\Db_p$ be the $p$-dimensional \emph{duplication matrix}, that is such that~$\Db_p({\rm vech}\,\Ab)={\rm vec}\,\Ab$ for any~$p\times p$ symmetric matrix~$\Ab$, where~${\rm vech}\,\Ab$ is the $p(p+1)/2$-vector obtained from~${\rm vec}\,\Ab$ by restricting to the entries from the lower-triangular part of~$\Ab$. Consider then~$\Wb:={\rm vech}\,\Zb=\Db_p^-({\rm vec}\,\Zb)$, where~$\Db_p^-=(\Db_p'\Db_p)^{-1}\Db_p'$ is the Moore-Penrose inverse of~$\Db_p$; see page~57 in~\cite{MagNeu2007}. Clearly, $\Wb$ has density
\begin{eqnarray*}
{\bf w}\mapsto g({\bf w})
& \!\! = \!\! &
\Big(\frac{1}{2\pi}\Big)^{p(p+1)/2}
\exp\Big(-\frac{1}{2}\, {\bf w}' (\Db_p^-({\bf I}_{p^2}+ {\bf K}_p)(\Db_p^-)')^{-1} {\bf w} \Big)
\\[2mm]
& \!\! = \!\! &
\Big(\frac{1}{2\pi}\Big)^{p(p+1)/2}
\exp\Big(-\frac{1}{4}\, {\bf w}' \Db_p'\Db_p {\bf w} \Big)
,
\end{eqnarray*}
where we used the identity~${\bf K}_{p}\Db_p=\Db_p$; see again page~57 in~\cite{MagNeu2007}. The resulting density for~$\Zb$ is therefore
$$
\zb \mapsto f(\zb)=g({\rm vech}\,\zb)
=
\Big(\frac{1}{2\pi}\Big)^{p(p+1)/2}
\exp\Big(-\frac{1}{4} ({\rm tr}\,\zb^2) \Big)
.
$$
The density of~${\pmb \ell}$ in~(\ref{explidensl}) then follows from Lemma~2.3 of \cite{Tyl83b}. 
\cqfd
\vspace{4mm}


%
%
%
%
%
%

 The proof of Lemma~\ref{Lemeigenvectors} requires the following linear algebra result.
 
 \begin{Lem} 
\label{cofactors}
Let~$\Ab$ be a $p\times p$ matrix and~$\lambda$ be an eigenvalue of~$\Ab$. Assume that~$\vb=(C_{11},\ldots,C_{1p})'\neq {\bf 0}$, where~$C=(C_{ij})$ stands for the cofactor matrix of~$\Ab-\lambda\mathbf{I}_p$. Then~$\vb$ is an eigenvector of~$\Ab$ associated with eigenvalue~$\lambda$. 
\end{Lem}

 \noindent {\bf Proof of Lemma~\ref{cofactors}}.
For any~$j=1,\ldots,p$, denote as~$(\Ab-\lambda \mathbf{I}_p)_j$ the $j$th row of~$\Ab-\lambda\mathbf{I}_p$ (left as a row vector). Since~$\lambda$ is an eigenvalue of~$\Ab$, we then have~$(\Ab-\lambda \mathbf{I}_p)_1\vb=\det(\Ab-\lambda \mathbf{I}_p)=0$. Now, for~$j=2,\ldots,p$,  
 $$
(\Ab-\lambda \mathbf{I}_p)_j\vb
=
\det
\left(
\begin{array}{c}
(\Ab-\lambda \mathbf{I}_p)_j\\[0mm]	
(\Ab-\lambda \mathbf{I}_p)_2\\[0mm]	
\vdots\\[0mm]	
(\Ab-\lambda \mathbf{I}_p)_p\\[0mm]	
\end{array}
\right)
=
0
,
$$
since this is the determinant of a matrix with (at least) twice the same row. We conclude that~$(\Ab-\lambda \mathbf{I}_p)\vb=\mathbf{0}$. Since~$\vb\neq {\bf 0}$, this shows that~$\vb$ is an eigenvector of~$\Ab$ associated with eigenvalue~$\lambda$.   
\cqfd
\vspace{4mm}


 \noindent {\bf Proof of Lemma~\ref{Lemeigenvectors}}.
First note that since 
 \begin{eqnarray*}
{\bf E}_n= \hat{\Gamb}_n\pr \Gamb
=
\left( 
\begin{array}{cc} E_{n,11} & {\bf E}_{n,12} \\
 {\bf E}_{n,21} & {\bf E}_{n,22} \end{array} 
 \right)
\end{eqnarray*}
is an orthogonal matrix, we easily obtain that
\begin{equation} 
\label{12eq}
{\bf E}_{n,21}= -\frac{1}{E_{n,11}} {\bf E}_{n,22} {\bf E}_{n,12}\pr 
,
\end{equation} 
\begin{equation} 
\label{14eq}
{\bf E}_{n,22} {\bf E}_{n,22}' = \mathbf{I}_{p-1} - {\bf E}_{n,21} {\bf E}_{n,21}'
\end{equation} 
and
\begin{equation} 
\label{13eq}
E_{n,11} {\bf E}_{n,12}' = - {\bf E}_{n,22}\pr {\bf E}_{n,21}
.
\end{equation} 
We start with the proof of~(i)--(ii). 
As in the proof of Lemma~\ref{Lemeigenvalues}, we have from \eqref{firsteq} that, still with the random vector~${\bf Z}_n$ defined in~(\ref{defZn}), $\ell_{n1}=\sqrt{n} (\hat{\lambda}_{n1}-(1+v_n))$ is the largest eigenvalue of
 \begin{eqnarray}
 \sqrt{n}\Gamb\pr {\bf S}_n \Gamb- \sqrt{n}(1+v_n) {\bf I}_p 
  &\!\!=\!\!& {\bf Z}_n+ {\sqrt n}(\Lamb_n- (1+v_n){\bf I}_p) 
  \nonumber
  \\
  &\!\!=\!\!& {\bf Z}_n-{\rm diag}(0, \sqrt{n}v_n,\ldots,\sqrt{n}v_n)
.
\label{equZlam1}
 \end{eqnarray}
One readily checks that~$\wb_{n1}:=\Gamb\pr \hat{\thetab}_{n1}=(E_{n,11}, {\bf E}_{n,12})\pr$ is a corresponding unit eigenvector. Using the same notation as in the proof of Lemma~\ref{Lemeigenvalues}, Lemma~\ref{cofactors} yields that~$\wb_{n1}$ is proportional to the vector of cofactors associated with the first row of
 \begin{equation}
 \label{equZlam1bis}	
{\bf M}_{n,1}:=\left(\begin{array}{cc} {Z}_{n,11}-\ell_{n1} &  {\bf Z}_{n,21}' \\
 {\bf Z}_{n,21} & {\bf Z}_{n,22} - \sqrt{n}v_n {\bf I}_{p-1} - \ell_{n1} {\bf I}_{p-1} \end{array} \right)
,
\end{equation}
or equivalently, of 
\begin{eqnarray*}
\frac{1}{(\sqrt{n} v_n)^{p-1}}
\,
{\bf M}_{n,1}
=
\left(\begin{array}{cc} {Z}_{n,11}-\ell_{n1} &  {\bf Z}_{n,21}' \\ [1mm]
\frac{1}{\sqrt{n} v_n} {\bf Z}_{n,21} &  \frac{1}{\sqrt{n} v_n} {\bf Z}_{n,22} -  {\bf I}_{p-1} - \frac{1}{\sqrt{n} v_n} \ell_{n1} {\bf I}_{p-1} \end{array} \right).
\end{eqnarray*}
Since $\ell_{n1}$ is $O_{\rm P}(1)$ (Lemma~\ref{Lemeigenvalues}) and so are~${\bf Z}_{n,21}$ and ${\bf Z}_{n,22}$, we obtain that~$(E_{n,11},{\bf E}_{n,12})'=(1,0,\ldots,0)'+o_{\rm P}(1)$ (recall that~$E_{n,11}>0$ almost surely) and that~$\sqrt{n} v_n {\bf E}_{n,12}'=O_{\rm P}(1)$. Since ${\bf E}_{n,22}$ is bounded, it directly follows from \eqref{12eq} that $\sqrt{n} v_n {\bf E}_{n,21}=O_{\rm P}(1)$. In view of~(\ref{14eq}), we then obtain that ${\bf E}_{n,22} {\bf E}_{n,22}' - \mathbf{I}_{p-1}$ is~$o_{\rm P}(1)$. Now, by using~(\ref{13eq}), we have
\begin{eqnarray} 
\label{thething}
{\bf Z}_{n,21}
&\!\!\!=\!\!\!& 
\ell_{n1} E_{n,11} {\bf E}_{n,12}' + 
{\bf E}_{n,22}' {\rm diag}(\ell_{n2},\ldots,\ell_{np}) {\bf E}_{n,21} 
\nonumber 
\\[1mm]
& & 
\hspace{12mm} 
+ \sqrt{n} (1+v_n)  E_{n,11} {\bf E}_{n,12}\pr
+ \sqrt{n} {\bf E}_{n,22}\pr{\bf E}_{n,21}
\nonumber 
\\[1mm] 
&\!\!\!=\!\!\!& 
{\bf E}_{n,22}' \; {\rm diag} (\ell_{n2}-\ell_{n1},\ldots,\ell_{np}-\ell_{n1}) {\bf E}_{n,21}  - \sqrt{n} v_n  {\bf E}_{n,22}\pr {\bf E}_{n,21}
,
\end{eqnarray}
which yields~$\sqrt{n}v_n{\bf E}_{n,22}\pr{\bf E}_{n,21} = - {\bf Z}_{n,21} +o_{\rm P}(1)$. Since Lemma~\ref{LemLL} entails that~${\bf Z}_{n,21}$ is asymptotically {${\cal N}({\bf 0}, (1+v){\bf I}_{p-1})$ in case~(i) and ${\cal N}({\bf 0},{\bf I}_{p-1})$ in case~(ii)}, the asymptotic normality result for~$\sqrt{n}r_n{\bf E}_{n,22}\pr{\bf E}_{n,21}$ follows, which, in view of~(\ref{13eq}), also establishes the one for~$\sqrt{n}r_n{\bf E}_{n,12}'$. 

We turn to the proof of~(iii)--(iv). As above, $\wb_{n1}=\Gamb\pr \hat{\thetab}_{n1}=(E_{n,11}, {\bf E}_{n,12})\pr$ is the unit eigenvector associated with the largest eigenvalue~$\ell_{n1}=\sqrt{n} (\hat{\lambda}_{n1}-(1+v_n))$ of~(\ref{equZlam1}), or equivalently, with the largest eigenvalue~$\tilde{\ell}_{n1}=\ell_{n1}+\sqrt{n}v_n=\sqrt{n} (\hat{\lambda}_{n1}-1)$ of
 \begin{equation}
\label{shdd}
{\bf Z}_n+{\rm diag}(\sqrt{n}v_n,0,\ldots,0)
.
\end{equation}
Similarly,~$\wb_{nj}:=\Gamb\pr \hat{\thetab}_{nj}
={\bf E}_n'\eb_j
$, $j=2,\ldots,p$, where~$\eb_j$ stands for the $j$th vector of the canonical basis of~$\R^p$, are the unit eigenvectors associated with the $p-1$ smallest eigenvalues~$\ell_{n2}=\sqrt{n} (\hat{\lambda}_{n2}-1),\ldots,\ell_{np}=\sqrt{n} (\hat{\lambda}_{np}-1)$ of~(\ref{shdd}). Consequently, the joint distribution of~$\wb_{nj}$, $j=1,\ldots,p$ --- that is, the joint distribution of the columns of~${\bf E}_n'$ --- converges weakly to the joint distribution of the unit eigenvectors (associated with eigenvalues in decreasing order, and with the signs fixed as in the statement of the theorem) of
$$
{\bf Z}+\lim_{n\to\infty} {\rm diag}(\sqrt{n}v_n,0,\ldots,0)
$$
(recall that, in cases~(iii)--(iv), $\Zb_n$ converges weakly to the random matrix~$\Zb$ defined in Lemma~\ref{Lemeigenvalues}). This establishes the result. 
 \cqfd
\vspace{4mm}


\section{Technical proofs for Section~3}
\label{SecSupNull}


We start with the proof of Theorem~\ref{CorolAnd}. 
\vspace{3mm}

\noindent {\bf Proof of Theorem~\ref{TheorHPVNull}}. In this proof, all stochastic convergences are as~$n\to\infty$ under~${\rm P}_{\thetab_1^0,r_n,v}$. Since $\tilde{\thetab}_j\pr \Sigb_n \thetab_{1}^0=0$ for~$j=2,\ldots,p$, the test statistic in~(\ref{HPV2010}) rewrites
$$
Q_{\rm HPV}
=
 \frac{n}{\hat{\lambda}_{n1}}  \sum_{j=2}^p \hat{\lambda}_{nj}^{-1} 
 \big(\tilde{\thetab}_{nj}\pr {\bf S}_n \thetab_{1}^0 \big)^2
 =
 \,
 {\bf N}_n\pr {\bf N}_n
 ,
$$
 with 
$$
 {\bf N}_n
 := 
 \frac{1}{\sqrt{\hat{\lambda}_{n1}}}
 \,
 \Big(
\hat{\lambda}_{n2}^{-1/2} 
\,
  \tilde{\thetab}_{n2}\pr  \sqrt{n}({\bf S}_n- \Sigb_n) \thetab_{1}^0
,
\ldots
,
\hat{\lambda}_{np}^{-1/2} 
\,
  \tilde{\thetab}_{np}\pr  \sqrt{n}({\bf S}_n- \Sigb_n) \thetab_{1}^0
 \Big)'
 .
$$
Since the $\tilde{\thetab}_{nj}$'s are unit vectors
\vspace{-.5mm}
  and since $ \sqrt{n}({\bf S}_n- \Sigb_n)$ is $O_{\rm P}(1)$ as soon as~$r_n=O(1)$ (Lemma~\ref{LemLL}), we have that $\tilde{\thetab}_{nj}\pr  \sqrt{n}({\bf S}_n- \Sigb_n) \thetab_{1}^0
$ is $O_{\rm P}(1)$ for any~$j=2,\ldots,p$. Therefore, Lemma \ref{Lemeigenvalues}  entails that, still with~$\delta$ being equal to one if~$r_n\equiv 1$ and to zero if~$r_n=o(1)$, 
$$ 
{\bf N}_n 
=
\frac{1}{\sqrt{1+ \delta v}} 
\Big(\tilde{\thetab}_{n2}\pr  \sqrt{n}({\bf S}_n- \Sigb_n) \thetab_1^0, \ldots,  \tilde{\thetab}_{np}\pr  \sqrt{n}({\bf S}_n- \Sigb_n) \thetab_1^0\Big)+o_{\rm P}(1)
.
$$
Using the fact that~$\sum_{j=2}^p \tilde{\thetab}_{nj} \tilde{\thetab}_{nj}\pr={\bf I}_{p}- \thetab_{1}^0\thetab_{1}^{0\prime}$ is idempotent, it follows that 
 \begin{eqnarray} 
Q_{\rm HPV}
 &\!\!=\!\!& 
 \frac{n}{1+\delta v} \sum_{j=2}^p \,  \thetab_{1}^{0\prime} ({\bf S}_n- \Sigb_n) \tilde{\thetab}_{nj}\tilde{\thetab}_{nj}\pr ({\bf S}_n- \Sigb_n) \thetab_{1}^0 +o_{\rm P}(1) 
 \nonumber 
 \\[1mm]
 &\!\!=\!\!& 
 \frac{n}{1+\delta v} \, \thetab_{1}^{0\prime} ({\bf S}_n- \Sigb_n) ({\bf I}_{p}- \thetab_{1}^{0}\thetab_{1}^{0\prime}) ({\bf S}_n- \Sigb_n)\thetab_{1}^0 +o_{\rm P}(1)
\label{asympeqproofpre}
\\[3mm]
 &\!\!=\!\!& 
 {\bf V}_n\pr {\bf V}_n +o_{\rm P}(1)
 ,
 \nonumber 
 \end{eqnarray}
where we let 
 $
 {\bf V}_n:=(1+\delta v)^{-1/2} \, ({\bf I}_{p}- \thetab_{1}^{0}\thetab_{1}^{0\prime}) \sqrt{n}({\bf S}_n- \Sigb_n)\thetab_{1}^{0}
.
 $
Since Lemma~\ref{LemLL} entails that~$\Vb_n$ is asymptotically normal with mean zero and covariance matrix ${\bf I}_{p}- \thetab_{1}^{0}\thetab_{1}^{0\prime}$ and since~${\bf I}_{p}- \thetab_{1}^{0}\thetab_{1}^{0\prime}$ is idempotent with trace~$p-1$, the result then follows from, e.g., Theorem~9.2.1 in \cite{Ra71}. 
\cqfd
 \vspace{4mm}


\noindent {\bf Proof of Theorem~\ref{TheorAndNull}}. In this proof, all stochastic convergences are as~$n\to\infty$ under~${\rm P}_{\thetab_1^0,r_n,v}$. Before proceeding, note that~(\ref{Andstat}) allows to write
\begin{eqnarray} 
\lefteqn{
\hspace{-59mm} 
Q_{\rm A}
=
  \sum_{j=2}^p \Bigg(\frac{\sqrt{n} (\hat{\lambda}_{n1}-(1+r_n v))-\sqrt{n}(\hat{\lambda}_{nj}-1)}{\sqrt{\hat{\lambda}_{n1}\hat{\lambda}_{nj}}}  \, \hat{\thetab}_{nj}\pr \thetab_1^0
 +
  \frac{\sqrt{n}r_nv}{\sqrt{\hat{\lambda}_{n1}\hat{\lambda}_{nj}}} \, \hat{\thetab}_{nj}\pr \thetab_1^0 \Bigg)^{\!2}  
}
  \nonumber
   \\[2mm]
& & 
\hspace{-50mm} 
=:
 \sum_{j=2}^p \big( R_{nj}+S_{nj} \big)^2
.
\label{Andstat23}
 \end{eqnarray} 

Let us start with the proof of Parts~(i)--(ii) of the theorem. Parts~(i)--(ii) of Lemmas~\ref{Lemeigenvalues}--\ref{Lemeigenvectors} imply that, for any~$j=2,\ldots,p$, $R_{nj}=o_{\rm P}(1)$ and~$S_{nj}=O_{\rm P}(1)$, so that
$$
Q_{\rm A}
=
  \sum_{j=2}^p \Bigg(  \frac{\sqrt{n}r_nv}{\sqrt{\hat{\lambda}_{n1}\hat{\lambda}_{nj}}} \, \hat{\thetab}_{nj}\pr \thetab_1^0 \Bigg)^2  
 +o_{\rm P}(1)
$$
in regimes~(i)--(ii). Consequently,
$$
Q_{\rm A}
=
  \frac{v^2}{1+\delta v}
  \,
  \big( \sqrt{n}r_n {\bf E}_{n,21}\big)' \big( \sqrt{n}r_n{\bf E}_{n,21} \big)
+o_{\rm P}(1)
,
$$
where we used the fact that~$\sqrt{n}r_n {\bf E}_{n,21}=O_{\rm P}(1)$. Therefore, we may write
 \begin{eqnarray*} 
Q_{\rm A}
&\!\!=\!\!&
  \frac{v^2}{1+\delta v}
  \,
  \big( \sqrt{n}r_n {\bf E}_{n,22}\pr{\bf E}_{n,21}\big)'
  \big( \sqrt{n}r_n {\bf E}_{n,22}\pr{\bf E}_{n,21}\big)
\\[2mm]
& & 
\hspace{3mm} 
+
  \frac{v^2}{1+\delta v}
  \,
  \big( \sqrt{n}r_n {\bf E}_{n,21}\big)'
  (\mathbf{I}_{p-1}-{\bf E}_{n,22}{\bf E}_{n,22}\pr)
  \big( \sqrt{n}r_n {\bf E}_{n,21}\big)
+o_{\rm P}(1),
 \end{eqnarray*} 
and the result follows from Lemma~\ref{Lemeigenvectors}(i)--(ii). 
\vspace{3mm}

Let us turn to the proof of Part~(iii). Using Lemmas~\ref{Lemeigenvalues}--\ref{Lemeigenvectors},
the decomposition~(\ref{Andstat23}) readily yields that, in regime~(iii), $Q_{\rm A}$ converges weakly to
$$
\sum_{j=2}^p 
\,
(\ell_1-\ell_j+v)^2 
\big(
(\wb_1(v),\ldots,\wb_p(v))')_{j1}
\big)^{2} 
=
\sum_{j=2}^p 
\,
(\ell_1-\ell_j+v)^2 
\big(w_{j1}(v)\big)^{2} 
, 
$$
where~$\ell_1$ is the largest eigenvalue of~${\bf Z}-{\rm diag}(0, v,\ldots,v)$, $\ell_2\geq \ldots\geq \ell_p$ are the $p-1$ smallest eigenvalues of~${\bf Z}+{\rm diag}(v, 0, \ldots, 0)$, and~$\wb_j(v)=(w_{j1}(v),\ldots,w_{jp}(v))'$ is the unit eigenvector associated with the $j$th largest eigenvalue of~${\bf Z}+{\rm diag}(v,0,\ldots,0)$ and satisfying~$w_{j1}(v)>0$ almost surely (inspection of the proofs of Lemmas~\ref{Lemeigenvalues}--\ref{Lemeigenvectors} reveals that the dependence between the~$\ell_j$'s and the~$\wb_j(v)$'s is through the fact that these are computed from the same random matrix~$\Zb$). The result then follows from the fact that~$\ell_1(v):=\ell_1+v$, $\ell_2(v):=\ell_2,\ldots,\ell_p(v):=\ell_p$ are the eigenvalues (in decreasing order) of~${\bf Z}+{\rm diag}(v,0,\ldots,0)$. Since the proof of Part~(iv) of the result follows exactly along the same lines by taking $v=0$ in the proof of Part~(iii), this establishes the result.
\cqfd
\vspace{4mm}


\noindent {\bf Proof of Corollary~\ref{CorolAnd}}. 
Consider the bivariate case~$p=2$ and regime~(iv). Applying Lemmas~\ref{Lemeigenvalues}--\ref{Lemeigenvectors} to the decomposition~(\ref{Andstat23}) readily yields that
\begin{equation}
	\label{QAasreg4}
		Q_{\rm A}
\stackrel{\mathcal{D}}{\to} 
\big(\ell_1-\ell_2\big)^2 E_{21}^2
,
\end{equation}
where~${\pmb\ell}=(\ell_1,\ell_2)'$ and~$E_{21}$ are mutually independent and have the distributions described in Lemma~\ref{Lemeigenvalues}(iv) and Lemma~\ref{Lemeigenvectors}(iv), respectively (mutual independence
\vspace{-.7mm}
  follows from \cite{And63}, once it is seen that the asymptotic joint distribution of~${\pmb\ell}_n=(\ell_{n1},\ell_{n2})'$ and~$E_{n,21}=\hat{\thetab}_{n2}\pr \thetab_1^0$ is the same in regime~(iv) as in the spherical case associated with~$v=0$). In particular,~$Z:=E_{21}^2$ is the squared of the lower left entry of a random matrix~${\bf E}$ that has an invariant Haar distribution over the collection of $2\times 2$ orthogonal matrices; see the remark below Lemma~\ref{Lemeigenvectors}. It follows that~$Z$ is ${\rm Beta}(\frac{1}{2},\frac{1}{2})$, hence has density
$
z\mapsto f_{Z}(z)
=
1/(\pi \sqrt{z(1-z)})
.
$
As for~${\pmb\ell}$, it has density 
$$
(\ell_1,\ell_2)' 
\mapsto 
f_{(\ell_1,\ell_2)}(\ell_1,\ell_2)
=
\frac{1}{4\sqrt{2\pi}} \, (\ell_1 - \ell_2) e^{-(\ell_1^2+\ell_2^2)/4}
\mathbb{I}[\ell_1\geq \ell_2] 
,
$$ 
where~$\mathbb{I}[A]$ is the indicator function of~$A$. Letting~$U_1 = (\ell_1-\ell_2)/\sqrt{2}$ and~$U_2 = (\ell_1+\ell_2)/\sqrt{2}$, we thus have that~$(U_1,U_2)'$ has density
\begin{eqnarray*}
\lefteqn{	
(u_1,u_2)'\mapsto f_{(U_1,U_2)}(u_1,u_2) = f_{(\ell_1,\ell_2)}\big((u_1+u_2)/\sqrt{2},(-u_1+u_2)/\sqrt{2}\big)
}
\\[2mm]
& & \hspace{3mm} 
=
\frac{1}{4\sqrt{2\pi}} \, \sqrt{2} u_1 e^{-(u_1^2+u_2^2)/4} 
\mathbb{I}[u_1\geq 0] 
=
\frac{1}{4\sqrt{\pi}} \,  u_1 e^{-(u_1^2+u_2^2)/4} 
\mathbb{I}[u_1\geq 0] 
.
\end{eqnarray*}
Therefore, $U_1=(\ell_1-\ell_2)/\sqrt{2}$ has density
\begin{eqnarray*}
\lefteqn{	
\hspace{-1mm} 
f_{U_1}(u_1)
=
\mathbb{I}[u_1\geq 0] 
\int_{-\infty}^\infty
\frac{1}{4\sqrt{\pi}} \,  u_1 e^{-(u_1^2+u_2^2)/4} 
\,du_2
}
\\[2mm]
& & \hspace{-5mm} 
=
\frac{1}{2}
u_1 e^{-u_1^2/4} 
\mathbb{I}[u_1\geq 0] 
\int_{-\infty}^\infty
\frac{1}{\sqrt{2\pi}\times\sqrt{2}} \, e^{-u_2^2/4} 
\,du_2
=
\frac{1}{2}
u_1 e^{-u_1^2/4}
\mathbb{I}[u_1\geq 0] 
.
\end{eqnarray*}
We conclude that~$Y=(\ell_1-\ell_2)^2=2U_1^2$ has density
$$
y\mapsto f_Y(y)
=
\frac{1}{2\sqrt{2y}} f_{U_1}(\sqrt{y}/\sqrt{2})
=
\frac{1}{8} 
 e^{-y/8}
\mathbb{I}[y\geq 0] 
.
$$
Finally, we obtain that the asymptotic distribution of~$Q_{\rm A}$ in the regime considered, that coincides with the distribution of~$YZ$ (see~(\ref{QAasreg4})), has density 
\begin{eqnarray}
\lefteqn{
\hspace{-.5mm} 
q\mapsto f_Q(q)
=
\mathbb{I}[q>0] 
\int_0^1 \frac{1}{z} f_Y(q/z) f_Z(z)\, dz
=
\frac{1}{8\pi}
\mathbb{I}[q>0] 
\int_0^1  \frac{e^{-q/(8z)}}{z^{3/2}(1-z)^{1/2}}\, dz
}
\nonumber
\\[2mm]
& &
\hspace{4mm} 
\,
=
\frac{1}{8\pi}
\bigg[
-
\frac{4\sqrt{2\pi}}{\sqrt{q}} 
e^{-q/8}
\Phi\bigg( \sqrt{\frac{q(1-z)}{4z}} \bigg)
\bigg]_0^1
\mathbb{I}[q>0] 
=
\frac{1}{\sqrt{8\pi q}} e^{-q/8} \mathbb{I}[q>0] 
,
\label{integraltocompute}
\end{eqnarray}
which, as was to be proved, is the pdf of the~$4\chi^2_1$ distribution. To justify the integral computation in~(\ref{integraltocompute}), note that, denoting respectively as~$\Phi$ and~$\phi$ the cumulative distribution function and probability density function of the standard normal distribution, 
\begin{eqnarray*}
\lefteqn{
\frac{d}{dz} \Phi\bigg( \sqrt{\frac{q(1-z)}{4z}} \bigg)
=
\phi\bigg( \sqrt{\frac{q(1-z)}{4z}} \bigg)
\frac{1}{2} \sqrt{\frac{4z}{q(1-z)}}
\bigg( \frac{-4qz-4q(1-z)}{16z^2} \bigg)
}
\\[2mm]
& & 
\hspace{-7mm} 
=
-
\frac{1}{4} 
\phi\bigg( \sqrt{\frac{q(1-z)}{4z}} \bigg)
\sqrt{\frac{q}{z^3(1-z)}}
=
-
\frac{1}{4\sqrt{2\pi}} 
\exp\bigg(- \frac{q(1-z)}{8z} \bigg)
\sqrt{\frac{q}{z^3(1-z)}}
,
\end{eqnarray*}
so that
$$
-
\frac{4\sqrt{2\pi}}{\sqrt{q}} 
\exp\bigg(-\frac{q}{8} \bigg)
\frac{d}{dz} \Phi\bigg( \sqrt{\frac{q(1-z)}{4z}} \bigg)
=
\exp\bigg(- \frac{q}{8z} \bigg)
\sqrt{\frac{1}{z^3(1-z)}}
,
$$
that is,
$$
 \int  \frac{e^{-q/(8z)}}{z^{3/2}(1-z)^{1/2}}\, dz
 =
-
\frac{4\sqrt{2\pi}}{\sqrt{q}} 
\exp\bigg(-\frac{q}{8} \bigg)
\Phi\bigg( \sqrt{\frac{q(1-z)}{4z}} \bigg)
.
\vspace{-9.5mm}
$$
\cqfd
\vspace{6mm}


\section{Technical proofs for Section~4}
\label{SecSupLAN}


We start with the proof of Theorem~\ref{TheorLAN}.
\vspace{3mm}

\noindent {\bf Proof of Theorem~\ref{TheorLAN}}. 
Clearly,
$
\Sigb_n 
=
(1+r_n v) \thetab_1 \thetab_1\pr
+
({\bf I}_p-\thetab_1\thetab_1\pr) 
$
has determinant~$1+r_n v$ and inverse matrix
$
\Sigb_n^{-1}
=
(1+r_n v)^{-1} \thetab_1 \thetab_1\pr
+
({\bf I}_p-\thetab_1\thetab_1\pr) 
=
{\bf I}_p- (r_nv/(1+r_n v))\thetab_1\thetab_1\pr
$.
Thus, letting~$\Sb_{n0}:=\frac{1}{n}\sum_{i=1}^n \Xb_{ni}\Xb_{ni}'$, the hypothesis~${\rm P}_{\thetab_1,r_n,v}$ is associated with the density (with respect to the Lebesgue measure~$\mu$ over~$(\R^p)^n$)
\begin{eqnarray*}
\frac{d{\rm P}_{\thetab_1,r_n,v}}{d\mu}
&\!\!\!=\!\!\!&
\frac{(2\pi)^{-np/2}}{(1+r_n v)^{n/2}} 
\exp\Big( \!-\frac{1}{2} \, \sum_{i=1}^n \Xb_{ni}' 
\big[{\bf I}_p-(r_nv/(1+r_n v)) \thetab_1\thetab_1\pr\big]
\Xb_{ni}
\Big)
\\[2mm]
&\!\!\!=\!\!\!&
\frac{(2\pi)^{-np/2}}{(1+r_n v)^{n/2}} 
\exp\Big( \!
-\frac{1}{2} \, \sum_{i=1}^n \|\Xb_{ni}\|^2
+\frac{nr_n v}{2(1+r_n v)} 
\thetab_1' \Sb_{n0} \thetab_1 \Big)
,
\end{eqnarray*}
which yields (see~(\ref{deflikeratio}))
\begin{eqnarray*}
	\Lambda_n
&\!\!\!=\!\!\!&
\frac{nr_n v}{2(1+r_n v)}
\Big( 
(\thetab_1^0+\nu_n\taub_n)' \Sb_{n0} (\thetab_1^0+\nu_n \taub_n) 
-
\thetab_1^{0\prime} \Sb_{n0} \thetab_1^0 
\Big)
\\[2mm]
&\!\!\!=\!\!\!&
 \frac{n r_n v \nu_n}{1+r_n v}\,  \taub_n\pr\Sb_{n0} \thetab_1^0
 + 
 \frac{n r_n v \nu_n^2}{2 (1+r_nv)}\, \taub_n\pr\Sb_{n0} \taub_n
 .
\end{eqnarray*}

Throughout the rest of the proof, we work under~${\rm P}_{\thetab_1^0,r_n,v}$, and, accordingly, 
$
\Sigb_n 
=
{\bf I}_p+ r_n v\, \thetab_1^0\thetab_1^{0\prime}
$. 
Since 
$
n(\Sb_{n0}-\Sb_{n})
=
n\bar{\Xb}_n \bar{\Xb}_n'
$
is~$O_{\rm P}(1)$, we obtain that, in cases~(i)--(iv), 
$$
	\Lambda_n
=
 \frac{n r_n v \nu_n}{1+r_n v}\,  \taub_n\pr\Sb_{n} \thetab_1^0
 + 
 \frac{n r_n v \nu_n^2}{2 (1+r_nv)}\, \taub_n\pr\Sb_{n} \taub_n
+
o_{\rm P}(1)
.
$$
Now, $\taub_n\pr \Sigb_n \thetab_1^0= (1+r_n v) \taub_n\pr \thetab_1^0=- \frac{1}{2}(1+r_n v) \nu_n \|\taub_n\|^2$ (see~\eqref{tau}), so that
\begin{eqnarray*} 
\hspace{-0mm} 
\Lambda_n
&\!\!\!\!=\!\!\!\!&
 \frac{n r_n v \nu_n}{1+r_n v}\, \taub_n\pr ({\bf S}_n- \Sigb_n) \thetab_1^0
 +  \frac{n r_n v \nu_n^2}{2 (1+r_n v)}\, \taub_n\pr{\bf S}_n \taub_n
 - \frac{n r_n v \nu_n^2}{2}  \| \taub_n\|^2 
+
o_{\rm P}(1)
 \nonumber
  \\[1mm] 
&\!\!\!\!=\!\!\!\!&
 \frac{n r_n v \nu_n}{1+r_n v}\, \taub_n\pr ({\bf S}_n- \Sigb_n) \thetab_1^0
 +  \frac{n r_n v \nu_n^2}{2 (1+r_n v)}\, \taub_n\pr({\bf S}_n- (1+r_n v) {\bf I}_p) \taub_n
+
o_{\rm P}(1)
 .
\end{eqnarray*}
The identity~$(1+r_n v) {\bf I}_p= \Sigb_n+ r_n v ({\bf I}_p- \thetab_1^0 \thetab_1^{0\prime})$  then provides
\begin{eqnarray} 
\Lambda_n
 &\!\!\!=\!\!\!&
\frac{n r_n v \nu_n}{1+v_n} \, \taub_n\pr ({\bf S}_n- \Sigb_n) \thetab_1^0
 +  \frac{n r_n v \nu_n^2}{2 (1+r_n v)} \, \taub_n\pr({\bf S}_n- \Sigb_n) \taub_n
 \nonumber 
 \\[1mm]
 & & 
 \hspace{23mm} 
 -\frac{n r_n^2 v^2 \nu_n^2}{2 (1+r_n v)} \, \taub_n\pr ({\bf I}_p- \thetab_1^0 \thetab_1^{0\prime}) \taub_n 
+
o_{\rm P}(1)
\nonumber
 \\[1mm]
 &\!\!\!=\!\!\!&
\frac{n r_n v \nu_n}{1+r_n v} \, \taub_n\pr ({\bf S}_n- \Sigb_n) \thetab_1^0
 +  \frac{n r_n v \nu_n^2}{2 (1+r_n v)} \, \taub_n\pr({\bf S}_n- \Sigb_n) \taub_n
 \nonumber
 \\[1mm]
 & & 
 \hspace{13mm} 
 -\frac{n r_n^2 v^2 \nu_n^2}{2 (1+r_n v)} \, \| \taub_n\|^2+\frac{n r_n^2 v^2 \nu_n^4}{8 (1+r_n v)} \, \| \taub_n\|^4 
+
o_{\rm P}(1)
,
\label{step3}
 \end{eqnarray}
where we used~\eqref{tau} again. We can now consider the cases~(i)--(iv).

Let us start with cases~(i)--(ii). Take then $\nu_n=1/(\sqrt{n} r_n)$ (recall that~\mbox{$r_n\equiv 1$} in case~(i)) and let~$\delta=1$ (resp., $\delta=0$) if case~(i) (resp., case~(ii)) is considered. The facts that, in both cases~(i)--(ii), $\sqrt{n}({\bf S}_n- \Sigb_n)=O_{\rm P}(1)$ (Lemma~\ref{LemLL}),  $\taub_n\pr \thetab_1^0=- \frac{1}{2} \nu_n \|\taub_n\|^2=o(1)$ and~$\sqrt{n}r_n \nu_n^2=o(1)$ yield
\begin{eqnarray*}
\lefteqn{
\Lambda_n 
=
 \frac{v}{1+\delta v}\, \taub_n\pr \sqrt{n}({\bf S}_n- \Sigb_n) \thetab_1^0-\frac{v^2}{2(1+\delta v)} \| \taub_n\|^2+o_{\rm P}(1) 
}
\\[1mm]
& & 
\hspace{-7mm} 
=
\frac{v}{1+\delta v}\, \taub_n\pr ({\bf I}_{p}-\thetab_1^0 \thetab_1^{0\prime})\sqrt{n}({\bf S}_n- \Sigb_n) \thetab_1^0
+
\frac{v}{1+\delta v} (\taub_n\pr \thetab_1^0) \thetab_1^{0\prime}\sqrt{n}({\bf S}_n- \Sigb_n) \thetab_1^0 
\\[1mm]
& &
\hspace{13mm}  
-\frac{v^2}{2(1+\delta v)}  \taub_n\pr ({\bf I}_{p}-\thetab_1^0 \thetab_1^{0\prime}) \taub_n-\frac{v^2}{2(1+\delta v)}  (\taub_n\pr \thetab_1^0)^2 + o_{\rm P}(1) 
\\[1mm]
& & 
\hspace{-7mm} 
=
\frac{v}{1+\delta v}\, \taub_n\pr ({\bf I}_{p}-\thetab_1^0 \thetab_1^{0\prime})\sqrt{n}({\bf S}_n- \Sigb_n) \thetab_1^0 -\frac{v^2}{2(1+\delta v)}  \taub_n\pr ({\bf I}_{p}-\thetab_1^0 \thetab_1^{0\prime}) \taub_n+o_{\rm P}(1)
.
\end{eqnarray*}
Now, Lemma~\ref{LemLL} implies that
$$
({\bf I}_{p}-\thetab_1^0 \thetab_1^{0\prime})\sqrt{n}({\bf S}_n- \Sigb_n) \thetab_1^0
=
(\thetab_1^{0\prime}\otimes ({\bf I}_{p}-\thetab_1^0 \thetab_1^{0\prime}))
{\rm vec}(\sqrt{n}({\bf S}_n- \Sigb_n))
$$
is asymptotically normal with mean zero and covariance matrix~$(1+\delta v) 
({\bf I}_{p}-\thetab_1^0 \thetab_1^{0\prime})$, which establishes the result in cases~(i)--(ii).

Let us turn to cases~(iii)--(iv). The facts that, in both cases, we have $\sqrt{n}({\bf S}_n- \Sigb_n)=O_{\rm P}(1)$ and $\sqrt{n}r_n=O(1)$ implies that, with~$\nu\equiv 1$,
\eqref{step3} becomes 
\begin{eqnarray} 
\lefteqn{
\Lambda_n
= n r_n v \taub_n\pr ({\bf S}_n- \Sigb_n) \thetab_1^0
 +  \frac{n r_n v}{2} \taub_n\pr({\bf S}_n- \Sigb_n) \taub_n
}
\nonumber
\\[1mm]
& & 
\hspace{33mm} 
 -\frac{n r_n^2 v^2}{2} \| \taub_n\|^2+\frac{n r_n^2 v^2}{8} \| \taub_n\|^4 
+o_{\rm P}(1)
.
 \label{LAQproof}
 \end{eqnarray}
In case~(iv), where~$\sqrt{n}r_n=o(1)$, it directly follows that~$\Lambda_n=o_{\rm P}(1)$. In case~(iii), where~$r_n=1/\sqrt{n}$, (\ref{LAQproof}) yields the announced stochastic quadratic expansion of~$\Lambda_n$. Finally, still in case~(iii), Lemma~\ref{LemLL} implies that, if~$(\taub_n)\to\taub$, then
$
\taub_n\pr \sqrt{n}({\bf S}_n- \Sigb_n) \big( \thetab_1^0+\frac{1}{2} \taub_n\big)
$
is asymptotically normal with mean zero and variance~$\|\taub\|^2 - \frac{1}{4}\|\taub\|^4$,
which establishes the result.
\cqfd
\vspace{3mm}


\noindent {\bf Proof of Theorem~\ref{HPVnonull12}}. 
Since optimality of the test~$\phi_{\rm HPV}$ was established above the statement of the theorem, we only have to derive the asymptotic non-null distribution of~$Q_{\rm HPV}$ under~${\rm P}_{\thetab_1^0+ \taub_n/(\sqrt{n} r_n),r_n,v}$. 
\vspace{-.4mm}
 Recall that, under~${\rm P}_{\thetab_1^0,r_n,v}$, 
$
\sqrt{n}({\bf I}_p-\thetab_1^0\thetab_1^{0\prime})({\bf S}_n-\Sigb_n)\thetab_1^0
$
is asymptotically normal with mean zero and covariance  matrix~$(1+\delta v) 
({\bf I}_{p}-\thetab_1^0 \thetab_1^{0\prime})$; see the proof of Theorem~\ref{TheorLAN}(i)--(ii). A routine application of the Le Cam third lemma thus provides that, under~${\rm P}_{\thetab_1^0+ \taub_n/(\sqrt{n} r_n),r_n,v}$, with~$(\taub_n)\to\taub$, $
\sqrt{n}({\bf I}_p-\thetab_1^0\thetab_1^{0\prime})({\bf S}_n-\Sigb_n)\thetab_1^0
$
is asymptotically normal with mean~$
v \, ({\bf I}_p-\thetab_1^0\thetab_1^{0\prime}) \taub 
$
 and covariance matrix~$(1+\delta v) 
({\bf I}_{p}-\thetab_1^0 \thetab_1^{0\prime})$.  Since contiguity implies that the asymptotic equivalence (see~(\ref{asympeqproofpre}))
$$
Q_{\rm HPV}
=
 \frac{n}{1+\delta v} \, \thetab_{1}^{0\prime} ({\bf S}_n- \Sigb_n) ({\bf I}_{p}- \thetab_{1}^{0}\thetab_{1}^{0\prime}) ({\bf S}_n- \Sigb_n)\thetab_{1}^0 +o_{\rm P}(1)
$$
also holds under~${\rm P}_{\thetab_1^0+ \taub_n/(\sqrt{n} r_n),r_n,v}$, we obtain that, under the same sequence of hypotheses, $Q_{\rm HPV}$ is asymptotically non-central chi-square with~$p-1$ degrees of freedom and with non-centrality parameter
$
(v \, ({\bf I}_p-\thetab_1^0\thetab_1^{0\prime}) \taub)'
((1+\delta v)({\bf I}_p-\thetab_1^0\thetab_1^{0\prime}))^-
(v \, ({\bf I}_p-\thetab_1^0\thetab_1^{0\prime}) \taub)
,
$
which establishes the result (note that~(\ref{tau}) here implies that~$\thetab_1^{0\prime}\taub=0$). 
\cqfd
 \vspace{4mm}


\noindent {\bf Proof of Theorem~\ref{HPVnonull34}}. 
Since the claims in regime~(iv) are trivial, we only show the result in regime~(iii). The sequences of hypotheses~${\rm P}_{\thetab_1^0,1/\sqrt{n},v}$ and~${\rm P}_{\thetab_1^0+ \taub_n,1/\sqrt{n},v}$ are mutually contiguous (which follows, as already mentioned, from the Le Cam first lemma). The Le Cam third lemma then yields that, under~${\rm P}_{\thetab_1^0+ \taub_n,1/\sqrt{n},v}$, with~$(\taub_n)\to\taub$, $
\sqrt{n}({\bf I}_p-\thetab_1^0\thetab_1^{0\prime})({\bf S}_n-\Sigb_n)\thetab_1^0
$
is asymptotically normal with mean
$$
v (1+\thetab_1^{0\prime}\taub)
({\bf I}_{p}-\thetab_1^0 \thetab_1^{0\prime}) \taub
=
v \Big(1-\frac{1}{2} \| \taub\|^2\Big) \Big(\taub+\frac{1}{2} \| \taub\|^2 \thetab_1^0\Big)
$$
and covariance matrix~${\bf I}_{p}-\thetab_1^0 \thetab_1^{0\prime}$, where we used~\eqref{tau}.  As in the proof of the previous theorem, contiguity yields that the asymptotic equivalence (which, in~(\ref{asympeqproofpre}), also holds in regime~(iii))
$$
Q_{\rm HPV}
=
 \frac{n}{1+\delta v} \, \thetab_{1}^{0\prime} ({\bf S}_n- \Sigb_n) ({\bf I}_{p}- \thetab_{1}^{0}\thetab_{1}^{0\prime}) ({\bf S}_n- \Sigb_n)\thetab_{1}^0 +o_{\rm P}(1)
$$
extends to~${\rm P}_{\thetab_1^0+ \taub_n,1/\sqrt{n},v}$, with~$(\taub_n)\to\taub$, and we may therefore conclude that, under this sequence of hypotheses, $Q_{\rm HPV}$ is asymptotically non-central chi-square with~$p-1$ degrees of freedom and with non-centrality parameter
$$
\Big[
v \Big(1-\frac{1}{2} \| \taub\|^2\Big) \Big(\taub+\frac{1}{2} \| \taub\|^2 \thetab_1^0\Big)
\Big]'
({\bf I}_p-\thetab_1^0\thetab_1^{0\prime})^-
\Big[
v \Big(1-\frac{1}{2} \| \taub\|^2\Big) \Big(\taub+\frac{1}{2} \| \taub\|^2 \thetab_1^0\Big)
\Big]
.
$$
This establishes the result since a direct computation shows that this non-centrality parameter coincides with the one in~(\ref{ncpHPV}). 
\cqfd
  \vspace{4mm}


\section{Technical proofs for Section~5} 
\label{SecSupellipt}


The results stated in Section~\ref{sec:ellipt} require the following elliptical versions of Lemmas~\ref{Lemeigenvalues}--\ref{Lemeigenvectors}. 
\vspace{3mm}

\begin{Lem} 
\label{Lemeigenvaluesellipt}
Fix a unit $p$-vector~$\thetab_1$, $v> 0$, a bounded  positive real sequence~$(r_n)$, and~$f\in\mathcal{F}$. 
Let~$\Zb_f(v)$ be a $p\times p$ random matrix such that
$$
{\rm vec}({\bf Z}_f(v))\sim {\cal N}\big({\bf 0}, (1+\kappa_p(f)) ({\bf I}_{p^2}+ {\bf K}_p)(\Lamb(v))^{\otimes 2}+\kappa_p(f) ({\rm vec}\,\Lamb(v)) ({\rm vec}\,\Lamb(v))\pr\big)
,
$$ 
with~$\Lamb(v):={\rm diag}(1+v,1,\ldots,1)$, and let~$\Zb_{f,22}(v)$ be the matrix obtained from~$\Zb_f(v)$ by deleting its first row and first column. Write~$\Zb_f:=\Zb_f(0)$ and~$\Zb_{f,22}:=\Zb_{f,22}(0)$. Then, under~${\rm P}_{\thetab_1,r_n,v,f}$, 
$$
{\pmb \ell}_n:=\big(\sqrt{n}(\hat{\lambda}_{n1}- (1+r_nv)), \sqrt{n}(\hat{\lambda}_{n2}-1), \ldots,  \sqrt{n}(\hat{\lambda}_{np}-1) \big)\pr
\stackrel{\mathcal{D}}{\to}
{\pmb \ell}=(\ell_1, \ldots,\ell_p)'
,
$$
where~${\pmb \ell}$ is as follows:
\begin{itemize}
\item[(i)] if $r_n\equiv 1$, then~$\ell_1=(\Zb_f(v))_{11}$ $($hence is normal with mean zero and variance~$(2+3\kappa_p(f))(1+v)^2)$ and $\ell_2\geq  \ldots\geq \ell_p$ are the eigenvalues of~${\bf Z}_{f,22}(v)$;
\vspace{1mm} 
\item[(ii)] if $r_n$ is $o(1)$ with $\sqrt{n} r_n  \to \infty$, then~$\ell_1=(\Zb_f)_{11}$ $($hence is normal with mean zero and variance~$2+3\kappa_p(f))$ and $\ell_2\geq  \ldots\geq \ell_p$ are the eigenvalues of~${\bf Z}_{f,22}$;
\vspace{1mm} 
\item[(iii)] if $r_n= 1/\sqrt{n}$, then~$\ell_1$ is the largest eigenvalue of~${\bf Z}_f-{\rm diag}(0, v,\ldots,v)$ and $\ell_2\geq \ldots\geq \ell_p$ are the $p-1$ smallest eigenvalues of~${\bf Z}_f+{\rm diag}(v, 0, 
\linebreak
\ldots, 0)$;
\vspace{1mm} 
\item[(iv)] If $r_n=o(1/\sqrt{n})$, then~${\pmb \ell}$ is the vector of eigenvalues of~${\bf Z}_f$ $($in decreasing order$)$, hence has density 
\begin{eqnarray}
\label{explidenslellipt}	
	\lefteqn{
\hspace{16mm} 
(\ell_1, \ldots, \ell_p)' 
\mapsto 
b_{p,f}
\exp
\bigg(
-\frac{1}{4(1+\kappa_p(f))}\, 
\bigg\{
\bigg( \sum_{j=1}^p \ell_j^2 \bigg)
}
\\[2mm]
& & 
\hspace{-10mm} 
-
\frac{\kappa_p(f)}{(p+2)\kappa_p(f)+2} 
\bigg(\sum_{j=1}^p \ell_j\bigg)^2
\bigg\}
\bigg)
\bigg( \prod_{1\leq k<j \leq p} (\ell_k- \ell_j) \bigg)
\,
\mathbb{I}[\ell_1\geq\ldots\geq \ell_p]
,
\nonumber
\end{eqnarray}
where~$b_{p,f}$ is a normalizing constant.
\end{itemize}
\end{Lem}
\vspace{1mm}

Recall that in the corresponding Gaussian result, namely Lemma~\ref{Lemeigenvalues}, $\ell_1$ and~$(\ell_2, \ldots,\ell_p)'$ are mutually independent in cases~(i)--(ii). In Lemma~\ref{Lemeigenvaluesellipt}(i)
\linebreak
--(ii), this independence holds if and only if~$\kappa_p(f)=0$, that is, if and only if the underlying elliptical kurtosis coincides with the multinormal one.  
\vspace{3mm}

{\bf Proof of Lemma~\ref{Lemeigenvaluesellipt}}. With the exception of~(\ref{explidenslellipt}), the proof of this lemma simply follows by replacing the random matrix~$\Zb(v)$ (resp.,~$\Zb$) by~$\Zb_f(v)$  (resp.,~$\Zb_f$) in the proof of Lemma~\ref{Lemeigenvalues}. In particular, note that $(\Zb_f(v))_{11}\sim \mathcal{N}(0,(2+3\kappa_p(f))(1+v)^2)$ in case~(i), whereas $(\Zb_f)_{11}\sim \mathcal{N}(0,2+3\kappa_p(f))$ in case~(ii). Deriving the explicit density in  (\ref{explidenslellipt}) can also be bone by using the same approach as in Lemma~\ref{Lemeigenvalues} but more changes are needed to obtain the result, and we will therefore be more explicit for this part of the proof. Let~$\Wb_f:={\rm vech}\,\Zb_f=\Db_p^-({\rm vec}\,\Zb_f)$, where~$\Db_p^-=(\Db_p'\Db_p)^{-1}\Db_p'$ is the Moore-Penrose inverse of the duplication matrix~$\Db_p$. By definition of~$\Zb_f$, the random vector~$\Wb_f$ has density~$
{\bf w}\mapsto g({\bf w})$, with
$$
g({\bf w})
=
a_{p,f}
\exp\Big(-\frac{1}{2}\, {\bf w}' 
\big(
\Db_p^-
\{
(1+\kappa_p(f)) ({\bf I}_{p^2}+ {\bf K}_p)
+
\kappa_p(f) \Jb_p
\}
(\Db_p^-)'
\big)^{-1} 
{\bf w} \Big)
,
$$
where~$a_{p,f}$ is a normalizing constant and where we let~${\bf J}_p=({\rm vec}\,{\bf I}_p) ({\rm vec}\,{\bf I}_p)'$. By using the identities~${\bf K}_p{\bf D}_p={\bf D}_p$ and $\Db_p^-({\rm vec}\,\mathbf{I}_p)=\Db_p'({\rm vec}\,\mathbf{I}_p)={\rm vech}\,\mathbf{I}_p$, it is easy to check that 
\begin{eqnarray*}
\lefteqn{
\big
(\Db_p^-
\{
(1+\kappa_p(f)) ({\bf I}_{p^2}+ {\bf K}_p)
+
\kappa_p(f) \Jb_p
\}
(\Db_p^-)'
\big)^{-1}
}
\\[2mm]
& & 
\hspace{3mm} 
=
\frac{1}{2(1+\kappa_p(f))}\, 
\Db_p'
\bigg\{
\frac{1}{2} ({\bf I}_{p^2}+ {\bf K}_p)
-
\frac{\kappa_p(f)}{2(1+\kappa_p(f))+p \kappa_p(f)} 
\Jb_p
\bigg\}
\Db_p
.
\end{eqnarray*}
 Using the identities~$({\rm vec}\,\Ab)'({\rm vec}\,\Bb)={\rm tr}[\Ab'\Bb]$ and~${\bf K}_p({\rm vec}\,\Ab)={\rm vec}(\Ab')$, the resulting density for~$\Zb_f$ is therefore
\begin{eqnarray*}
\lefteqn{
\zb \mapsto f(\zb)=g({\rm vech}\,\zb)
}
\\[2mm]
& & 
\hspace{3mm} 
=
a_{p,f}
\exp\Big(-\frac{1}{4(1+\kappa_p(f))}\, 
({\rm vec}\,\zb)'
\Big\{
{\bf I}_{p^2}
-
\frac{\kappa_p(f)}{(p+2)\kappa_p(f)+2} 
\Jb_p
\Big\}
{\rm vec}\,\zb \Big)
\\[2mm]
& & 
\hspace{3mm} 
=
a_{p,f}
\exp\Big(
-\frac{1}{4(1+\kappa_p(f))}\, 
\Big\{
({\rm tr}\,\zb^2)
-
\frac{\kappa_p(f)}{(p+2)\kappa_p(f)+2} 
({\rm tr}\,\zb)^2
\Big\}
\Big)
.
\end{eqnarray*}
As in the proof of Lemma~\ref{Lemeigenvalues}, the density of~${\pmb \ell}$ in~(\ref{explidenslellipt}) then follows from Lemma~2.3 of \cite{Tyl83b}. 
\cqfd
\vspace{4mm}

\begin{Lem} 
\label{Lemeigenvectorsellipt}
Fix a unit $p$-vector~$\thetab_1$, $v> 0$, a bounded  positive real sequence~$(r_n)$, and~$f\in\mathcal{F}$. 
Let~$\Zb_f$ be a $p\times p$ random matrix such that
$$
{\rm vec}({\bf Z}_f)\sim {\cal N}\big({\bf 0}, (1+\kappa_p(f)) ({\bf I}_{p^2}+ {\bf K}_p)+\kappa_p(f) ({\rm vec}\,{\bf I}_p) ({\rm vec}\,{\bf I}_p)\pr\big)
.
$$ 
Let ${\bf E}_f(v):=(\wb_1(v),\ldots,\wb_p(v))'$, where~$\wb_j(v)=(w_{j1}(v),\ldots,w_{jp}(v))'$ is the unit eigenvector associated with the $j$th largest eigenvalue of~${\bf Z}_f+{\rm diag}(v,0,
\linebreak
\ldots,0)$ and such that~$w_{j1}(v)>0$ almost surely. Extending the definitions to the case~$v=0$, write~${\bf E}_f:={\bf E}_f(0)$. Then, we have the following 
under~${\rm P}_{\thetab_1,r_n,v,f}$:
\begin{itemize} 
\item[(i)]  
if $r_n\equiv 1$, then~$E_{n,11}=1+o_{\rm P}(1)$, ${\bf E}_{n,22}{\bf E}_{n,22}\pr=\mathbf{I}_{p-1}+o_{\rm P}(1)$, $\sqrt{n} {\bf E}_{n,21}=O_{\rm P}(1)$, and both~$\sqrt{n}{\bf E}_{n,22}\pr{\bf E}_{n,21}$ and~$\sqrt{n}{\bf E}_{n,12}'$ are asymptotically normal with mean zero and covariance matrix~$v^{-2}(1+v)(1+ \kappa_p(f)) {\bf I}_{p-1};$
\vspace{1mm} 
\item[(ii)]   
if~$r_n$ is $o(1)$ with $\sqrt{n} r_n \to \infty$, then~$E_{n,11}=1+o_{\rm P}(1)$, ${\bf E}_{n,22}{\bf E}_{n,22}\pr=\mathbf{I}_{p-1}+o_{\rm P}(1)$, $\sqrt{n} r_n {\bf E}_{n,21}=O_{\rm P}(1)$, and both~$\sqrt{n}r_n{\bf E}_{n,22}\pr{\bf E}_{n,21}$ and \linebreak
$\sqrt{n}r_n{\bf E}_{n,12}'$ are asymptotically normal with mean zero and covariance matrix~$v^{-2} (1+ \kappa_p(f)){\bf I}_{p-1};$
\vspace{1mm} 
\item[(iii)]  
if~$r_n=1/\sqrt{n}$, then~${\bf E}_n$ converges weakly to~${\bf E}_f(v);$
\vspace{1mm} 
\item[(iv)] if~$r_n=o(1/\sqrt{n})$, then ${\bf E}_n$ converges weakly to~${\bf E}_f$, or equivalently, to the random matrix~${\bf E}$ from Lemma~\ref{Lemeigenvectors}.   
\end{itemize}
\end{Lem}

\noindent {\bf Proof of Lemma~\ref{Lemeigenvectorsellipt}}.
The proof of~(i)--(ii) follows by proceeding exactly as in the proof of Lemma~\ref{Lemeigenvectors}(i)--(ii), once it is seen that  Lemma~\ref{TheorHPVNullellipt} implies that, under~${\rm P}_{\thetab_1,r_n,v,f}$, ${\bf Z}_{n,21}$ is asymptotically ${\cal N}({\bf 0}, (1+v)(1+ \kappa_p(f)){\bf I}_{p-1})$ in case~(i) and ${\cal N}({\bf 0},(1+ \kappa_p(f)) {\bf I}_{p-1})$ in case~(ii); here, ${\bf Z}_{n,21}$ refers to the bottom left block (see~(\ref{partitionZ})) of the random matrix~${\bf Z}_n$  in~(\ref{defZn}). 
As for the proof of~(iii)--(iv), it also follows as in the corresponding parts of the proof of Lemma~\ref{Lemeigenvectors} as~${\bf Z}_{n}$, under~${\rm P}_{\thetab_1,r_n,v,f}$, converges weakly to the random matrix~${\bf Z}_{f}$ defined in the statement of the present result (note that the equality in distribution of~${\bf E}_f$ and~${\bf E}$ results from the fact that their respective antecedents~${\bf Z}_f$ and~${\bf Z}$ both have a spherically invariant distribution; see Section~2 in \cite{Tyl83b}). 
\cqfd
\vspace{3mm}

We end this appendix by proving Theorem~\ref{TheorHPVNullellipt}.
\vspace{2mm}

\noindent {\bf Proof of Theorem~\ref{TheorHPVNullellipt}}. Unless mentioned otherwise, all stochastic convergences 
\vspace{-1mm}
 in this proof are as $\ny$ under ${\rm P}_{\thetab_1^0,r_n,v, f}$.  
Clearly, Lemma~\ref{Lemeigenvaluesellipt} entails that $\hat{\lambda}_{n1}=1+ \delta v+o_{\rm P}(1)$ and $\hat{\lambda}_{np}=1+o_{\rm P}(1)$ for $j=2, \ldots, p$ (still with the quantity~$\delta$ defined at the beginning of the appendix). Therefore, the same arguments as in the proof of Theorem~\ref{TheorHPVNull} imply that
\begin{equation}
	\label{ghaj}
\frac{Q_{\rm HPV}\n}{1+ \kappa_p(f)}
 =
 \,
 {\bf V}_{n,f}\pr  {\bf V}_{n,f} + o_{\rm P}(1)
 ,
\end{equation}
 with
$$
 {\bf V}_{n,f}
 := 
 \frac{1}{\sqrt{(1+ \kappa_p(f))(1+ \delta v)}} 
 \,
 ({\bf I}_{p}- \thetab_{1}^{0}\thetab_{1}^{0\prime}) \sqrt{n}({\bf S}_n- \Sigb_n)\thetab_{1}^{0}.
 $$
By using the fact that~$((\thetab_{1}^{0})' \otimes ({\bf I}_{p}- \thetab_{1}^{0}\thetab_{1}^{0\prime})) ({\rm vec}\,{\bf I}_p)={\bf 0}$, it follows from Lemma~\ref{LemLLellipt} that~${\bf V}_{n,f}$ is asymptotically normal with mean zero and covariance matrix ${\bf I}_{p}- \thetab_{1}^{0}\thetab_{1}^{0\prime}$. Since~${\bf I}_{p}- \thetab_{1}^{0}\thetab_{1}^{0\prime}$ is idempotent with trace~$p-1$,~(\ref{ghaj}) and Theorem~9.2.1 in \cite{Ra71} then entail that
\begin{equation}
	\label{tlzodg}
\frac{Q_{\rm HPV}\n}{1+ \kappa_p(f)}
\stackrel{\mathcal{D}}{\to}
\chi^2_{p-1}
.
\end{equation}
In particular, $Q\n_{\rm HPV}$ is $O_{\rm P}(1)$. Therefore, if~$\hat{\kappa}_p\n$ is a weakly consistent estimator of~$\kappa_p(f)$, then
$$
Q^{(n)\dagger}_{\rm HPV}
=
\frac{Q_{\rm HPV}\n}{1+ \hat{\kappa}_p\n}
=
\frac{Q_{\rm HPV}\n}{1+ \kappa_p(f)}
+o_{\rm P}(1)
,
$$
which, in view of~(\ref{tlzodg}),
\vspace{-1mm}
  would prove that $Q^{(n)\dagger}_{\rm HPV}$ is asymptotically~$\chi^2_{p-1}$ and would also show that $
Q_{\rm HPV}^{(n)\dagger}=Q_{\rm HPV}\n+o_{\rm P}(1)
$ in the multinormal case (since~$\kappa_p(\phi)=0$). 

We thus conclude the proof by showing that~$\hat{\kappa}_p\n$ converges to~$\kappa_p(f)$ in probability. With~$\Sigb_n=\sigma_n({\bf I}_p+ r_n v\, \thetab_1\thetab_1\pr)$, the sample average and empirical covariance matrix of the random vectors~$\Yb_{ni}:=\Sigb_n^{-1/2}{\Xb}_{ni}$, $i=1,\ldots,n$, are $\bar{\Yb}_n:=\Sigb_n^{-1/2}\bar{\Xb}_n$ and ${\Sb}_{n,\Yb}:=\Sigb_n^{-1/2}\Sb_n\Sigb_n^{-1/2}$, respectively (here, $\Sigb_n^{-1/2}$ stands for the inverse of the symmetric positive-definite square root of~$\Sigb_n$). Hence, the estimator in~(\ref{hatkur}) rewrites 
$$
\hat{\kappa}_p\n
=
\frac{1}{np(p+2)} 
\sum_{i=1}^{n} 
\,
( ({\Yb}_{ni}-\bar{\Yb}_n)'\Sb_{n,\Yb}^{-1}({\Yb}_{ni}-\bar{\Yb}_n) )^2
-1
.
$$
Now, the triangular array~$\Yb_{ni}$, $i=1,\ldots,n$,  $n=1,2,\ldots$ contains random vectors that are independent and identically distributed, not only in columns but also in rows (the distribution of~$\Yb_{ni}$ 
\vspace{-.5mm}
 does not depend on~$n$). The weak consistency of~$\hat{\kappa}_p\n$ for~$\kappa_p(f)$ in the present triangular array setup thus directly follows from the corresponding result in the standard non-triangular case; see, e.g., page~103 of \cite{And2003}.  
\cqfd


\bibliographystyle{imsart-nameyear.bst} 
\bibliography{Paper.bib}           
\vspace{3mm} 


\end{document}